\documentclass[11pt]{article}
\usepackage{amsmath}
\usepackage{latexsym} 
\usepackage{enumitem}
\usepackage{multirow}
\input{amssym.def} 
\input{amssym}
\usepackage{xcolor}
\usepackage{hyperref}
\date{} 
\newfont{\fra}{eufm10 scaled 1095} 
\newfont{\Bb}{msbm10 scaled 1095} 
\newfont{\Bbs}{msbm10 scaled 1395}
\newfont{\Bbg}{msbm10 scaled 1280} 
\newcommand\CC{{\mbox{\Bb C}}}
\newcommand\RR{{\mbox{\Bb R}}} 
\newcommand\RRR{{\mbox{\Bbs R}}} 
\newcommand\NN{{\mbox{\Bb N}}} 
\newcommand\ZZ{{\mbox{\Bb Z}}}

\newcommand\Z{{\Bbb Z}}
\newcommand\R{{\Bbb R}}
\newcommand\HH{{\mbox{\Bb H}}}

\newcommand\fg{{\frak{g}}} 
\newcommand\fh{{\frak h}}

\newcommand\fl{{\frak l}} 
\newcommand\fm{{\frak m}} 
\newcommand\fn{{\frak n}} 
\newcommand\fa{{\frak a}} 
 
\newcommand\fd{{\frak d}}

\newcommand\fk{{\frak k}} 
\newcommand\fp{{\frak p}}
\newcommand\fr{{\frak r}} 
\newcommand\fs{{\frak s}} 
 
\newcommand\fu{{\frak u}} 
\newcommand\fz{{\frak z}} 

\newcommand\bH{{\mathsf{H}}}
\newcommand\bN{{\mathsf{N}}}
\newcommand\bP{{\mathsf{P}}}
\newcommand\ph{\varphi} 
\newcommand\eps{\varepsilon}

\newcommand{\fsl}{\mathop{{\frak s \frak l}}} 
\newcommand{\fsu}{\mathop{{\frak s \frak u}}} 
 
\newcommand{\fso}{\mathop{{\frak s \frak o}}}
\newcommand{\fsp}{\mathop{{\frak s \frak p}}}
\newcommand{\fiso}{\mathop{{\frak i \frak s \frak o}}}
\newcommand{\End}{\mathop{{\rm End}}} 
\newcommand{\Aff}{\mathop{{\rm Aff}}} 
\newcommand{\Aut}{\mathop{{\rm Aut}}} 
\newcommand{\Iso}{\mathop{{\rm Iso}}}
\newcommand{\GL}{\mathop{{\rm GL}}} 
\newcommand{\SL}{\mathop{{\rm SL}}} 
\newcommand{\PSL}{\mathop{{\rm PSL}}}
\newcommand{\SO}{\mathop{{\rm SO}}} 
\newcommand{\SU}{\mathop{{\rm SU}}}
 
\newcommand{\grO}{{\rm O}}

\newcommand{\id}{\mathop{{\rm id}}} 
\newcommand{\ad}{\mathop{{\rm ad}}} 
\newcommand{\tr}{\mathop{{\rm tr}}} 
\newcommand{\Ad}{\mathop{{\rm Ad}}}

\newcommand{\Ker}{\mathop{{\rm ker}}}

\newcommand{\diag}{\mathop{{\rm diag}}} 
 
\newcommand{\Span}{\mathop{{\rm span}}} 
\newcommand{\pro}{{\rm pr}} 
\newcommand\ip{{\langle\cdot \,,\cdot \rangle}} 
 
\newcommand\lb{{[\cdot\,,\cdot]}} 
 
\newcommand\dd{\fd_{\alpha,\gamma}(\fl,\theta_\fl,\fa)} 
\newcommand\la{\langle}
\newcommand\ra{\rangle}
\newcommand\proof{{\sl Proof. }} 
\newcommand{\qed}{\hspace*{\fill}\hbox{$\Box$}\vspace{2ex}} 
\newcommand{\benur}{\begin{enumerate}[label=(\roman*)]}
\newtheorem{theo}{Theorem}[section] 
\newtheorem{pr}[theo]{Proposition} 
\newtheorem{de}[theo]{Definition} 
\newtheorem{re}[theo]{Remark} 
\newtheorem{co}[theo]{Corollary} 
\newtheorem{lm}[theo]{Lemma} 
\begin{document} 
\title{Pseudo-Riemannian symmetric spaces of signature $(2,2)$}
\author{Ines Kath\thanks{Institut f\"ur Mathematik und Informatik 
der Universit\"at Greifswald,
Walther-Rathenau-Str.\,47, D-17489 Greifswald, Germany, 
ines.kath@uni-greifswald.de, corresponding author} \ and Matti Lyko\thanks{Institut f\"ur Mathematik und Informatik
der Universit\"at Greifswald,
Walther-Rathenau-Str.\,47, D-17489 Greifswald, Germany, 
matti.lyko@uni-greifswald.de}} 
\maketitle 
\begin{abstract}
\noindent We study all four-dimensional simply-connected indecomposable non-semisimple pseu\-do-Riemannian symmetric spaces whose metric has signature $(2,2)$. We present models and compute their isometry groups.  We solve the problem of the existence or non-existence of compact quotients by properly acting discrete subgroups of the isometry group. This continues and completes earlier work by Maeta. 
\end{abstract}
\tableofcontents
\section{Introduction}
While Riemannian symmetric spaces are well known and intensely studied, much less is known about pseudo-Riemannian ones. In the same way as Riemannian simply-connected symmetric spaces, pseudo-Riemannian ones can be described purely algebraically by their associated infinitesimal objects, so-called symmetric triples. For this description one uses the transvection group of the symmetric space, which is a particular subgroup of the isometry group. A symmetric triple consists of the Lie algebra $\hat \fg$ of the transvection group, an $\ad$-invariant inner product $\ip$ on $\hat \fg$ and an involution $\theta$ of $\hat \fg$ satisfying some compatibility conditions. 

Pseudo-Riemannian symmetric spaces that have a semisimple transvection group were classified by Berger \cite{Berger}. These symmetric spaces will be referred to as semisimple. The classification problem for pseudo-Riemannian symmetric spaces that are not semisimple is rather involved and is only solved for small index of the metric. Lorentzian symmetric spaces with a solvable transvection group were classified by Cahen and Wallach \cite{CW}. 
For signature $(2,n)$, first classification results were obtained by Cahen and Parker \cite{CP1, CP2}. However, the results in \cite{CP1} concerning the case of a solvable transvection group were incomplete. This was observed by Neukirchner who revised the classification \cite{N}. In \cite{KOsymm}, a new and more structural approach to the classification problem was developed. It was shown that every indecomposable non-semisimple symmetric triple has the structure of a uniquely determined balanced quadratic extension.  Extensions of this kind can be described by a quadratic cohomology set. For small index of the metric, or under additional assumptions on the metric, this cohomology set can be computed explicitly. This leads to a classification of symmetric triples in these situations. In \cite{KOsymm} this method was applied to the classification in the case of signature $(2,n)$, where it turned out that also Neukirchner’s classification was not quite correct. 

We have already explained that the classification of simply-connected symmetric spaces is accomplished through the classification of symmetric triples. As a consequence, the classification is presented in a purely algebraic form. Here, we aim to derive, in certain cases, a geometric form from this algebraic representation. We want to find suitable models of the symmetric spaces itself and describe their geometry. We will focus on four-dimensional spaces with metrics of signature (2,2). However, the same methods work also for certain families of higher-dimensional spaces. 

As already noted, the classification of indecomposable non-semisimple symmetric triples is based on their description as balanced quadratic extensions. More precisely,  a symmetric triple $(\hat\fg,\theta,\ip)$ of this kind is a balanced quadratic extension of a Lie algebra with involution $(\fl,\theta_\fl)$ by a semi-simple orthogonal $(\fl,\theta_\fl)$-module. For signature $(2,2)$, there are five possible isomorphism classes for $(\fl,\theta_\fl)$. This leads to five families of symmetric triples, where one of these families will be divided into two subfamilies which differ by the $(\fl,\theta_\fl)$-module used in the extension. The corresponding families of simply-connected symmetric spaces will be denoted by $X_1, X_2, N,Y,Z$ and $Z'$, together with a specification of parameters.  The spaces of type $X_1$ and $X_2$ have a solvable non-nilpotent transvection group. Their symmetric triples are balanced quadratic extensions of $(\fl,\theta_\fl)=(\RR,-\id)$. The transvection group of spaces of type $Y$ is also solvable and non-nilpotent. For these spaces $\fl$ is isomorphic to $\fso(2)\ltimes \RR^2$ or to $\fso(1,1)\ltimes \RR^2$.  On the other hand, the spaces of type $N$ have a nilpotent transvection group, with their symmetric triples being balanced quadratic extensions of $(\RR^2,-\id)$. The transvection groups of the spaces of type $Z$ and $Z'$ are not solvable but have a non-trivial Levi factor. Here $\fl$ is isomorphic to  $\fsu(2)$ or $\fsl(2,\RR)$ for the $Z$-spaces and to $\fsl(2,\RR)$ with a different involution for the $Z'$-spaces.

Our goals in this paper are 
\begin{itemize}
    \item[\bf (1)] to find explicit models for the non-semisimple simply-connected symmetric spaces of signature (2,2),
    \item[\bf (2)] to determine their isometry group, 
    \item[\bf (3)] to decide whether these spaces admit compact quotients by discrete subroups of their isometry group.
\end{itemize}

In the following, we will explain these items in more detail, and we will describe our results. 

\paragraph*{(1) Explicit models} From the symmetric triple $(\hat \fg, \theta,\ip)$ we obtain the associated simply-connected symmetric space $X$ as a homogeneous space $\hat G/\hat G_+$. The inner product $\ip$ gives us a $\hat G$-invariant metric $g$ on $X$. So, in principle, the symmetric space $X$ is known. Here we want to give an explicit description of the quotient $\hat G/\hat G_+$ and the metric. 

For spaces in the families $X_1, X_2, Y$ we obtain metrics on $\RR^4$ that  behave similarly to Lorentzian plane waves. They are of the form $2dudv+\ip+H(u,x)du^2$, where $\ip$ is pseudo-Euclidean of signature $(1,1)$. There exists a non-zero parallel vector field $V$, which is unique up to multiplication by a scalar. The transvection group contains the 5-dimensional Heisenberg group, which acts transitively on each leaf of the integrable distribution $V^\perp$.

The spaces of type $N$ have a simple realisation by a metric on $\RR^4$. Since the action of the isometry group looks rather complicated in this model, we will give two alternate descriptions of $N$ as an extrinsic symmetric space in $\RR^5$. However, we do not study extrinsic symmetric spaces systematically here. 

For the $Z$- and $Z'$-spaces we also give a description using extrinsic symmetric spaces. The spaces of type $Z$ can be embedded in $\RR^6$ and interpreted as the tangent bundle of the standard sphere $S^2\subset \RR^3$ or the tangent bundle of the Riemannian hyperbolic plane $H^2\subset \RR^{1,2}$, depending on the parameters of the space $Z$. The spaces in the $Z'$-family are isometric to the universal cover of an extrinsic symmetric space in $\RR^6$, which can be understood as the tangent bundle of the pseudo-Riemannian sphere $S^{1,1}\subset\RR^{1,2}$. 

In order to obtain an embedding as an extrinsic symmetric space in the cases mentioned above, the main tool is the use of pseudo- and para-Hermitian structures. For this purpose we generalise the well-known fact that (Riemannian) Hermitian symmetric spaces admit a realisation as extrinsic symmetric spaces. In Section~\ref{Sext}, we will prove the following: Let $X$ be a symmetric space and $(\hat\fg,\theta,\ip)$ its symmetric triple. If $X$ admits a pseudo-Hermitian structure, then (up to a covering) $X$ can be embedded extrinsically symmetric in $(\hat\fg,\ip)$. If $X$ admits a para-Hermitian structure, then $X$ can be embedded extrinsically symmetric in $(\hat\fg,-\ip)$ (again up to a covering).

\paragraph*{(2) Isometry groups} Isometry groups of Cahen-Wallach spaces, i.e.\ of solvable Lorentzian symmetric spaces are  determined in \cite{KOcq}. The same method can be applied to symmetric spaces of signature $(2,n)$. This is done for some of the spaces of signature $(2,2)$ in \cite{M}. In this paper we will also compute the isometry groups of the remaining spaces of signature $(2,2)$. The isometry group can be essentially larger than the transvection group. For the spaces $X_1(\eps,-\eps,1)$, $\eps=\pm1$, and the two spaces of type  $N$, this indeed happens. Their isometry group has a larger dimension than their transvection group. In all other cases that we consider here, the transvection group has finite index in the isometry group.

\paragraph*{(3) Compact quotients}

By a compact quotient of a homogeneous space $X = G/H$ we mean a compact quotient manifold $\Gamma \backslash X$, where $\Gamma\subset G$ is a discrete subgroup of $G$ acting properly and freely on $X$. Compact quotients are also called compact Clifford-Klein forms. Their existence or non-existence for a given homogeneous space has been intensively studied for many years.
A key aspect is whether the stabiliser $H$ is compact; if it is, then the action of any discrete subgroup of $G$ on $X$ is proper. Hence studying compact quotients of $X$ is essentially equivalent to studying discrete uniform subgroups (also called cocompact lattices) of $G$. However, when dealing with homogeneous spaces with a non-compact stabiliser, the action of a discrete group $\Gamma \subset G$ on $X$ is not necessarily proper. In such cases, the question of the existence of compact quotients becomes a more complicated problem and only partial answers are known.

Let us recall some important results, focusing on the case where $G/H$ is a symmetric space. If the metric is Riemannian, then $H$ is compact. Using this fact Borel showed that every simply-connected Riemannian symmetric space has a compact Clifford-Klein form by proving that a connected semisimple Lie group always contains a discrete uniform subgroup \cite{Bo}. This solves the problem in the Riemannian situation. Now, considering pseudo-Riemannian symmetric spaces, the isotropy group $H$ is in general not compact. Additionally, the isotropy representation may not be completely reducible.  
Even for semisimple symmetric spaces the existence problem of compact quotients is far from being solved. An exception to this are the semisimple groups, which admit discrete uniform subgroups and thus compact quotients. Several partial results were obtained by T.\;Kobayashi. He studied homogeneous spaces $G/H$ of reductive type. He characterised proper actions on spaces of this kind (see \cite{Ko1}) and found sufficient conditions for existence and non-existence of compact quotients (see e.g. \cite{Ko1, Ko2}). These results can be applied to semisimple symmetric spaces, which yields for example tables of symmetric spaces admitting compact quotients. For detailed information and further results, e.g. on tangential symmetric spaces see \cite{KY}. Now let us turn to non-semisimple symmetric spaces. We will say that a pseudo-Riemannian symmetric space is indecomposable if it is not a non-trivial product of pseudo-Riemannian symmetric spaces. Lorentzian symmetric spaces which are indecomposable but not semisimple are called Cahen-Wallach spaces. They have a solvable transvection group.
The question of which Cahen-Wallach spaces have compact quotients is studied in \cite{KOcq}. There are already some results for non-semisimple symmetric spaces with a metric of signature $(2,n)$. Namely, Maeta \cite{M} studied indecomposable spaces with a metric of signature $(2,2)$ whose transvection group is solvable. Among those whose transvection group is not nilpotent, he determined all spaces that have compact quotients by a discrete subgroup of the isometry group. These are exactly two spaces. In the notation of our paper they correspond to $X_1(\eps,-\eps,1)$ for $\eps=\pm1$. One of these spaces, $X_1(1,-1,1)$, is isometric to the hyperbolic oscillator group with its biinvariant metric, which is known to have cocompact lattices. For the two spaces with a nilpotent transvection group, e.g. those of type $N$, he proved that they do not admit compact quotients by discrete subgroups of the transvection group. As noted above, the isometry group of these spaces is larger than the transvection group. It is not solvable and therefore the methods used by Maeta do not apply to this case. So for spaces of type $N$ the problem of the existence of compact quotients by a discrete group of isometries remained open. 

In the present paper we solve the existence problem of compact quotients by isometries for all remaining symmetric spaces of signature $(2,2)$, i.e. for those with a nilpotent or a non-solvable transvection group. We use classical methods like syndetic hulls and the theory of subgroups of $\PSL(2,\RR)$. It turns out that all these spaces do not admit compact quotients:

\begin{theo}
   If $X$ is an indecomposable non-semisimple symmetric space of signature $(2,2)$, then $X$ admits a compact quotient by a discrete subgroup of the isometry group if and only if $X$ is isometric to one of the spaces $X_1(\eps,-\eps,1)$, $\eps=\pm1$.
\end{theo}

Although we defined quotients as quotients by discrete subgroups that act freely, the freeness condition is not important for the spaces considered here. We will apply Selberg's Lemma to show that if $\Gamma$ is a discrete subgroup of the isometry group of a symmetric space $X$ of signature $(2,2)$ acting properly and cocompactly, then $\Gamma$ contains a finite index subgroup that acts freely.

\section{The transvection and the isometry group}
\label{S2}
Let $X$ be a simply-connected semi-Riemannian symmetric space and $\Iso(X)$ its isometry group. We choose a base point $o\in X$. For each $p\in X$ let $\theta_p\in\Iso(X)$ be the point symmetry at $p$. Locally around $p$, this map is the geodesic reflection about $p\in X$. The transvection group $\hat G$ of $X$ is the normal subgroup of $\Iso(X)$ generated by all compositions $\theta_p\circ\theta_q$ for $p,q\in X$.  
It is the smallest subgroup of $\Iso(X)$ that acts transitively on $X$ and is invariant under conjugation by  $\theta_o$\,. In particular, we obtain a diffeomorphism $\hat G /\hat G_+ \to X$, where $\hat G_+\subset\hat G$ denotes the stabiliser of $o$. We denote by $\theta: \Iso(X) \to \Iso(X)$ the conjugation by $\theta_o$. It induces an involutive isomorphism of the Lie algebra $\fiso(X)$ of $\Iso(X)$, and of the Lie algebra $\hat\fg$ of the transvection group $\hat G$. We denote these involutions by $\theta$ as well. For each $\theta$-invariant subspace $\fh$ of $\fiso(X)$ we denote by $\fh_+$ and $\fh_-$ the eigenspaces of $\theta$ with eigenvalues $1$ and $-1$, respectively. In particular, $\fh=\fh_+\oplus\fh_-$. The subspace $\fm:=\fiso(X)_-$ can be naturally identified with the tangent space $T_oX$ at the base point. Indeed, if $\pi$ denotes the projection $\pi:\Iso(X)\to X$, $g\mapsto g(o)$, and $e$ denotes the identity in $\Iso(X)$, then $d\pi_e:\fm\to T_oX$ is an isomorphism. The Lie algebra $\hat \fg$ equals $\hat\fg=\hat\fg_+\oplus\hat\fg_-=[\fm,\fm]\oplus\fm\subset\fiso(X)$, and $\hat\fg_+$ coincides with the Lie algebra of the stabiliser $\hat G_+$. By the natural identification $T_oX\cong \fm$ the metric on $X$ defines a (non-degenerate indefinite) scalar product on $\fm$, which extends to an $\ad(\hat\fg)$-invariant scalar product $\ip$ on $\hat \fg$ for which $\theta$ is an isometry. The triple $(\hat\fg, \theta, \ip)$ is called a \emph{symmetric triple}. It determines $X$ completely \cite[Section 2]{CP2}.   

Let us first recall how we recover the transvection group from the symmetric triple.  Denote by $\tilde G$ the simply-connected Lie group with Lie algebra $\hat \fg$ and by $\tilde G_+$ its connected subgroup with Lie algebra $\hat\fg_+$. Then $X=\tilde G/\tilde G_+$. However, $\tilde G$ is not necessarily equal to the transvection group since, in general, it does not act effectively. Let $N\subset\tilde G$ denote the set of elements that act trivially on $X$. 

\begin{lm}
   The group $N$ is discrete and normal in $\tilde G$. It equals the intersection of $\tilde G_+$ and the center of $\tilde G$. 
\end{lm}
\proof An element $g_0\in \tilde G$ belongs to $N$ if and only if $g_0\in g\tilde G_+ g^{-1}$ for all $g\in\tilde G$. Thus $N=\bigcap_{g\in\tilde G} g\tilde G_+ g^{-1}$. In particular, $N$ is normal in $\tilde G$ and contained in $\tilde G_+$. It is discrete since its Lie algebra $\fn$ satisfies $[\fn,\hat \fg_-]\subset\fn\cap\hat \fg_-=0$, which implies $\fn=0$ since $\hat \fg_+$ acts faithfully on $\hat \fg_-$. As a normal discrete subgroup of $\tilde G$ it is contained in the centre $Z(\tilde G)$. Hence $N\subset \tilde G_+\cap Z(\tilde G)$. The converse inclusion is trivially satisfied. 
\qed

We put $\hat G:=\tilde G/N$ and $\hat G_+ := \tilde G_+/N$. Then, $X=\hat G/\hat G_+$, and $\hat G$ is the transvection group of $X$. 

Let us now describe how we can determine the isometry group of $X$ from the transvection group, which is a closed normal subgroup of the isometry group \cite[Section 1]{CP2}. The method to do this is well-known, but few explicit references can be found in the literature, see e.g. \cite[Proposition 1.5]{N}, \cite[Proof of Proposition 2.6]{KOcq}, or \cite[Proposition 2.15, Lemma 2.17]{M}. Here we present these results in a way adapted to our purposes.

We denote by $\Iso(X)_+\subset\Iso(X)$ the stabiliser of $o\in X$.
\begin{lm}\label{Piso}
The isometry group of $X$ is isomorphic to $(\hat G\rtimes \Iso(X)_+)/\hat G_+$.

Suppose that there is a subgroup $P_0$ of $\Iso(X)_+$ such that $\Iso(X)_+=P_0\cdot\hat G_+$ and $\hat G_+\cap P_0=\{e\}$, then $\Iso(X)\cong\hat G\rtimes P_0$, where 
$P_0$ acts on $\hat G$ by conjugation. 
\end{lm}
\proof 
The group $\Iso(X)$ acts on $\hat G$ by conjugation. The homomorphism $\hat G\rtimes \Iso(X)_+ \rightarrow \Iso(X)$, $(g,p)\mapsto gp$ is surjective. Its kernel equals $\{(g^{-1},g)\in \hat G\rtimes \Iso(X)_+\mid g\in \hat G \}$ $\cong \hat G_+$, where the latter isomorphism is given by $(g^{-1},g)\mapsto g$. Hence we have $\Iso(X)\cong (\hat G\rtimes \Iso(X)_+)/\hat G_+$, where $g\in \hat G_+$ acts by $(\hat g,p)\mapsto(\hat g pg^{-1} p^{-1}, pg)$. 
As for the second statement, the map $\hat G\rtimes P_0\to \Iso(X)$, $(g,p_0) \mapsto g p_0$ is an isomorphism. \qed

We denote by $\Aut(\hat \fg,\theta,\ip)$ the group of automorphisms of the symmetric triple $(\hat \fg,\theta,\ip)$, i.e., of isometric Lie algebra automorphisms of $\hat \fg$ that commute with $\theta$. Similarly, let $\Aut(\hat G, \theta,\ip)$ be the group of Lie group automorphisms of $\hat G$ that commute with the involution $\theta$ of $\hat G$ and preserve the biinvariant metric on $\hat G$ associated with $\ip$. Since $\hat \fg_+=[\fm,\fm]$, an element $\phi$ of  $\Aut(\hat \fg,\theta,\ip)$ is determined by its restriction to $\fm\subset\hat\fg$.

\begin{lm} \label{Liso2}
The following maps are isomorphism of Lie groups:
\begin{eqnarray}
\Psi_1:&&\Iso(X)_+ \longrightarrow \Aut(\hat \fg,\theta,\ip), \quad  f\longmapsto \phi, \ \phi|_\fm=df_o, \label{iso1}\\
\Psi_2:&&\Aut(\hat G, \theta,\ip) \longrightarrow \Aut(\hat \fg,\theta,\ip),\quad \ph\longmapsto \ph_*=d\ph_e, \\
\Psi:&& \Iso(X)_+ \longrightarrow  \Aut(\hat G, \theta,\ip),\quad f\longmapsto C_f,\ C_f(g)=fgf^{-1}.
\end{eqnarray}
Moreover, $\Psi=\Psi_2^{-1}\circ\Psi_1$, and the inverse of $\Psi$ is given by $\ph\longmapsto f_\ph,\, f_\ph(g\hat G_+)=\ph(g)\hat G_+$. 
\end{lm}
\proof The map $\Psi_1$ is well defined. Indeed, $df_o:\fm\to\fm$ preserves $\ip|_{\fm\times\fm}$ and the Riemann curvature tensor, which is given by $R(X,Y)Z=[[X,Y],Z]$. Hence it uniquely extends to an automorphism $\phi$ of $(\hat \fg,\theta,\ip)$. Obviously, $\Psi_1$ is injective.

Clearly, $\Psi_2$ is an injective homomorphism. We show that it is also surjective. We consider the universal cover $\tilde G$ of $\hat G$. Let $\tilde \theta$ be the involution of $\tilde G$ corresponding to $\theta:\hat\fg\to\hat\fg$ and $\ip^\sim$ the biinvariant metric on $\tilde G$ corresponding to $\ip$ on $\hat\fg$. Then $\Aut(\hat \fg,\theta,\ip)$ is isomorphic to $\Aut(\tilde G, \tilde\theta,\ip^\sim)$. Furthermore, each element $\tilde \phi$ of $\Aut(\tilde G, \tilde\theta,\ip^\sim)$ descends to an automorphism $\phi$ of $\hat G =\tilde G/N$. Indeed, since $\tilde\phi$ commutes with $\tilde \theta$, it preserves $\tilde G_+$. Since it also preserves the centre of $\tilde G$, it preserves $N$. Of course, $\phi$ is an isometry on $(\hat G,\ip)$ and commutes with $\theta:\hat G\to\hat G$ since this is true on the Lie algebra level. Consequently, $\Psi_2$ is an isomorphism.

Let us first verify that $\Psi=\Psi_2^{-1}\circ\Psi_1$. We have to show that $\Psi_2(\Psi(f))=\Psi_1(f)\in\Aut(\hat \fg,\theta,\ip)$ holds for any $f\in\Iso(X)_+$. As explained above, it is sufficient to verify $\Psi_2(\Psi(f))|_\fm=\Psi_1(f)|_\fm$. Since $\Psi_2(\Psi(f))$ and $\Psi_1(f)$ preserve $\fm$, this is equivalent to $d\pi_e\Psi_2(\Psi (f))=d\pi_e\Psi_1(f)$, where $\pi$ denotes the projection $\pi:\hat G\to X$, $g\mapsto g(o)$. We have $d\pi_e\Psi_2(\Psi(f))=d\pi_e\circ (dC_f)_e=d(\pi\circ C_f)_e=d( f\circ\pi)_e$ since $f(o)=o$. Because of $d(f\circ \pi)_e=df_o\circ d\pi_e=d\pi_e\Psi_1(f)$ we are done (note that in (\ref{iso1}) we already used the identification $d\pi_e:\fm \to T_oX$). The map defined in the last line of the lemma is a right inverse of $\Psi$. Indeed, the fact that $\Psi(f_\ph)(g)=f_\ph g f_\ph^{-1}\in \hat G$, and the identity $f_\ph g f_\ph^{-1}(g_0 \hat G_+) =f_\ph(g(f_\ph^{-1} (g_0\hat G_+)))=f_\ph(g\ph^{-1} (g_0)\hat G_+)=\ph(g)g_0\hat G_+$ for all $g_0\in\hat G$, prove that $\Psi(f_\ph)(g)=\ph(g)$ holds for all $g\in G_+$. This shows that $\Psi$ is surjective, which also implies the surjectivity of $\Psi_1$. Since we already know that $\Psi$ is injective, it is an isomorphism and the right inverse discussed above equals $\Psi^{-1}$.   \qed

Under the isomorphism $\Psi$, the action of $\Iso(X)_+$ on $\hat G$ by conjugation turns into the natural action of $\Aut(\hat G, \theta,\ip)$ on $\hat G$. We want to use the isomorphisms from Lemma \ref{Liso2} to reformulate Lemma \ref{Piso} such that it becomes more suitable for the computation of the isometry group.

\begin{co} \label{co2.4}
Assume that there exists a subgroup $\bar P_0\subset\Aut(\hat\fg,\theta,\ip)$ such that $\Aut(\hat\fg,\theta,\ip)=\bar P_0 \cdot \Ad(\hat G_+)$ and $\bar P_0 \cap \Ad(\hat G_+)=\{e\}$. Then 
$\Iso(X)=\hat G \rtimes \tilde P_0$
for $\tilde P_0:=\Psi_2^{-1}(\bar P_0)$.\end{co}
\proof Note that $\Psi_1(\hat G_+)=(\Psi_2\circ\Psi)(\hat G_+)=\Ad(\hat G_+)$. Thus we can apply Lemma \ref{Piso} to $P_0:=\Psi_1^{-1}(\bar P_0)$.  Now we use the isomorphism $\Psi$ to identify $P_0$  with $\tilde P_0$. \qed
\section{Extrinsic symmetric spaces} \label{Sext}
In this section we review the notion of an extrinsic symmetric space, and prove that every pseudo-Hermitian and every para-Hermitian symmetric space is extrinsic symmetric up to a covering. We will describe an explicit construction for obtaining an extrinsic symmetric space from the symmetric triple of such spaces. Some of the symmetric spaces of signature (2,2) are pseudo-Hermitean or para-Hermitean, and realizing them as extrinsic symmetric spaces will provide us with explicit models of those spaces.

Let $(V,\ip)$ be a pseudo-Euclidean vector space and $M\hookrightarrow V$ be a non-degenerate submanifold. Here `non-degenerate' means that the restriction of $\ip$ to each tangent space $T_xM$, $x\in M$, is non-degenerate, where we understand $T_xM$ as a linear subspace of $V$. For $x\in M$, we denote by ${\cal N}_x M$ the normal space at $x$, this space being understood as an affine subspace of $V$.  The reflection $s_{x}$ of $V$ at ${\cal N} _xM$ is the affine map $s_x:V\rightarrow V$ defined by 
\begin{equation}\label{Es}
ds_{x}|_{T_xM}=-\id,\quad s_{x}|_{{\cal N}_xM}=\id.
\end{equation} 
A non-degenerate submanifold  $M\hookrightarrow V$ is called {\em extrinsic symmetric} or {\em symmetric submanifold} if $s_x(M)=M$ for all $x\in M$. Such manifolds were introduced and first studied by \cite{Fe} in the Euclidean context. Later on the notion was generalised to the case of a pseudo-Euclidean ambient space by Naitoh \cite{Nai}, and Kim and Eschenburg \cite{EK}. In the pseudo-Euclidean case, the infinitesimal object that is related to an extrinsic symmetric space is a so-called weak extrinsic symmetric triple. A detailed study of the correspondence between extrinsic symmetric spaces and their infinitesimal objects can be found in \cite{Kext}. Here, however, we only need the following stronger version of these infinitesimal objects and their role in determining an extrinsic symmetric space.
 
An {\em extrinsic symmetric triple} $(\tilde\fg,\ip\tilde{},\Phi)$ consists of a metric Lie algebra $(\tilde\fg,\ip\tilde{}\,)$ and a pair $\Phi=(D,\tilde\theta)$, where $\tilde\theta\in\Aut(\tilde\fg)$ is an isometric involution, and $D\in\fso(\tilde\fg)$ is an anti-symmetric derivation 
satisfying $D\tilde\theta=-\tilde\theta D$, $D^3=-D$  and $
[\fk_-,\fk_-]=\fk_+$, where 
$$\fk=\{X\in\tilde\fg \mid \tilde\theta(X)=X\},\quad \fk_+:=\fk\cap\Ker D,\quad \fk_-:=\{X\in\fk\mid D^2(X)=-X\}. $$

Let $(\tilde\fg,\ip\tilde{},\Phi)$ be an extrinsic symmetric triple with $\Phi=(D,\tilde\theta)$. 
We consider $$V:=\{X\in \tilde\fg\mid \tilde\theta(X)=-X\}$$ together with the restriction of $\ip\tilde{}$ to $V$ as a pseudo-Euclidean space. The isometry group of $V$ equals $\Iso(V)=\grO(V)\ltimes V$.  Let $\phi$ be the homomorphism from $\fk$ to the Lie algebra $\frak{iso}(V)=\fso(V)\ltimes V$ of $\Iso(V)$ defined by
$$\phi(u)=\big((\ad u)|_{V},-D(u)\big),\quad u\in\fk.$$ 
Note that $\phi$ is injective. 
Let $K\subset \Iso(V)$ be the connected Lie subgroup with Lie algebra $\phi(\fk)$.
Then the submanifold
$$M_{\tilde\fg,\Phi}:=K(0) =\{k(0)\mid k\in K\} \subset V$$ 
is extrinsic symmetric in $V$
and $K$ is the \emph{extrinsic transvection group}, i.e. the subgroup of $\Iso(V)$ generated by all compositions $s_x\circ s_y$ for $x,y\in M_{\tilde\fg,\Phi}$ \cite[Prop.\,4.5]{Kext}. Of course it is also an ordinary symmetric space. The Lie algebra of its (ordinary) transvection group coincides with $\fk$, see \cite[Remark 4.6]{Kext}. The decomposition $\fk=\fk_+\oplus\fk_-$ defines an involution $\tau$ by $\tau|_{\fk_+}=\id$, $\tau|_{\fk_-}=-\id$. Together with the restriction of $\ip\tilde{}$, these objects constitute the symmetric triple $(\fk,\tau,\ip\tilde{}\,|_\fk)$ of the (ordinary) symmetric space $M_{\tilde\fg,\Phi}$.

If $D=\ad (\xi)$ is an inner derivation, then we can also proceed as follows. Let $K'\subset \SO(V)$ be the connected Lie subgroup with Lie algebra $\{\ad(u)|_V\mid u\in\fk\}\subset\fso(V)$. 
Then the submanifold \begin{equation}\label{ext'}
M_{\tilde\fg,\xi}:= K'( \xi)  =\{k'(\xi)\mid k'\in K'\}  \subset V
\end{equation}
is extrinsic symmetric in $V$. 
It differs from $M_{\tilde\fg,\Phi}$ by a translation by $\xi$.

Let $X$ be a pseudo-Hermitian symmetric space and and $(\fg, \theta,\ip)$ its symmetric triple. The K\"ahler structure on $X$ corresponds to an antisymmetric derivation $J$ of $\fg$ such that $\theta J=J\theta$ and $J^2|_{\fg_-}=-\id$ and $J|_{\fg_+}=0$. We will show that, up to a covering map, $X$ is extrinsic symmetric in $V:=\fg$. This is well known for (Riemannian) Hermitian symmetric spaces, see for example \cite{Esch} and references therein. As we shall see, this is also true in the pseudo-Riemannian case, where $J$ is not necessarily an inner derivation. Moreover, we will obtain a similar statement for para-Hermitian spaces.

\begin{pr}\label{APext}
Let $X$ be a pseudo-Hermitian symmetric space and $(\fg,\theta,\ip)$ its symmetric triple. Let $J\in\frak{der}(\fg)$ be the derivation corresponding to the K\"ahler structure on $X$. If $K$ denotes the connected subgroup of $\SO(\fg)\ltimes \fg$ with Lie algebra $\{(\ad(u),-Ju)\mid u\in\fg\}\subset \fso(\fg)\ltimes \fg$, then $M:=K(0)$ is extrinsic symmetric in $(\fg,\ip)$, and the universal covers of $M$ and $X$ are isometric.

Suppose that the derivation $J$ is given by $J=\ad(\xi)$ for some $\xi\in\fg$. Let $K' \subset \SO(\fg)$ be the connected subgroup with Lie algebra $\ad(\fg)\subset \fso(\fg)$. Then, $M':= K'(\xi)$ is extrinsic symmetric in $(\fg,\ip)$, and the universal covers of $M'$ and $X$ are isometric.
\end{pr}

\proof
From the given data we will construct an extrinsic symmetric triple. We define $\tilde \fg$ as the Lie algebra direct sum $\fg\oplus\fg$, and endow $\tilde\fg$ with the inner product $\ip\tilde{}=\frac{1}{2}\big(\ip\oplus\ip \big)$, and with the involution $\tilde\theta\in \Aut(\tilde\fg)$ defined by $\tilde\theta(u,v)=(v,u)$. Lastly, define the derivation $D$ via $D(u,v):=(Ju,-Jv)$. The triple $(\tilde\fg,\ip\tilde{},\Phi)$ with $\Phi=(D,\tilde\theta)$ constitutes an extrinsic symmetric triple. Then, $V=V_+\oplus V_-$ for $V_\pm=\{(u,-u)\mid u\in\fg_\pm\}$, and $\fk=\fk_+\oplus\fk_-$ for $\fk_\pm=\{(u,u)\mid u\in\fg_\pm\}$.
 
The subspace $V$ is isometrically isomorphic to $\fg$ via the map $V \to \fg, \, (u,-u)\mapsto u$. Under this isomorphism, the group $K\subset \SO(\fg)\ltimes\fg$ defined in the proposition becomes the connected subgroup of $\SO(V)\ltimes V$ with Lie algebra $\{(\ad(\bar u),-D\bar u)\mid \bar u=(u,u)\in\fk\} \subset \fso(V)\ltimes V$. Therefore, the extrinsic symmetric space $M_{\tilde\fg,\Phi}$ is mapped to $M:=K(0)$ by the linear isometry $V\cong \fg$. Thus, $M$ is extrinsic symmetric in $(\fg,\ip)$. The symmetric triple of the symmetric space $M$ is equal to $(\fk,\tau,\ip\tilde{}\,|_\fk)$, with the involution $\tau$ determined by the decomposition $\fk=\fk_+\oplus\fk_-$. This symmetric triple is isomorphic to the symmetric triple $(\fg,\theta,\ip)$ of $X$ via $\fk\to\fg, \, (u,u)\mapsto u$. Thus the symmetric spaces $M$ and $X$ have the same universal cover. 

Moreover, if $J=\ad(\xi)$, then $D=\ad(\xi,-\xi)$. The extrinsic symmetric space $M_{\tilde\fg,(\xi,-\xi)}$ is mapped to $M':=K'(\xi)$ via the same linear isometry $V\cong\fg$, again implying that $M'$ is extrinsic symmetric and has the same symmetric triple as $X$.
\qed

Another structure on pseudo-Riemannian manifolds related to K\"ahler structures are \emph{para-K\"ahler structures}, also called \emph{bi-Lagrangian structures}. A para-K\"ahler structure on a pseudo-Riemannian manifold is a parallel field $J$ of skew-symmetric endomorphisms with $J^2=\id$. The two eigenspaces of $J$ are isotropic, i.e., the metric restricted to these spaces is zero. A symmetric space $X$ with para-K\"ahler structure is also called \emph{para-Hermitian}. If $X$ is a para-Hermitian symmetric space with symmetric triple $(\fg,\theta,\ip)$, then the para-K\"ahler structure on $X$ corresponds to an anti-symmetric derivation $J$ of $\fg$ such that $\theta J=J\theta$ and $J^2|_{\fg_-}=\id$ and $J|_{\fg_+}=0$. We will show that, up to a covering map, $X$ is extrinsic symmetric in $(\fg,-\ip)$. 

\begin{pr} \label{APext2}
    Let $X$ be a para-Hermitian symmetric space and $(\fg,\theta,\ip)$ its symmetric triple. Let $J\in\frak{der}(\fg)$ be the derivation corresponding to the para-K\"ahler structure on $X$. If $K$ is the connected subgroup of $\SO(\fg)\ltimes \fg$ with Lie algebra $\{(\ad(u),-Ju)\mid u\in\fg\}\subset \fso(\fg)\ltimes \fg$, then $M:=K(0)$ is extrinsic symmetric in $(\fg,-\ip)$, and the universal covers of $M$ and $X$ are isometric.

Suppose that the derivation $J$ is given by $J=\ad(\xi)$ for some $\xi\in\fg$. Let $K' \subset \SO(\fg)$ be the connected subgroup with Lie algebra $\ad(\fg)\subset \fso(\fg)$. Then, $M':= K'(\xi)$ is extrinsic symmetric in $(\fg,-\ip)$, and the universal covers of $M'$ and $X$ are isometric.
\end{pr}

\proof
Analogously to the case of pseudo-Hermitian symmetric spaces, we will construct an extrinsic symmetric triple out of the given data. Now $\tilde\fg=\fg_{\Bbb C}$ equals the complexification of $\fg$. The inner product $\ip\tilde{}$ on $\tilde \fg$ is equal to the real part of the complex bilinear extension of $\ip$. The involution $\tilde \theta$ equals the complex conjugation on $\fg_{\Bbb C}$. Thus $\fk=\fg$ and $V=i\fg$. 
The complex linear extension of $J$ to $\tilde \fg$ is also called $J$. We define $D:=iJ$. Then $\fk_{\pm}=\fg_\pm$ and $V_\pm= i\fg_\pm$.
The triple $(\tilde\fg,\ip\tilde{},\Phi)$ with $\Phi=(D,\tilde\theta)$ constitutes an extrinsic symmetric triple.

Here, $V$ is isometrically isomorphic to $(\fg,-\ip)$ via the map $V \to \fg, \, iu\mapsto u$. Under this isomorphism $K\subset \SO(\fg)\ltimes \fg$ becomes the connected subgroup with Lie algebra $\{(\ad(u),-Ju)\mid u\in\fg\}\subset \fso(V)\ltimes V$, and $K'\subset \SO(\fg)$ the connected subgroup with Lie algebra $\ad(\fg)\subset \fso(V)$. As above one concludes that the extrinsic symmetric space $M_{\tilde \fg,\Phi}$ is isometric to $M:=K(0)$, and $M_{\tilde\fg,\xi}$ to $M':=K'(\xi)$, that $M$ and $M'$ are also extrinsic symmetric, and that their symmetric triples are isomorphic to the one of $X$, meaning that they have the same universal cover as $X$.
\qed 
\section{Symmetric spaces of signature (2,2)}
\subsection{Review of the classification result}\label{S31}
The structure of pseudo-Riemannian symmetric spaces was studied in \cite{KOsymm}.  There, a functorial assignment was constructed that associates with a pseudo-Riemannian symmetric
space $X$ a triple consisting of a Lie algebra with involution,
a semisimple orthogonal module of the Lie algebra with involution, and a quadratic cohomology class of this module.
This led to a classification scheme for indecomposable non-semisimple pseudo-Riemannian symmetric spaces. In the case of metrics of index two this scheme gives an explicit classification in the form of a list. In the following we will specialise this list to signature $(2,2)$. Let us first recall the notion of a quadratic extension.

Let $(\fl,\theta_\fl)$ be a Lie algebra with involution and let $(\rho,\fa,\ip_\fa,\theta_\fa)$ be an orthogonal $(\fl,\theta_\fl)$-module, which means that
\begin{enumerate}
\item $\rho$ is a representation of
the Lie algebra
$\fl$ on the finite-dimensional real vector space $\fa$,
\item   $\ip_\fa$ is
a scalar product on $\fa$ which satisfies
$\langle \rho(l)a,a'\rangle_\fa + \langle a, \rho(l)a'\rangle_\fa =0$
for all $l\in\fl$ and $a,a'\in \fa$,
\item  $\theta_\fa$ is an
involutive isometry of $\fa$ such that
 $\theta_\fa \circ \rho(\theta_\fl(L))  =\rho(L)\circ \theta_\fa$
for all $L\in\fl$.
\end{enumerate} We consider the vector space $$\fd:=\fl^*\oplus\fa\oplus\fl$$ and define an inner product $\ip$ and an involutive endomorphism $\theta$ on $\fd$ by
\begin{eqnarray*}
\langle z+a+l,z'+a'+l'\rangle&:=& \langle a,a'\rangle_\fa
+z(l') +z'(l) \\
\theta(z+a+l)&:=& \theta_\fl ^*(z)+\theta_\fa(a)+\theta_\fl(l)
\end{eqnarray*}
for $z,\,z'\in \fl^*$, $a,\,a'\in \fa$ and
$l,\,l'\in \fl$. Now we choose a 2-form $\alpha$ on $\fl$ with values in $\fa$ and a 3-form
$\gamma$ on $\fl$ with values in $\RR$ such that $(\alpha,\gamma)$ is a $(\theta_\fl,\theta_\fa)$-invariant quadratic cocycle, that is, $d\alpha=0$, 
$d\gamma=\textstyle{\frac12}\langle\alpha \wedge\alpha\rangle$ and $\theta_\fa\circ \theta_\fl^*\alpha =\alpha$, $\theta_\fl^*\gamma=\gamma$. The map $\la\cdot\wedge\cdot\ra$ used in the condition for $d\gamma$ denotes the composition of the usual wedge product and $\ip_\fa$. 
Then the bilinear map
$\lb:\fd\times\fd\rightarrow \fd$ defined by $[\fl^*,\fl^*\oplus\fa] =0$ and
\begin{eqnarray}
\ [l,l'] &=& \gamma(l,l',\cdot) +\alpha(l,l')+[l,l']_\fl \label{lb1}\\
\ [l,a] &=& \rho(l)a - \langle a,\alpha(l,\cdot)\rangle\\
\ [l,z]& = & \ad {}^*(l)(z) \label{lb3}\\
\ [a,a']&=&\langle\rho(\cdot)a,a'\rangle
\end{eqnarray}
for $z\in \fl^*$, $a,\,a'\in \fa$ and
$l,\,l'\in \fl$
is a Lie bracket and $\dd:=(\fd,\theta,\ip)$ is a symmetric triple.

Whenever we consider a vector space $W$ with involution $\theta$, we will denote by $W_\pm$ the $(\pm 1)$-eigenspace of $\theta$. For elements $w\in W$, we will write $w=w_++w_-$ where $w_+\in W_+$, and $w_-\in W_-$.

\begin{pr}\label{class}
If $(\hat\fg,\ip,\theta)$ is a symmetric triple associated
with an indecomposable non-semisimple symmetric
space of signature $(2,2)$, then it is isomorphic to
$\dd$ for exactly one of the data in the following list (which contains only data giving rise to such triples):
\begin{enumerate}
\item $\fl=\RR^1=\RR\cdot e_1,\ \theta_\fl=-\id_\fl$,
\begin{itemize}
\item[(i)] $\fa=\RR^4$, $\theta_\fa=\diag(1,1,-1,-1)$, $\ip_\fa=\diag(\eps_1,\eps_2,-1,1)$, $\eps_1,\eps_2\in\{1,-1\}$, \\
$\rho(e_1)(a_1,a_2,a_3,a_4)=(\eps_1 a_3, -\eps_2 \lambda a_4, a_1,\lambda a_2)$, $\lambda\in\RR_{>0}$,\\
$\alpha=0,\ \gamma=0$;
\item[(ii)] $\fa=\RR^4$, $\theta_\fa=\diag(1,1,-1,-1)$, $\ip_\fa=\diag(-1, 1, -1,1)$,  \\
$\rho(e_1)(a_1,a_2,a_3,a_4)=(-\nu a_3+a_4,\, a_3+\nu a_4,\, \nu a_1 +a_2,\,  a_1- \nu a_2)$, $\nu\in\RR_{>0}$,\\
$\alpha=0,\ \gamma=0$;
\end{itemize}
\item $\fl=\RR^2=\Span\{e_1,e_2\} ,\ \theta_\fl=-\id_\fl $,
\begin{enumerate}
\item[] $\fa=\RR=\RR \cdot a_0,\  \theta_\fa=\id_\fa,\   \la a_0,a_0 \ra=\kappa,\  \kappa=\pm1,$\\
$\rho=0 $, \\
$\alpha(e_1,e_2)=a_0$,  $\ \gamma=0$; 
\end{enumerate}
\item$\fl=\{[e_1,e_2]_\fl=e_3,\, [e_1,e_3]_\fl=-\eps e_2\},$ $\ \eps\in\{1,-1\}$,\\[0.5ex]
$\cong \left\{ \begin{array}{ll} \fso(2)\ltimes \RR^2, &\mbox{if } \eps=1,\\  \fso(1,1)\ltimes \RR^2, &\mbox{if } \eps=-1,\end{array}\right. $ \\[0.5ex]
$\theta_\fl=\diag(-1,-1,1)$, 
\begin{itemize}
\item[] $\fa=0$, \\
$\alpha=0,\ \gamma(e_1,e_2,e_3)=\kappa, \ \kappa\in\{1,-1\}$;
\end{itemize}
\item $\fl=\{[e_1,e_2]_\fl=e_3,\ [e_1,e_3]_\fl=-e_2,\ [e_2,e_3]_\fl=\eps e_1\}$, $\ \eps\in\{1,-1\}$ \\[0.5ex] 
$\cong \left\{ \begin{array}{ll} \fsu(2), &\mbox{if } \eps=1,\\  \fsl(2,\RR), &\mbox{if } \eps=-1,\end{array}\right. $ \\[0.5ex] 
$\theta_\fl=\diag(1,-1,-1)$, 
\begin{enumerate}
\item[] $\fa=0$, \\
$\alpha=0,\  \gamma(e_1,e_2,e_3)=  c,\  c\in\RR$. 
\end{enumerate}
\item $\fl=\fsl(2,\RR)=\{[e_1,e_2]_\fl=e_3,\ [e_1,e_3]_\fl=-e_2,\ [e_2,e_3]_\fl=- e_1\}$, \\[0.5ex]
$\theta_\fl=\diag(-1,1,-1)$, 
\begin{enumerate}
\item[] $\fa=0$, \\
$\alpha=0$,
$\gamma(e_1,e_2,e_3)=c,\ c\in\RR$.
\end{enumerate}
\end{enumerate}
\end{pr}

\proof We reduce \cite[Theorem 7.10]{KOsymm} to the case of signature $(2,2)$ using the notation from \cite{KOsymm}. 

If $\fl=\RR$,  then $\fa_-$ is two-dimensional and the restriction of $\ip_\fa$ to $\fa_-$ has signature $(1,1)$. So we have the following five cases:
\begin{itemize}
    \item $\fa=\fa_+^{0,1}\oplus\fa_-^{1,0}  \oplus \fa_+^{1,0}\oplus\fa_-^{0,1}$,  $\rho=\tilde\rho_1^-  \oplus \tilde\rho_\lambda^+$ for $\lambda\in\RR_{>0}$,  
    \item $\fa=\fa_+^{0,1}\oplus\fa_-^{1,0}  \oplus \fa_+^{0,1}\oplus\fa_-^{0,1}$,  $\rho=\tilde\rho_1^-  \oplus\rho^+_\mu$ for $\mu\in\RR_{>0}$,  
\item $\fa=\fa_+^{1,0}\oplus\fa_-^{1,0}  \oplus \fa_+^{1,0}\oplus\fa_-^{0,1}$,
$\rho=\rho_1^-  \oplus \tilde\rho_\lambda^+$ for $\lambda\in\RR_{>0}$,
\item $\fa=\fa_+^{1,0}\oplus\fa_-^{1,0}  \oplus \fa_+^{0,1}\oplus\fa_-^{0,1}$,
$\rho=\rho_1^-  \oplus \rho^+_\mu$ for $\mu\in\RR_{>0}$,
\item $\fa=\fa_+^{1,1}\oplus\fa_-^{1,1}$,
$\rho=\rho''_{1,\nu}$, $\nu\in\RR_{>0}$. 
\end{itemize}
The first four cases merge into item 1\,(i) of the proposition. The fifth case is item 1\,(ii). Note that contrary to the claim in \cite{KOsymm} the symmetric triples for the parameters $\nu$ and $-\nu$ are isomorphic.

If $\dim \fl_-=2$, then $\fl$ is one of the Lie algebras $\RR^2$, $\fn(2)\cong \fso(2)\ltimes\RR^2$, $\fl=\fr_{3,-1}\cong\fso(1,1)\ltimes\RR^2$, the 3-dimensional Heisenberg algebra $\fh(1)$, $\fsu(2)$ or $\fsl(2,\RR)$, where $\fsl(2,\RR)$ appears with two different involutions. In all cases $\fa_-$ must be trivial to obtain a four-dimensional symmetric space. This excludes $\fl=\fh(1)$. For $\fl=\RR^2$, the condition $\fa_-=0$ reduces the possibilities for $\fa$ in item (2) of \cite[Theorem 7.10]{KOsymm} to cases 2. and 3. for $p=q=0$. These two cases result in item 2 of the proposition.
For the remaining $\fl$, \cite[Theorem 7.10]{KOsymm} shows that the condition $\fa_-=0$ implies $\fa=0$. The cases $\fl=\fso(2)\ltimes \RR^2$ and $\fl=\fso(1,1)\ltimes \RR^2$ merge into item 3 of the proposition. For $\fl\in\{\fsu(2), \fsl(2,\RR) \}$, we obtain the following. Cases (6) and (7) of \cite[Theorem 7.10]{KOsymm} give item 4 of the proposition. Finally, case (8) of \cite[Theorem 7.10]{KOsymm} gives item 5. above. \qed 

 Each of the triples in Proposition \ref{class} corresponds to a simply-connected symmetric space. Thus we obtain five families of symmetric spaces.

\begin{de}
Let us denote the simply-connected symmetric spaces associated with the symmetric triples occuring on the above list by 
\begin{enumerate}
\item 
\begin{itemize}
\item[(i)] \ $X_1(\eps_1,\eps_2,\lambda)$,  $\ \eps_1,\eps_2\in\{1,-1\}$, $\lambda\in\RR_{>0}$,
\item[(ii)] \ $X_2(\nu)$,   $\ \nu\in\RR_{>0}$, 
\end{itemize}
\item $N(\kappa)$, $\ \kappa\in\{1,-1\},$
\item $Y(\eps,\kappa)$, $\ \eps, \kappa\in\{1,-1\}$,
\item $Z(\eps,c)$, $\ \eps\in\{1, -1\}$, $c\in \RR$,
\item $Z'(c)$, $\ c\in \RR$,
\end{enumerate}
according to their item in the list and the parameters on which they depend.
\end{de}
The classification in \cite{KOsymm} is based on the fact that every symmetric triple $(\fg,\ip,\theta)$ without simple ideals can be endowed with the structure of a quadratic extension. In general, this can be done in several ways. There is, however, a canonical choice. This is characterised by the additional condition that  $\fl^*$ coincides with the canonical isotropic ideal of $\fg$. A quadratic extension with this property is called \emph{balanced}. The canonical isotropic ideal is uniquely determined by the Lie algebra structure of $\fg$ \cite[Section~5]{KOsymm}. In particular it is invariant under automorphisms  of $\fg$. This leads to the following remark, which will be useful later in the determination of the isometry group of the symmetric spaces under consideration.
\begin{re}\label{balanced}{\rm 
The data in the Theorem~\ref{class} are given in such a way that the canonical isotropic ideal of $\fg$ equals $\fl^*$. In particular, each automorphism of $(\fg,\ip,\theta)$ preserves $\fl^*$ and its orthogonal complement $\fl^*\oplus\fa$.  }
\end{re}

In \cite[Corollary 7.12]{KOsymm} all symmetric triples $(\hat\fg,\theta,\ip)$ associated
with an indecomposable pseudo-Hermitian symmetric space of (real) signature $(2,2q)$, $q\ge0$, were determined. In \cite[Theorem 4.1]{KOesi} this result was improved and also the K\"ahler structures were classified. Recall that the K\"ahler structure corresponds to an anti-symmetric derivation $J$ on $\hat\fg$ which satisfies $J\theta=\theta J$, $J^2|_{\hat\fg_-}=-\id_{\hat\fg_-}$ and $J|_{\hat\fg_+}=0$. If we reduce the list in \cite[Theorem 4.1]{KOesi} to signature $(2,2)$, we obtain the following result.

\begin{pr} \label{PHerm}
If $X$ is a simply-connected indecomposable non-semisimple symmetric
space of signature $(2,2)$, then it is pseudo-Hermitian if and only if it is isometric to $N(\kappa)$ for $\kappa=\pm1$ or to $Z(\eps,c)$ for $\eps=\pm1$ and $c\in\RR$. 

If $(\hat \fg,\theta,\ip)$ is the associated symmetric triple given as in Proposition~\ref{class}, then up to an automorphism of $(\hat \fg,\theta,\ip)$ the derivation $J$ corresponding to the complex structure is equal to 
\begin{enumerate}
    \item $(-J_\fl^*)\oplus0\oplus J_\fl:\fl^*\oplus\fa\oplus\fl\to\fl^*\oplus\fa\oplus\fl$ for $J_\fl$ defined by $J_\fl(e_1)=e_2$, $J_\fl(e_2)=-e_1$ if $X=N(\kappa)$,
    \item $(-J_\fl^*)\oplus J_\fl:\fl^*\oplus\fl\to\fl^*\oplus\fl$ for $J_\fl=\ad_\fl(e_1)$ if $X=Z(\eps,c)$.
\end{enumerate}
\end{pr}

Let us turn to para-Hermitian symmetric spaces. Recall that a para-K\"ahler structure on a symmetric space associated with $(\hat \fg, \theta,\ip)$ is equivalent to an anti-symmetric derivation $J$ on $\hat\fg$ which satisfies $J\theta=\theta J$, $J^2|_{\hat\fg_-}=\id_{\hat\fg_-}$ and $J|_{\hat\fg_+}=0$.

\begin{pr} {\rm \cite[Proposition 4.5]{KOesi}} \label{Ppara}
If $X$ is a simply-connected indecomposable non-semisimple symmetric
space of signature $(2,2)$, then it admits a para-K\"ahler structure if and only if it is isometric to $N(\kappa)$ for $\kappa=\pm1$ or to $Z'(c)$ for some $c\in\RR$. 

If $(\hat \fg,\theta,\ip)$ is the associated symmetric triple given as in Proposition~\ref{class}, then up to an automorphism of $(\hat \fg,\theta,\ip)$ the derivation $J$ corresponding to the para-K\"ahler structure is equal to 
\begin{enumerate}
    \item $(-J_\fl^*)\oplus0\oplus J_\fl:\fl^*\oplus\fa\oplus\fl\to\fl^*\oplus\fa\oplus\fl$ for $J_\fl=\diag(1,-1)$ with respect to the basis $e_1,e_2$ of $\fl=\RR^2$ if $X=N(\kappa)$,
    \item $(-J_\fl^*)\oplus J_\fl:\fl^*\oplus\fl\to\fl^*\oplus\fl$ for $J_\fl=\ad_\fl(e_2)$ if $X=Z'(c)$.
\end{enumerate}

\end{pr}

\subsection{Preparatory calculations}\label{S32}
We will see that the transvection groups for the symmetric spaces $X_1(\eps_1,\eps_2,\lambda)$, $X_2(\nu)$ and $Y(\eps,\kappa)$ have a common structure. They are extensions of a Heisenberg group by the real line.  In this section we will consider more general extensions of this type.  

Let $\fa$ be a finite-dimensional vector space, and $\omega: \, \fa \wedge \fa\to\RR$ a 2-form on $\fa$. Moreover, let $L\in\fsp (\fa,\omega)$.
\begin{de}\label{Heis}
We denote by $H(\omega)$ the Heisenberg group $H(\omega)=\RR\times \fa$ with group multiplication
\begin{equation}\label{Heism}
(z,a)\cdot (z',a')=(z+z'+\textstyle \frac12 \omega(a,a'),a+a')\end{equation}
and define the semidirect product $H(\omega)\rtimes\RR$, where $t\in\RR$ acts on $H(\omega)$ by 
\begin{equation}\label{Heisa}
t.(z,a)=(z,e^{tL}a) \; .
\end{equation}
\end{de}
Now, suppose $\hat \fg$ is a quadratic extension of $(\fl,\theta_\fl)=(\RR,-\id)$ by an orthogonal $(\fl,\theta_\fl)$-module $(\rho,\fa,\ip_\fa,\theta_\fa)$. Let $\RR$ be spanned by the vector $e_1$ and denote by $L$ the map $\rho(e_1)\in\fso (\fa)$. We define the 2-form $\omega$ on $\fa$ by $\omega=\langle L\, \cdot \, , \cdot \rangle$. It satisfies $\omega(\fa_+,\fa_+)=\omega(\fa_-,\fa_-)=0$. We may form the Heisenberg group $H(\omega)=\RR\times\fa\cong \fl^*\times \fa$, and the semidirect product $\hat G= H(\omega) \rtimes \RR$ as in Definition \ref{Heis}. Then, the Lie algebra of $\hat G$ is equal to $\hat \fg$. \\
The elements of $\hat G$ are of the form $(z,a,t)$ for $(z,a)\in H(\omega)$ and  $t\in\RR$. When it is more convenient, we will write just $t$ for $(0,0,t)\in \hat G$ and $(z,a)$ for $(z,a,0)\in \hat G$. Then, $(z,a,t)$ may also be written as $t\cdot(z,a)$, and we will usually do so.

The set $\hat G_+:=\{(0,a,0)\mid a\in \fa_+\}$ is an abelian subgroup of $\hat G$ with Lie algebra $\hat\fg_+$.

\begin{pr}\label{prep}
The map
\begin{eqnarray*} \Phi:\ X=\hat G/\hat G_+ & \longrightarrow & \RR\times \fa_-\times\RR\\
t\cdot(z,a)\cdot \hat G_+& \longmapsto & (v,x,u):=(z+\textstyle{\frac12}\omega(a_+,a_-),a_-,t)
\end{eqnarray*}
is a diffeomorphism 
with inverse map
\begin{eqnarray*} \Phi^{-1}:\ \RR\times \fa_-\times\RR & \longrightarrow & X=\hat G/\hat G_+\\
(v,x,u)& \longmapsto &u\cdot (v,x)\cdot\hat G_+ .
\end{eqnarray*}
\end{pr}
\proof
The map $\Phi$ is defined such that $\Phi(t\cdot(z,a)\cdot \hat G_+)$ is the unique representative of the coset $t\cdot(z,a)\cdot \hat G_+ \in \hat G/\hat G_+$ that lies in $\RR\times\fa_-\times\RR\subset \hat G$. Therefore, $\Phi$ is a well-defined section of the projection $\hat G \to \hat G/\hat G_+$, and $\Phi$ is smooth. Then, $\Phi$ is a bijection onto $\Phi(\hat G/\hat G_+)=\RR\times\fa_-\times \RR$, and its inverse is the restriction of the projection $\hat G \to \hat G/\hat G_+$ to $\RR\times\fa_-\times\RR$, which is smooth. This proves the claim. 
\qed

Via the diffeomorphism $\Phi$, the left-action of $\hat G$ on $\hat G/\hat G_+$ gives rise to a left-action of $\hat G$ on $\RR\times \fa_-\times\RR$, given by 
\begin{equation} l_{\hat g}(v,x,u):=\Phi (\hat g \cdot \Phi^{-1}(v,x,u)) \label{transvSect}\end{equation} for $\hat g\in\hat G$. 
In particular, $(v_0,x_0,u_0)=l_{u_0\cdot(v_0,x_0)}(0,0,0).$ Next we compute the differential of $l_{u_0\cdot(v_0,x_0)}$ at $(0,0,0)$. Since
\begin{eqnarray*}
l_{u_0\cdot(v_0,x_0)}(v,x,u)&=&\Phi\big(u_0\cdot(v_0,x_0)\cdot u\cdot(v,x)\cdot \hat G_+\big)\\
&=&  \textstyle \big(v_0+v+\frac{1}{2}\omega (e^{-uL}x_0,x)+\frac{1}{2}\omega((e^{-uL}x_0+x)_+,(e^{-uL}x_0+x)_- ),\\ &&\quad(e^{-uL}x_0+x)_- ,u_0+u\big),\\
&=&\textstyle \big(v_0+v+\omega ((e^{-uL}x_0)_+,x)+\frac{1}{2}\omega((e^{-uL}x_0)_+,(e^{-uL}x_0)_-),\\ && \quad(e^{-uL}x_0)_-+x ,u_0+u\big),
\end{eqnarray*}
we obtain
\begin{equation}\label{diffl} dl_{u_0\cdot(v_0,x_0)}|_{(0,0,0)}=\left(\begin{array}{ccc}1&0& \frac{1}{2}\omega (x_0,Lx_0) \\
0&\id_{\fa_-}& 0\\
0&0&1\end{array}\right) \in \GL(\RR\times \fa_-\times\RR).
\end{equation}
The inner product on $\fg_-$ gives rise to a metric on $\hat G/\hat G_+$, which we pull back along $\Phi^{-1}$ to a metric on $\RR \times \fa_- \times \RR$. We denote that metric by $g$. In $(0,0,0)$, $g$ is given by the restriction of $\ip$ to $\fg_- \cong T_{(0,0,0)}\RR \times \fa_- \times \RR$, and in any point $(v_0,x_0,u_0)$, we obtain $g$ by left-translating $g_{(0,0,0)}$ via $l_{u_0\cdot(v_0,x_0)}$. Using equation (\ref{diffl}), we obtain
\begin{equation} 
g=2dudv+\ip_{\fa_-}+\langle L^2x,x\rangle_{\fa} \cdot du^2, \label{metric}
\end{equation} 
where $\ip_{\fa_-}$ denotes the restriction of $\ip_\fa$ to $\fa_-\times \fa_-$.

\subsection{The spaces $X_1(\eps_1,\eps_2,\lambda)$, $X_2(\nu)$ and $Y(\eps,\kappa)$}
\begin{re}{\rm
The space $X_1(1,-1,1)$ is equal to the hyperbolic oscillator group, also called split oscillator group or Boidol's group. This group admits a biinvariant metric of signature (2,2). For more information see \cite{GK}.
}\end{re}
Let $\omega_0$ be the standard symplectic form on $\fa=\RR^4$, i.e., 
\begin{equation}\label{standardomega}
\omega_0(a,a')=a_1a_3'-a_3a_1'+a_2a_4'-a_4a_2'.\end{equation}
We denote by ${\rm H}_5$ the `usual' 5-dimensional Heisenberg group $H(\omega_0)$ and consider the
semidirect product 
%\begin{equation}\label{GL}
 \[   G_L:={\rm H}_5\rtimes\RR,\]
%\end{equation} 
where the action of $\RR$ on ${\rm H}_5$ is given by (\ref{Heisa}) for a map $L\in\fsp(4,\RR)$. Recall that $\fa=\fa_+\oplus\fa_-=\RR^2\oplus\RR^2=\RR^4$ and that we denote by $a_+$ and $a_-$ the components of $a\in\fa$ in $\fa_+$ and $\fa_-$. In particular, 
\[ \omega_0(a,a')=\la a_+,a_-'\ra_2-\la a_-,a'_+\ra_2\]
holds for the standard scalar product $\ip_2$ on $\RR^2$. It is well known that the map 
\begin{eqnarray}
M:\quad {\rm H}_5 &\longrightarrow& \Aff(\RR^3) \label{MHeis} \\
(z,a)&\longmapsto & (A,b)=\left(\begin{array}{cc|c} 1&a_+& z+\frac12 \la a_+,a_-\ra_2 \nonumber\\
0&I_2&a_-\end{array}\right)
\end{eqnarray}
is an injective homomorphism.

In the following we give an explicit description of the transvection and isometry groups of the spaces $X_1(\eps_1,\eps_2,\lambda)$, $X_2(\nu)$ and $Y(\eps,\kappa)$. A slightly different description can already be found in \cite[Proposition~4.19]{M}.
\begin{pr}\label{P1}
Let $X$ be one of the spaces $X_1(\eps_1,\eps_2,\lambda)$ or $X_2(\nu)$ for $\eps_1,\eps_2\in\{1,-1\}$ and $\lambda,\nu\in\RR_{>0}$. 
\begin{enumerate} 
\item The transvection group $\hat G$ of $X$ is isomorphic to $G_L$ for $L\in\fsp(4,\RR)$ with
$$L(a)= \left\{ \begin{array}{l}
(-\eps_1 a_3, -\eps_2\lambda^2a_4,-a_1,a_2), \hspace{7em} \mbox{ if } X=X_1(\eps_1,\eps_2,\lambda), \\[1ex]
((\nu^2-1)a_3-2\nu a_4, -2\nu a_3-(\nu^2-1)a_4 ,-a_1, a_2), 
\mbox{ if } X=X_2(\nu).
\end{array}\right.$$
The stabiliser is equal to the abelian group $\hat G_+=\{0\}\times\fa_+\subset {\rm H}_5$.
\item If $X=X_1(\eps_1,\eps_2,\lambda)$, then the isometry group of $X$ is isomorphic to 
\[ \Iso(X)=\left\{ \begin{array}{ll}
\hat G\rtimes (\grO(1,1)\times\ZZ_2), & \mbox{if $\lambda=1$ and $\eps_1\not=\eps_2$,} \\
\hat G\rtimes (\grO(1)\times \grO(1)\times\ZZ_2), & \mbox{else} .
\end{array}\right.\]

Here $\grO(1)\times \grO(1)\cong \Z_2\times \Z_2$ is understood as a subgroup of $\grO(1,1)$ and $A\in \grO(1,1)$ acts on $\hat G$ by
$A.(z,a_++a_-,t)=(z, (A^\top)^{-1}a_++Aa_-, t)$, where $A^\top$ denotes the usual transposed of $A$. The remaining $\ZZ_2$-factor corresponds to the involution $\theta$.

If $X=X_2(\nu)$, then $\Iso(X)= \hat G\rtimes (\ZZ_2\times\ZZ_2)$.
One $\ZZ_2$-factor is generated by $\theta$. The other one acts by
$\delta.(z,a,t)=(z,\delta a ,t)$ for $\delta\in \{\pm1\}=\ZZ_2$.
\end{enumerate}
\end{pr}

\proof
The symmetric triple of the space $X$ is a quadratic extension of $(\fl,\theta)=(\RR,-\id)$ by $\fa=\RR^4$. As was discussed in Section~\ref{S32}, its transvection group is given as the extension $H(\omega)\rtimes \RR$, where $\RR$ acts on $H(\omega)$ via $t.(z,a)=(z,e^{tL}a)$. In this case the map $L:\fa\to\fa$ and the 2-form $\omega$ on $\fa=\RR^4$ are given by
\begin{eqnarray*}
L(a) &=& (\eps_1 a_3, -\eps_2 \lambda a_4, a_1,\lambda a_2),\label{EL1}\\
\omega(a,a')&=&a_3a'_1-a_1a'_3-\lambda(a_4a'_2-a_2a'_4), \label{Eo1}
\end{eqnarray*}
if $X=X_1(\eps_1,\eps_2,\lambda)$ and 
\begin{eqnarray*}
L(a) &=& (-\nu a_3+a_4,\, a_3+\nu a_4,\, \nu a_1 +a_2,\,  a_1- \nu a_2),\label{EL2}\\
\omega(a,a')&=&(\nu a_3-a_4)a_1'+(a_3+\nu a_4)a'_2-(\nu a_1+a_2)a_3' +(a_1-\nu a_2)a'_4,  \label{Eo2}
\end{eqnarray*}
if $X=X_2(\nu)$, where $a=(a_1,\dots,a_4)$,  $a'=(a'_1,\dots,a'_4)\in \RR^4$. Later on we will change the basis to get $L$ and $\omega$ as claimed in the proposition.

Recall that $e_1$ spans $\fl=\RR\subset \hat \fg$ and consider
\begin{eqnarray*}
\bar P_0&:=&\{F\in\Aut(\hat\fg,\theta,\ip)\mid \pro_\fa F(e_1)=0\}\\
&=& \{F_{r,S}\mid S\in\grO(\fa),\ r=\pm1,\ S\theta_\fa=\theta_\fa S,\ SL=rLS\},
\end{eqnarray*}
where
\[F_{r,S}: \hat\fg\longrightarrow\hat\fg, \qquad z+a+l\longmapsto rz+Sa+rl\]
for $z\in\fl^*$, $a\in\fa$, $l\in\fl$. Obviously, $\bar P_0$ is a subgroup of $\Aut(\hat\fg,\theta,\ip)$. We will show that $\Aut(\hat\fg,\theta,\ip)=\bar P_0 \cdot \Ad(\hat G_+)$ holds. Take $F\in\Aut(\hat\fg,\theta,\ip)$. Recall that $F$ preserves both $\fl^*$ and $\fa\oplus\fl^*$, see Remark~\ref{balanced}. We put $S:=\pro_\fa\circ F|_\fa\in \grO(\fa)$ and $a:=\pro_\fa F(e_1)$. Then $F\circ \Ad(L^{-1}S^{-1}a)$ is in $\bar P_0$ since 
\begin{eqnarray*}
\pro_\fa (F\circ \Ad(L^{-1}S^{-1}a))(e_1) &=& \pro_\fa F (\pro_{\fa\oplus\fl}\Ad(L^{-1}S^{-1}a)(e_1))\\
&=& \pro_\fa F (e_1-S^{-1}a)\\
&=& a-S(S^{-1}(a))\ =\ 0.
\end{eqnarray*}
Hence $F\in \bar P_0 \cdot \Ad(\hat G_+)$ holds, which proves the claim. Furthermore, the intersection $\bar P_0 \cap \Ad(\hat G_+)$ is trivial. 

The group $\bar P_0$ is generated by $\theta$ and all automorphisms $F_{1,S}$ for some $S$ satisfying 
\begin{equation}\label{ES}
S\in \grO(\fa),\ SL=LS,\ S\theta_\fa=\theta_\fa S.
\end{equation}
For $X=X_1(\eps_1,\eps_2,\lambda)$, this condition is equivalent to $S_\pm:=S|_{\fa_\pm}\in \grO(\fa_\pm)=\grO(1,1)$ and \[ S_- \left(\begin{array}{cc} 1&0\\0&\lambda\end{array} \right) =\left(\begin{array}{cc} 1&0\\0&\lambda\end{array} \right) S_+,\quad S_+ \left(\begin{array}{cc} \eps_1&0\\0&-\lambda\eps_2\end{array} \right) =\left(\begin{array}{cc} \eps_1&0\\0&-\lambda\eps_2\end{array} \right) S_-.\]
In particular, $S_-\in \grO(1,1)$ commutes with $\diag(\eps_1,-\lambda^2 \eps_2)$. If $\lambda=1$ and $\eps_1=-\eps_2$, then this is satisfied for arbitrary $S_-\in\grO(1,1)$. If not, then the condition is equivalent to $S_-=\diag(\delta_1,\delta_2)$ for $\delta_1,\delta_2=\pm1$. Vice versa, every such matrix $S_-$ defines a unique map $S$ that satisfies (\ref{ES}). Thus $\bar P_0$ is isomorphic to $\grO(1,1)\times \ZZ_2$ if $\lambda=1$ and $\eps_1=-\eps_2$. Otherwise, $\bar P_0$ is isomorphic to $\grO(1)\times\grO(1)\times\ZZ_2\cong\ZZ_2\times\ZZ_2\times\ZZ_2$. In both cases, $A\in\grO(1,1)$ acts on $\hat \fg$ by
$A.(z+a+l)=z+Aa_++Aa_-+l.$

For $X=X_2(\nu)$, the condition   (\ref{ES}) is equivalent to $S_\pm:=S|_{\fa_\pm}\in \grO(1,1)$ and \[ S_- \left(\begin{array}{cc} \nu&1\\1&-\nu\end{array} \right) =\left(\begin{array}{cc} \nu&1\\1&-\nu\end{array} \right) S_+,\quad S_+ \left(\begin{array}{cc} -\nu&1\\1&\nu\end{array} \right) =\left(\begin{array}{cc} -\nu&1\\1&\nu\end{array} \right) S_-.\] 
Hence $S_-\in\grO(1,1)$ commutes with {\small
$\left(\begin{array}{cc} a&-b\\b&a\end{array} \right)$}, where $a=1-\nu^2$ and $b=-2\nu\not=0$, which is only possible for $S_-=\pm \id_{\fa_-}$. Hence, $\bar P_0$ is isomorphic to $\ZZ_2\times\ZZ_2$, where one $\ZZ_2$-factor is generated by $\theta$ and the other one by the map $z+a+l\mapsto z-a+l$.

Now we change the coordinates in $\fa$ in order to get the claimed expression for $L$ and~$\omega$. 

For $X=X_1(\eps_1,\eps_2,\lambda)$, we put $\hat a=(\hat a_1,\dots,\hat a_4)= (-a_1,\lambda a_2, a_3, a_4)$. Then $\omega(a,a')=\omega_0(\hat a,\hat a')$. Moreover, we obtain 
$\hat L(\hat a)=(-\eps_1 \hat a_3, -\eps_2\lambda^2\hat a_4,-\hat a_1,\hat a_2)$
for the matrix $\hat L$ of $L$ with respect to the new basis. Finally, we rename $(\hat a_1,\dots,\hat a_4)$ back to $(a_1,\dots,a_4)$ and $\hat L$ to $L$.  For $\lambda=1$ and $\eps_1=-\eps_2$, we have also to determine the transformation of the action of $\grO(1,1)$ on $\fa$.  Take $A\in \grO(1,1)$. Then $A.(z+a+l)=z+Aa_++Aa_-+l$ becomes $A.(z+a+l)=z+\hat A a_++Aa_-+l$ for $\hat A=\diag(-1,1) A \diag (-1,1)= (A^\top)^{-1}$.

For $X=X_2(\nu)$, we use the coordinate transformation $\hat a=(\hat a_1,\dots,\hat a_4)=({-\nu a_1-a_2,}\linebreak a_1-\nu a_2, a_3, a_4)$ and proceed as above.

By applying Corollary \ref{co2.4} we obtain the assertion of the proposition.
\qed

\begin{pr} \label{IsoY}
\begin{enumerate}
\item The transvection group $\hat G$ of $Y(\eps,\kappa)$ is isomorphic to $G_L$ for $L\in\fsp(4,\RR)$ with
\begin{equation} \label{LY}
L(a)=(a_4,a_3-\kappa a_4,\eps\kappa a_1-\eps a_2,-\eps a_1)
\end{equation}
The stabiliser is equal to the abelian group $\hat G_+=\{0\}\times\fa_+\subset {\rm H}_5$.
\item The isometry group of $Y(\eps,\kappa)$ is equal to $\hat G\rtimes (\ZZ_2\times\ZZ_2)$. One $\ZZ_2$-factor is generated by $\theta$. The other one acts by
$\delta.(z,a,t)=(z,\delta a ,t)$ for $\delta\in \{\pm1\}=\ZZ_2$.
\end{enumerate}
\end{pr}

\proof
ad 1: First, we show that the transvection algebra $\hat{\fg}$ is isomorphic to a semi-direct algebra of a Heisenberg group with the real line. This allows us to apply the results of our preparatory calculations from Section \ref{S32}. 

We consider the Lie algebra $\hat\fg$ as the quadratic extension of 
\begin{equation}\label{lY}
\fl=\{[e_1,e_2]_\fl=e_3,\ [e_1,e_3]_\fl=-\eps e_2\}
\end{equation}
 by $\fa=0$ with $\gamma(e_1,e_2,e_3)=\kappa\in\{1,-1\}$, see Section~\ref{S31}. We denote by $\sigma^1,\sigma^2,\sigma^3\in\fl^*$ the dual basis of $e_1,e_2,e_3$. We define a new basis of $\hat \fg$:
\[b_1:=\sigma^1,\ b_2:=- e_3+\textstyle \frac{\kappa}2 \sigma^3,\ b_3:=\sigma^3,\ b_4:=\eps \sigma^2,\ b_5:=-e_2+\frac{\eps\kappa}2\sigma^2,\ b_6:=e_1.\]
The non-vanishing commutators between these basis vectors are
\begin{eqnarray}
&&[b_2,b_4]=b_1,\ [b_3,b_5]=b_1, \label{HY1}\\
&& [b_6,b_2 ]= -\eps b_5+\eps \kappa b_4,\ [b_6,b_3]= -\eps b_4,\ [b_6,b_ 4]=b_3 ,\ [b_6,b_ 5]=b_2-\kappa b_3 .\label{HY2}
\end{eqnarray}
Let $\hat\fa$ denote the subspace spanned by $b_2,\dots,b_5$, and $\ip_{\hat \fa}$ denote the restriction of the inner product on $\hat \fg$ to $\hat \fa$. Furthermore, let $L \in \fso (\hat \fa)$ be the restriction of $\ad(b_6)$ to $\hat \fa$. Let $\hat \fl$ denote $\RR \cdot b_6\cong \RR$. We make $\hat \fa$  an orthogonal $(\hat \fl , -\id_{\hat \fl})$-module by introducing the representation $\rho: \, \hat \fl \to \fso(\hat \fa)$, $l=te_1\mapsto tL$, and the involution $\theta_{\hat \fa}$ which is the restriction of the involution on $\hat \fg$ to $\hat\fa$. 
This shows that $\hat \fg$ can also be considered as the quadratic extension of $(\hat \fl,- \id_{\hat \fl})$ by $(\rho,\hat \fa,\ip_{\hat \fa},\theta_{\hat \fa})$ with $\hat\gamma=0$ and $\hat\alpha=0$.

Now, $\hat \fg$ has the form of the Lie algebras considered in Section \ref{S32} (in that section, $\hat \fa$ was denoted by $\fa$, and $\hat \fl$ was denoted by $\fl$). Thus the transvection group $\hat G$ is isomorphic to the semidirect product  
$H(\omega)\rtimes \RR,$
where $H(\omega)$ is the Heisenberg group defined by the 2-form $\omega: \, \hat \fa \wedge \hat \fa \to \RR$, $\omega(a,a') :=  \langle La , a'\rangle_{\hat \fa}$ and the action of $\RR$ on $H(\omega)$ is given as in Definition \ref{Heis}. If we identify $\hat a$ with $\RR^4$ using the basis $b_2,\dots,b_5$, then (\ref{HY1}) and (\ref{HY2}) imply that $L$ is given by (\ref{LY}) and $\omega$ equals the standard symplectic form $\omega_0$ defined in~(\ref{standardomega}).

ad 2: Here we use the description of $\fg$ as a balanced quadratic extension as given in Proposition~\ref{class}. Let $A$ be an automorphism of $(\fg,\theta,\ip)$. Then $A$ preserves the centre of $\fg$. Furthermore, $\fl^*_+:=\fl^*\cap \hat \fg_+$ and $\fl^*_-:=\fl^*\cap \hat\fg_-$ are invariant under $A$ since $\fl^*$ is the canonical isotropic ideal of $\hat\fg$, see Remark \ref{balanced}. If $\sigma^1,\sigma^2,\sigma^3$ denotes the basis of $\fl^*$ that is dual to $e_1,e_2,e_3$, this implies that $A$ preserves $\RR\cdot \sigma_1$, $\RR\cdot \sigma^3$ and $\Span\{\sigma^1,\sigma^2\}$. By Eqs.~(\ref{lY}) and  (\ref{lb3}), the group $\{\Ad(\exp te_3)|_{\fl_+^*}\mid t\in\RR\}$ consists of linear maps $A_0$ satisfying $A_0(\sigma^1)=\sigma^1$ and $A_0(\sigma^2)=\sigma^2+s\sigma^1$, $s\in\RR$. Therefore, $A=A'A_0$, where $A_0\in \Ad(\hat G_+)$ and $A'|_{\fl^*}=\diag(a^{-1},b^{-1},c^{-1})$, $a,b,c\in\RR$, with respect to the basis $\sigma^1,\sigma^2,\sigma^3$. Since $A'$ is an isometry, we obtain $\bar A':=\pro_\fl A'|_\fl=\diag(a,b,c)$. Since $A$ is a Lie algebra automorphism, also the map $\bar A'\in\End(\fl)$ is an automorphism and preserves $\gamma$. Using (\ref{lY}) and the fact that $\gamma$ is a non-vanishing 3-form on $\fl$, we obtain $a,b,c\in\{1,-1\}$  and $abc=1$. Now we consider the inner automorphisms $\Ad(\exp t\sigma^3)$. They are the identity on $\fl^*\oplus \RR e_3$ and map 
\[e_1\longmapsto e_1+t\sigma_2,\quad e_2\longmapsto e_2-t \sigma^1,\quad e_3\longmapsto e_3.\]
Hence $A'=A''A_1$ for some $A_1\in\Ad(\hat G_+)$ and $A''=\diag(a^{-1},b^{-1},c^{-1},a,b,c)$ with respect to the basis $\sigma^1,\sigma^2,\sigma^3,e_1,e_2,e_3$. Now we apply Corollary~\ref{co2.4} to \[\bar P_0=\{\diag(a^{-1},b^{-1},c^{-1},a,b,c)\mid a,b,c=\pm1, abc=1\}\cong \ZZ_2\times \ZZ_2\] and we obtain the assertion.
\qed

Recall, that elements of $\hat G=G_L={\rm H}_5\rtimes \RR$ are of the form $(z,a,t)$ for $(z,a)\in  {\rm H}_5=\RR\times\RR^4$ and  $t\in\RR$ and that we also write just $t$ for $(0,0,t)\in \hat G$ and $(z,a)$ for $(z,a,0)\in \hat G$. Then, $(z,a,t)=t\cdot(z,a)$ and we will usually write elements of $\hat G$ in the latter way.
\begin{pr} For $X=X_1(\eps_1,\eps_2,\lambda), \, X_2(\nu),$ or $ Y(\eps,\kappa)$ the map
\begin{eqnarray*} \Phi:\ X=\hat G/\hat G_+ & \longrightarrow & \RR^4\\
t\cdot(z,a)\cdot \hat G_+& \longmapsto & (v,x,u):=(z+\textstyle{\frac12}\omega_0(a_+,a_-),a_-,t)
\end{eqnarray*}
is a diffeomorphism. The action of the transvection group is given by 
\[t\cdot(z,a)\cdot (v,x,u) = \big(M(z,e^{-uL}a)(v,x),u+t\big)\]
for $t\cdot(z,a)\in \RR\ltimes{\rm H_5}$, where the map $M:{\rm H}_5\to\Aff(\RR^3)$ is defined in (\ref{MHeis}). In particular, the Heisenberg group ${\rm H}_5$ acts by affine transformations on each $u$-level.
\begin{itemize}
\item[(i)] For $X=X_1(\eps_1,\eps_2,\lambda)$, the metric is given by
$$
2dudv-dx_1^2+dx_2^2+(-\eps_1x_1^2-\eps_2\lambda^2x_2^2)du^2
$$
An isometry $(\delta_1,\delta_2,\delta_3)\in \grO(1)\times \grO(1)\times \ZZ_2$ acts by 
\[(\delta_1,\delta_2,\delta_3)\cdot (v,x,u) = (\delta_3 v,\delta_1\delta_3 x_1, \delta_2 \delta_3 x_2, \delta_3 u).\]
If, moreover, $\lambda=1$ and $\eps_1=-\eps_2$, then $A\in \grO(1,1)$ acts by
\[A\cdot (v,x,u)= (v, Ax, u).\]

\item[(ii)] For $X=X_2(\nu)$, the metric is given by
$$
2dudv-dx_1^2+dx_2^2+((\nu^2-1)(x_1^2-x_2^2)-4\nu x_1x_2)du^2
$$
An isometry  $(\delta_1,\delta_2)\in \ZZ_2\times \ZZ_2$ acts by 
\[(\delta_1,\delta_2)\cdot (v,x,u) = (\delta_2 v,\delta_1\delta_2 x, \delta_2 u).\]
\item[(iii)] For $X=Y(\eps,\kappa)$, the metric is given by
$$
2dudv -\eps(2dx_1dx_2+\kappa dx_2^2)+ (2x_1x_2 - \kappa x_2^2) du^2
$$
An isometry  $(\delta_1,\delta_2)\in \ZZ_2\times \ZZ_2$ acts by 
\[(\delta_1,\delta_2)\cdot (v,x,u) = (\delta_2 v,\delta_1\delta_2 x, \delta_2 u).\]
\end{itemize}
\end{pr}
\proof 
The symmetric space $X$ is of the form considered in Section \ref{S32}. That the map $\Phi$ is a diffeomorphism and the claimed form of the metric therefore follows from Proposition \ref{prep} and equation \eqref{metric} respectively. 
The image of $X$ under the diffeomorphism $\Phi$ may be identified with a subspace of the transvection group $\hat G$ via $\RR^4=\RR\times\RR^2\times\RR\cong \RR\times \fa_- \times \RR$  (or rather $\RR\times\hat\fa_-\times\RR$ in the case $X=Y(\eps,\kappa)$). Equation \eqref{transvSect} then describes the action of the transvection group on $X$ under this identification, while the action of the isometry group on $X$ follows from Propositions \ref{P1} and \ref{IsoY}.  
\qed

\begin{re}\label{PW}{\rm
The metrics of the families $X_1, X_2, Y$ behave similarly to Lorentzian plane waves. They are of the form $2dudv+\ip+H(u,x)du^2$, where $\ip$ is pseudo-Euclidean of signature $(1,1)$. There exists a non-zero parallel vector field $V$, which is unique up to multiplication by a scalar. Indeed, the holonomy representation is given by the adjoint representation of $\hat G_+$ on $\hat\fg_-$. If $X$ belongs to one of the families $X_1,X_2$ or $Y$, then the space of invariants of this representation is exactly $\fl_-^*\cong \RR$ or $\hat\fl_-^* \cong \RR$, respectively. We have also seen that the transvection group contains the 5-dimensional Heisenberg group, which acts transitively on each leaf of the integrable distribution~$V^\perp$.   
}\end{re}

\begin{pr}{\rm \cite{M}} The space $X_1(\eps_1,\eps_2,\lambda)$ admits a compact quotient if and only if $\lambda=1$ and $\eps_1\not=\eps_2$. Compact quotients can be obtained by discrete subgroups of the transvection group. The spaces $X_2(\nu)$ and $Y(\eps,\kappa)$ do not admit compact quotients by discrete subgroups of their isometry groups.
\end{pr}
\proof The spaces $X_1(\eps_1,\eps_2,\lambda)$, $X_2(\nu)$, and $Y(\eps,\kappa)$ correspond to the symmetric spaces of Case (B) in \cite{M}. In the notation of \cite{M} the transvection group of $X_1(1,-1,1)$ is equal to $G_{I_{1,1}, I_{1,1}}$ and the one of $X_1(-1,1,1)$ equals $G_{-I_{1,1}, I_{1,1}}$. The assertion now follows from \cite[Theorem 1.7]{M}. \qed

The work \cite{M} does not contain an explicit description of the compact quotients of $X_1(\pm 1,\mp 1,1)$. However, neither Maeta's work nor ours aims to provide such descriptions. It does not appear to be too difficult to obtain an explicit description of compact quotients of $X_1(1,-1,1)$ along the lines indicated in \cite{Malek}. There it is shown that compact quotients of $X_1(1,-1,1)$ are isometric to a quotient of the hyperbolic oscillator group (endowed with its bi-invariant metric) by a cocompact lattice or it is isometric to a quotient $\Gamma \backslash (\R\times{\rm H}_3)$, where $\Gamma$ is a cocompact lattice of $\R\times{\rm H}_3$ (with a suitable left-invariant metric). Lattices in the hyperbolic oscillator group are classified in \cite{GK}.

The following proposition shows that the conditions for the action of the discrete subgroup may be relaxed, it is not necessary to assume that it acts freely. We will show that the isometry groups of symmetric spaces of type $X_1$, $X_2$ and $Y$ are linear and apply Selberg's Lemma.
\begin{pr}\label{SolvSelberg}
Let $X$ be one of the spaces $X_1(\eps_1,\eps_2,\lambda)$, $X_2(\nu)$, and $Y(\eps,\kappa)$. Every discrete subgroup of the isometry group of $X$ acting properly and cocompactly on $X$ admits a finite index subgroup that also acts freely.  
\end{pr}
\proof We want to apply Selberg's Lemma. Let us therefore first prove that the isometry group of $X$ is virtually linear, i.e. that it contains a finite index subgroup for which there is an injective homomorphism into a matrix group. By Proposition \ref{P1}, the isometry group contains a connected finite index subgroup of the form ${\rm H}_5\rtimes_\rho Q$, where $Q=\RR$ or $Q=\RR\times\SO_0(1,1)$. In particular, $Q$ is a matrix group. We define a homomorphism ${\ph: \rm H}_5\rtimes_\rho Q \to {\rm H}_5\rtimes \Aut({\rm H}_5)$ by $\ph(h,q)= (h,\rho(q))$. Since the representation $\rho$ is faithful, $\varphi$ is injective. Now we use that ${\rm H}_5\rtimes \Aut({\rm H}_5)$ is linear. Indeed, the adjoint representation of this group on its Lie algebra $\tilde\fh$ is faithful since the centre of ${\rm H}_5\rtimes \Aut({\rm H}_5)$ is trivial. Thus it defines an injective homomorphism into $\GL(\tilde\fh)$. Consequently, 
$${\rm H}_5\rtimes_\rho Q \ni (h,q) \longmapsto (q,\Ad (\ph(h,q))) \in Q\times\GL(\tilde\fh) $$ 
is an injective homomorphism into a matrix group. 

Now, let $\Gamma$ be a discrete subgroup of the isometry group of $X$ acting properly and cocompactly. We may assume that $\Gamma$ is contained in the linear finite index subgroup of the isometry group. Then, Selberg's Lemma provides a torsion-free subgroup of $\Gamma$ with finite index, which then acts freely, see Lemma \ref{freeAction}.
\qed

\subsection{The spaces $N(\kappa)$}

The underlying Lie algebra of the symmetric triple is nilpotent and equals $\hat \fg=\fl^*\oplus \fa \oplus \fl=\RR^2\oplus\RR\oplus\RR^2$ as a vector space. As usual we write the elements of $\fl^*\oplus \fa \oplus \fl$ as $z+a+l$, with $z\in \fl^*$, $a\in\fa$, $l\in\fl$. The two-form $\alpha$ used to define the Lie bracket is equal to the canonical symplectic form on $\RR^2$, i.e. $\alpha(l,\hat l)=l_1\hat l_2-l_2\hat l_1$ for $l=(l_1,l_2)$, $\hat l=(\hat l_1,\hat l_2)\in\RR^2$. 

First we will determine the transvection group $\hat G$ of $N(\kappa)$ by integrating the Lie bracket of $\hat \fg$. The resulting form of the group multiplication has the advantage that the group $\SL^\pm(2,\RR)$ acts on $\hat G$ in a natural way. This leads to a workable description of the full isometry group. Later, in Corollary~\ref{coNext}, we will obtain an alternate description of $\hat G$ as a subgroup of stricly upper triangular matrices of rank 6. Then the group multiplication is easier to understand. However, the action of $\SL^\pm(2,\RR)$ takes a more complicated form.

It would be possible to write $\hat G$ in the form $({\rm H}_3\times\RR)\rtimes \RR$, where ${\rm H}_3$ denotes the three-dimensional Heisenberg group. Similar calculations as in Section~\ref{S32} would then again lead to a metric of the form $2dudv+\ip+H(u,x)du^2$, where $\ip$ is pseudo-Euclidean of signature $(1,1)$. However, doing so seems unnatural. Indeed, an analogous argument to that in Remark~\ref{PW} shows that the spaces $N(\kappa)$ admit a two-dimensional space of parallel vector fields, since $\fl_-^*=\RR^2$ is the space of invariants of $\Ad(\hat G_+)|_{\hat\fg_-}$. Writing $\hat G$ as $({\rm H}_3\times\RR)\rtimes \RR$ would then arbitrarily distinguish one of the parallel vector fields.

A description of the isometry group of $N(\kappa)$ can already be found in \cite[Proposition~4.16]{M}.
\begin{pr}\label{PN}\begin{enumerate}
\item The transvection group $\hat G$ of $N(\kappa)$ is isomorphic to the group 
$\RR^2\times\RR\times\RR^2$ with multiplication
\begin{equation} \label{GrpMult}   \textstyle(z,a,l)\cdot (\hat z,\hat a,\hat l\hspace{1pt}) =\left(z+\hat z+\frac \kappa 3 \alpha(l,\hat l\hspace{1pt})(l+\frac12 \hat l\hspace{1pt})+\kappa\hat a l,\,a+\hat a+\frac12 \alpha(l,\hat l\hspace{1pt}), \, l+ \hat l\right). 
\end{equation}
Under this isomorphism, the stabiliser 
$\hat G_+$ is equal to $\{(0,\hat a,0)\mid \hat a\in\RR\}$.
\item The isometry group of $N(\kappa)$ is isomorphic to $\hat G\rtimes \SL^\pm(2,\RR)$.  The action of $S\in\SL^\pm(2,\RR)$ on $\hat G$ is given by
$(z,a,l)\mapsto (|S|S z, |S|a, Sl)$. 
\end{enumerate} 
\end{pr}
\proof ad 1. For $z\in\fl\cong\RR^2$, we define $z^\flat:=\alpha(z,\cdot)$. Using a (pseudo-) unit vector we identify $\fa\cong \RR$. Moreover, we identify $\fl\cong\fl^*, z\mapsto -z^\flat$, thus $\fl^*\oplus\fa\oplus\fl\cong\fl\oplus\fa\oplus\fl$. Then
\begin{eqnarray}\label{lb2} 
[z+a+l\,,\, \hat z+\hat a+\hat l\,]&=& \kappa(\hat a\cdot l-a\cdot \hat l\hspace{1pt})+\alpha(l,\hat l\hspace{1pt})\ \in\ \fl\oplus\fa\subset \fl\oplus\fa\oplus\fl,\\
\la z+a+l\,,\,\hat z+\hat a+\hat l\,\ra&=&-\alpha(z,\hat l\hspace{1pt})-\alpha(\hat z,l)+\kappa a\hat a. \label{ipNk}
\end{eqnarray}
We describe the associated simply-connected Lie group. The underlying set is equal to  $\fl\times\fa\times\fl$.
If we put $h(z,a):=(z,a,0)$ and $\lambda(l):= (0,0,l)$, then \cite[Lemma 2.3]{gitter2} gives
\begin{eqnarray*}
h(z,a)\cdot h(\hat z,\hat a)&=&\textstyle{h(z+\hat z,a+\hat a)  }\\
\lambda(l)\cdot \lambda(\hat l\hspace{1pt}) &=& h\left(\textstyle{\frac\kappa3\alpha(l,\hat l\hspace{1pt})(l+\frac12 \hat l\hspace{1pt}),\frac12 \alpha(l,\hat l\hspace{1pt})}\right) \cdot \lambda(l+\hat l\hspace{1pt}) , \\
\lambda(l)\cdot h(z,a)\cdot \lambda(l)^{-1}&=&h(z+\kappa a l,a), 
\end{eqnarray*}
which implies \eqref{GrpMult}. The claimed expression for $\hat G_+$ is clear. The group $\hat G$ indeed acts effectively on $\hat G/\hat G_+$. Hence $\hat G$ is the transvection group of $N(\kappa)$.

ad 2. Let us first remark that
\[\Ad(h(0,\hat a))(z+a+l)=(z-\kappa \hat a l)+a+l\in\fl\oplus\fa\oplus \fl\] implies 
\begin{equation}\label{AdG} 
\Ad(\hat G_+)= \{ z+a+l\mapsto(z- rl)+a+l\mid r\in\RR \}.
\end{equation}
An endomorphism $F$ of $\fg$ is a Lie algebra automorphism preserving $\theta$ if and only if
\[F(z+a+l)=(|S|^{-1}\cdot Sz+T'l)+|S|a+Sl \in\fl\oplus\fa\oplus\fl\] 
for some $S\in\GL(\fl)\cong\GL(2,\RR)$ and $T'\in\End(\fl)$. 

If, in addition, $F$ is an isometry, then $|S|=\pm1$. In particular, $|S|^{-1}=|S|$.
Hence $F$ equals the composition $F_1\circ F_2$, where
\[F_1(z+a+l)=(z+Tl)+ a+l,\quad F_2(z+a+l)=|S| Sz+|S|a+Sl,\]
where $T=T'S^{-1}$. The map $F_2$ is an automorphism of $(\hat\fg,\theta,\ip)$, in particular it is an isometry. Hence $F_1$ is also an isometry. Consequently, $F_1(l)=Tl+l\in\fl\oplus\fl$ is isotropic, which implies $\alpha(Tl,l)=0$.
Thus, $Tl$ is a multiple of $l$, and since $T$ is linear, it follows that $T$ is a multiple of the identity. Consequently, $F_1$ is an element of $\Ad(\hat G_+)$ by \eqref{AdG}. Set
\begin{equation}\label{A0}
\bar P_0:=\{z+a+l\mapsto |S| Sz+|S|a+Sl\mid S\in{\rm SL}^{\pm}(2,\RR)\} \subset \Aut(\hat\fg,\theta,\ip).
\end{equation}
We have seen that $\Aut(\hat\fg,\theta,\ip)=\Ad(\hat G_+)\cdot\bar P_0$. Moreover, the intersection of $\bar P_0$ and $\Ad(\hat G_+)$ is trivial. Corollary \ref{co2.4} therefore yields that the isometry group of $N(\kappa)$ equals $\hat G \rtimes \tilde P_0$
for $\tilde P_0:=\Psi_2^{-1}(\bar P_0)$. Clearly, $\tilde P_0 \cong \SL^\pm(2,\RR)$ and the action is as claimed in the proposition.
\qed

\begin{pr}
The map
\begin{eqnarray*}
\Phi: \; N(\kappa)=\hat G/\hat G_+ & \longrightarrow & \RR^2\times \RR^2 \\
(z,a,l) \hat G_+&\longmapsto & (v,u):=(z-\kappa al, l)
\end{eqnarray*}
is a diffeomorphism. Under this diffeomorphism the metric on $N(\kappa)$ becomes
\[\textstyle -\frac \kappa 3(u_2du_1-u_1du_2)^2+2(du_1dv_2-du_2dv_1). \]
The action of the isometry group $\hat G\rtimes \SL^\pm(2,\RR)$ is given by
\[ \textstyle (z,a,l)\cdot S\cdot (v,u)=\left( |S|Sv+z- \frac \kappa 3 \alpha(l,Su)(Su+ \frac 12  l) -\kappa a(Su+l), Su+l\right).\]
\end{pr}

\proof
Two elements $(z,a,l),(\hat z,\hat a,\hat l) \in \hat G$ lie in the same coset in $\hat G/\hat G_+$ if and only if there exists an $a'\in\RR$ such that $(\hat z,\hat a,\hat l)=(z,a,l)(0,a',0)=(z+\kappa a' l, a+a',l)$. It follows that the map 
$$
\hat G / \hat G_+ \longrightarrow \hat G \; , \; (z,a,l)\hat G_+ \longmapsto (z-\kappa al,0,l)
$$
is a well-defined section of the quotient map $\hat G \to \hat G/\hat G_+$. In particular, 
it is injective, and thus a diffeomorphism onto its image $\RR^2\times\{0\}\times\RR^2\subset \hat G$. Identifying $\RR^2\times\{0\}\times\RR^2\cong \RR^2\times\RR^2$ gives the first claim. 

Let us now turn to the metric. The identity element of $\hat G$ is $e=(0,0,0)$. The metric of $\hat G/\hat G_+$ is left-invariant and at $e\hat G_+$ given by the restriction of $\ip$ to $\hat\fg_-\cong T_{e\hat G_+}(\hat G/\hat G_+)$. This metric is then pushed forward along $\Phi$ to obtain a metric on $\RR^2\times\RR^2$. According to Proposition~\ref{PN}, we have coordinates $(z_1,z_2,a,l_1,l_2)$ on $\hat G$. The coordinate vector fields $\partial_{z_1}, \partial_{z_2},\partial_{l_1}, \partial_{l_2}$ at the identity $e$ span $\hat\fg_-$. We use the same notation for their projection to $T_{e\hat G_+}(\hat G/\hat G_+)$. By \eqref{ipNk}, we have
$$
\langle \partial_{z_j},\partial_{z_k}\rangle =\langle \partial_{l_j},\partial_{l_k}\rangle=\langle \partial_{z_j},\partial_{l_j}\rangle =0  , \quad \langle \partial_{z_2},\partial_{l_1}\rangle=1=-\langle \partial_{z_1},\partial_{l_2} \rangle \; 
$$ 
at $e\hat G_+$.

On $\RR^2\times\RR^2$ with coordinates $(v_1,v_2,u_1,u_2)$, we have the vector fields $\partial_{v_j}$ and $\partial_{u_j}$, $j=1,2$. Fix a point $(v,u)\in\RR^2\times\RR^2$. We may write $(v,u)=\Phi((v,0,u)\hat G_+)$. Then, the metric at $(v,u)$ is determined by the requirement that $D:=d(\Phi\circ l_{(v,0,u)})$ be an isometry from $T_{e\hat G_+}(\hat G/\hat G_+)$ to $T_{(v,u)}(\RR^2\times\RR^2)$. From \eqref{GrpMult}, compute that
\begin{align*}
    D\partial_{z_j} &= \partial_{v_j} \\
    D\partial_{l_j} &= -(-1)^j\textstyle\frac{\kappa}{6}u_{j+1}(u_1\partial_{v_1}+u_2\partial_{v_2}) + \partial_{u_j} \; .
\end{align*}
The subscript in $u_{j+1}$ is to be taken modulo 2.

Combining the formulae for the inner product on $T_{e\hat{G}_+}(\hat G/\hat{G}_+)$ and the differential $D=d(\Phi\circ l_{(v,0,u)})$, the metric $g$ on $\RR^2\times\RR^2$ can be obtained. We have
$$
0=\langle \partial_{z_j},\partial_{z_k}\rangle =g(D\partial_{z_j},D\partial_{z_k}) =g(\partial_{v_j},\partial_{v_k}) \; ,
$$
thus
\[
    1=\langle \partial_{z_j},\partial_{l_k}\rangle = g(D\partial_{z_j},D\partial_{l_k})=g(\partial_{v_j},\partial_{u_k}) \, .
\]
Furthermore,
\[
    0= \langle \partial_{l_1},\partial_{l_1}\rangle = g(D\partial_{l_1},D\partial_{l_1}) =\textstyle\frac{\kappa}{3}u_2^2 g(\partial_{v_2},\partial_{u_1}) +g(\partial_{u_1},\partial_{u_1}) \,. \]  
Hence
$$
g(\partial_{u_1},\partial_{u_1}) = - \textstyle\frac{\kappa}{3}u_2^2 g(\partial_{v_2},\partial_{u_1}) = -\textstyle\frac{\kappa}{3}u_2^2 \; . 
$$
The metric coefficients $g(\partial_{u_1},\partial_{u_2})$ and $g(\partial_{u_2},\partial_{u_2})$ are computed analogously. This leads to the claimed expression of the metric.

To compute the action of the isometry group, consider $(v,u)=\Phi((v,0, u)\hat G_+)$ and take $S\in\SL^\pm(2,\RR)$. Then 
\[\Phi\big(S\cdot (v,0,u)\hat G_+\big)=\Phi\big((|S|Sv, 0, Su)\hat G_+\big)=(|S|Sv,Su)\]
and similarly
\begin{eqnarray*}\textstyle
    \Phi\big((z,a,l)\cdot (v,0, u)\hat G_+\big) &=& \Phi \left( {\textstyle v+z+ \frac \kappa 3 \alpha(l,u)(l+ \frac 12  u), a+\frac12\alpha(l,u), u+l}\right) \\ &=& \left( {\textstyle v+z- \frac \kappa 3 \alpha(l,u)(u+ \frac 12  l) -\kappa a(u+l), u+l}\right)
\end{eqnarray*}
by (\ref{GrpMult}).
\qed
\begin{re} \label{Nherm}
{\rm We know from Proposition \ref{PHerm} that $N(\kappa)$ is a pseudo-Hermitian symmetric space. We can find complex coordinates $\xi_1, \xi_2$ on $N(\kappa)=\CC^2$ such that 
\[\textstyle g=d\xi_1d\bar \xi_2+d\xi_2 d\bar \xi_1-\frac\kappa 2 |\xi_2|^2d\xi_2d\bar \xi_2.\]
A suitable coordinate transformation is given by $\xi_1=x_1+iy_1$, $\xi_2=u_1+iu_2$, where
$$\textstyle x_1=v_2+\frac \kappa{12} u_1 |u|^2,\ y_1=-v_1+\frac \kappa{12} u_2 |u|^2. $$
}\end{re}
Besides a Hermitian structure, $N(\kappa)$ also admits a para-Hermitian one, see Proposition~\ref{Ppara}. Each of these structures leads to a representation of $N(\kappa)$ as an extrinsic symmetric space. In the following, we will determine these extrinsic symmetric spaces, as well as their extrinsic transvection and extrinsic isometry groups.

\begin{co}\label{coNext}
The space $N(\kappa)$ is isometric to the  extrinsic symmetric space
\begin{equation}\label{Next}
M:=\{(\eta_1,-\textstyle\frac12 |\eta_2|^2,\eta_2)\mid \eta_1,\eta_2\in \CC\}\subset \CC\oplus\RR\oplus \CC,
\end{equation}
where the scalar product $\ip$ on $\CC\oplus\RR\oplus \CC$ with coordinates $\eta_1,x,\eta_2$ is given by 
\begin{equation}\label{gh}
d\eta_1d\bar\eta_2+d\bar\eta_1d\eta_2+\kappa dx^2.
\end{equation}
An isometry from the pseudo-Hermitian symmetric space $N(\kappa)=\CC^2$ onto this extrinsic symmetric space is obtained by 
\begin{eqnarray*}
    \iota:\ N(\kappa)=\CC^2 &\longrightarrow &M\subset \CC\oplus\RR\oplus \CC,\\
    (\xi_1,\xi_2)&\longmapsto &(\xi_1-\textstyle\frac\kappa4 |\xi_2|^2\xi_2, -\frac12|\xi_2|^2,\xi_2).
\end{eqnarray*}

If we identify $\CC\cong \RR^2$, then the extrinsic transvection group of this space becomes the subgroup
\[ \left\{\left.
\left(\begin{array}{ccc|c} I&\kappa X &-\frac \kappa 2 X X^\top +Y& Z\\
0&1 &- X^\top & -\frac12 X^\top X \\ 0&0&I&X\end{array} \right)\ \right| X,Z\in\RR^2,\ Y=\begin{pmatrix}0&-y\\y&0\end{pmatrix},\ y\in \RR\right\}\]
of the group of affine isometries of $(\CC\oplus\RR\oplus \CC,\ip)$ and is isomorphic to the (ordinary) transvection group $\hat G$ of $N(\kappa)$. The group of all extrinsic isometries is equal to
$\grO(2)\ltimes \hat G$,
where $\grO(2)$ is understood as $\{\diag(S,1,S)\mid S\in \grO(2)\}\subset \SO(\CC\oplus\RR\oplus \CC,\ip)$. 
\end{co}
\proof 
We already observed that $N(\kappa)$ is a pseudo-Hermitian symmetric space. Hence, up to a covering map, it can be embedded into $V:=\hat \fg$ as an extrinsic symmetric space, see Proposition  \ref{APext}. In the following we will apply the general procedure explained there to the special case of $N(\kappa)$. We consider $\hat\fg$ in the form described in Proposition~\ref{class}, i.e., $\hat\fg=\fl^*\oplus\fa\oplus\fl\cong \CC\oplus\RR\oplus\CC$ as a vector space endowed with the scalar product $(\ref{gh})$. The K\"ahler structure corresponds to the derivation $J=(-J_{\fl}^*)\oplus 0\oplus J_\fl$ of $\hat\fg$, where $J_\fl$ was introduced in Proposition~\ref{PHerm}. If we identify $\fl$ and $\fl^*$ with $\CC$, then $-J_\fl^*$ and $J_\fl$ are just the multiplication by $i$.
Let $K$ denote the connected subgroup of $\SO(V)\ltimes V$ with Lie algebra $\fk:=\{\phi(u):=(\ad(u),-Ju)\mid u\in\hat\fg\}$. Then $M:=K(0):=\{k(0)\mid k\in K\}$ is an extrinsic symmetric space with (extrinsic) transvection group $K$. According to Proposition \ref{APext}, up to a covering map, it is isometric to $N(\kappa)$ as a symmetric space. We will show that $M$ is simply-connected and therefore isometric to $N(\kappa)$.

Let us first compute the Lie algebra $\fk$. If we use the Lie bracket of $\hat \fg$ as given in Prop.~\ref{class}, we obtain that $\phi(z+a+l)$ equals
\[\left(\begin{array}{ccc|c} 0&\kappa X &Y& Z'\\
0&0 &- X^\top & 0 \\ 0&0&0&X\end{array} \right),\ X=\begin{pmatrix} l_2\\-l_1\end{pmatrix},\ Z'=\begin{pmatrix} z_2\\-z_1\end{pmatrix},\ Y=\begin{pmatrix}0&-\kappa a\\\kappa a&0\end{pmatrix}\ \]
%\[\phi(z+a+l)=\left(\begin{array}{ccccc|c}
%0&0&\kappa l_2&0&-\kappa a&z_2\\
%0&0&-\kappa l_1&\kappa a& 0&-z_1\\
%0&0&0&-l_2&l_1&0\\
%0&0&0&0&0&l_2\\
%0&0&0&0&0&-l_1
%\end{array}\right)\]
for $z=(z_1,z_2)\in\fl^*=\RR^2$, $a\in\fa=\RR$ and $l=(l_1,l_2)\in\fl=\RR^2$. To obtain $K$, we take the exponential of $\phi(z+a+l)$. This exponential has the claimed form for a suitable $Z\in\RR^2$.   

Having calculated $K$, it is obvious that the orbit $M=K(0)$ is the set given by (\ref{Next}). Obviously, $M$ is simply-connected thus isometric to $N(\kappa)$. It is a direct computation that the map $\iota$ is an isometry.

If $M$ is considered as an ordinary symmetric space, the Lie algebra of its (ordinary) transvection group is isomorphic to $\fk$, see Section~\ref{Sext}. Since, moreover, both $K$ and $\hat G$ are simply connected, they are isomorphic.

The group of all extrinsic isometries is contained in the normaliser of $K$ in $\grO(V)\ltimes V$. Since $K$ acts transitively on $M$, each extrinsic isometry can be written as a product $f\cdot f_0$, where $f$ is a transvection and $f_0$ is an extrinsic isometry that fixes the base point $0\in M$. In particular, $f_0$ is a linear map. Since it normalises $K$, $f_0$ preserves $W:=\CC\oplus 0\oplus 0\subset \CC\oplus \RR\oplus\CC$ and therefore also $W^\perp=\CC\oplus \RR\subset \CC\oplus\RR\oplus\CC$.
On the one-dimensional space $W^\perp/ W$ the isometry acts by $\pm1$. Hence $f_0$ has the block form
\[ f_0=\begin{pmatrix}
    T &*&*\\0&\pm1& X^\top\\0&0&S
\end{pmatrix}\]
with respect to $\CC\oplus\RR\oplus\CC$. Since $f_0$ is an isometry, we obtain $T=(S^\top)^{-1}$. Now we use that $f_0$ preserves $M$, which yields
\[\mp \textstyle\frac12 |\eta_2|^2+X^\top\eta_2=-\frac12 |S\eta_2|^2\]
for all $\eta_2\in \CC\cong\RR^2$. If we apply this identity to $\eta_2=X$, we obtain $X=0$. This implies that $S$ is an orthogonal map and that $f_0$ acts by $+1$ on $W^\perp /W$. Consequently,
\[ f_0=\begin{pmatrix}
    I &0&Y\\0&1& 0\\0&0&I
\end{pmatrix}\cdot \diag(S,1,S),\quad Y=\begin{pmatrix}0&-y\\y&0\end{pmatrix} \]
where we used that $f_0$ is an isometry in order to determine $Y$. Since the first matrix is in $K$, the assertion follows.
\qed

\begin{co} 
The space $N(\kappa)$ is isometric to the extrinsic symmetric space
\begin{equation}
M:=\{(p,q_1q_2,q)\mid p=(p_1,p_2),q=(q_1,q_2)\in \RR^2\}\subset \RR^2\oplus\RR\oplus \RR^2=\RR^5,
\end{equation}
where the scalar product $\ip$ on $\RR^5$ with coordinates $p_1,p_2,x,q_1,q_2$ is given by 
\begin{equation}
-2(dp_1dq_1+dp_2dq_2)-\kappa dx^2.
\end{equation}

The extrinsic transvection group of this space becomes the subgroup of the group of affine isometries of $(\RR^5,\ip)$ consisting of matrices
\[ 
\left(\begin{array}{ccc|c} I&\kappa X &-\frac \kappa 2 X X^\top +Y& Z\\
0&1 &- X^\top &  x_1x_2\\ 0&0&I& \bar X\end{array} \right), \] 
where
\[X,Z\in\RR^2,\ X=\begin{pmatrix} x_1\\x_2\end{pmatrix},\ \bar X=\begin{pmatrix} x_2\\x_1\end{pmatrix},\  Y=\begin{pmatrix}0&-y\\y&0\end{pmatrix},\ y\in \RR.\]
It is isomorphic to the (ordinary) transvection group $\hat G$ of $N(\kappa)$. The group of all extrinsic isometries is equal to
$\hat G\rtimes (\grO(1,1)\times\ZZ_2)$,
where $\grO(1,1)$ is understood as $\{\diag((S^{-1})^\top,1,S)\mid S\in \grO(1,1)\}\subset \SO(\RR^2\oplus\RR\oplus \RR^2,\ip)$  and $\ZZ_2$ is generated by $\diag(1,-1,-1,1,-1)\in\grO(\RR^5)$.
\end{co}

\proof
It was noted already that $N(\kappa)$ is para-Hermitian as well. This leads to another embeddeding of $N(\kappa)$ into $V=\hat\fg$ as an extrinsic symmetric space, see Proposition \ref{APext2}. The extrinsic symmetric space $M$ and its extrinsic transvection group can be determined analogously as in the proof of Corollary \ref{coNext}, where we now use Proposition~\ref{APext2}. We turn to the extrinsic isometry group. Again, any extrinsic isometry may be written as $f\cdot f_0$, where $f$ is a transvection, and $f_0$ is a linear isometry of the form
\[ f_0=\begin{pmatrix}
    (S^\top)^{-1} &*&*\\0&\pm1& X^\top\\0&0&S
\end{pmatrix}.\]
The isometry $f_0$ must preserve $M$, which yields the condition that
$$
\pm q_1q_2 + X^\top q = \tilde q_1 \tilde q_2 
$$
where $(\tilde q_1,\tilde q_2):= Sq$. The term $X^\top q$, which is linear in $q$, must therefore be equal to the term $\pm q_1q_2-\tilde q_1 \tilde q_2$, which is quadratic in $q$. This is only possible if both terms are zero. On the one hand, this implies $X^\top =0$. On the other hand, $S$ has to satisfy $b(Sq)=\pm b(q)$ for all $q$, where $b$ denotes the quadratic form $b(q):=q_1q_2$. Up to a factor and a change of basis, $b$ coincides with the standard Lorentz inner product on $\R^2$. In the case that $b(Sq)=b(q)$, $S$ must therefore lie in $\grO(1,1)$. If $b(Sq)=-b(q)$, then $S$ must be of the form $\diag(1,-1)S'$ for $S'\in \grO(1,1)$. The claim now follows.
\qed

We turn to the existence of compact quotients of $N(\kappa)$. It is again possible to relax the condition that a discrete subgroup of the isometry group acts  freely. Indeed, Selberg's lemma allows us to deduce freeness of the action (up to passing to a finite index subgroup) from its properness and cocompactness. 

\begin{pr}\label{cqN}  
There exists no discrete subgroup of the isometry group of $N(\kappa)$ that acts properly and cocompactly on $N(\kappa)$.
\end{pr}
In order to prove this, we consider the connected component $G$ of the isometry group of $N(\kappa)$.
We denote by $c_g:G\to G$ the conjugation by $g\in G$. Then
\begin{eqnarray}
%c_{\hat S}\big((z,a,l) S\big)&=&(\hat S a,s,\hat Sl)\,c_{\hat S}(S),\label{cN1}\\
c_{(0,0,\hat l\,)}\big((z,a,l) S\big)&=&\textstyle (*,\frac 12 \alpha(l+\hat l, l-S\hat l\,), l-S\hat l+\hat l\,)S,\label{cN2} \\
%c_{(0,\hat a,0)}\big((z,a,l) S\big)&=&(z-\kappa l\hat a, a,l)S,\label{cN3}\\
c_{(\hat z,0,0)}\big((z,a,l) S\big)&=&(z+\hat z-S\hat z,a,l)S,\label{cN4}
\end{eqnarray}
where the exact value of the first component $*$ in (\ref{cN2}) will not be needed later. %Differentiating (\ref{cN2}) and (\ref{cN3}), we obtain
%\begin{eqnarray}
%\Ad((0,\hat a,0))(z,a,l)&=&(z-\kappa \hat a l,a,l)\\
%\Ad((0,0,\hat l))(z,a,l)&=&\textstyle(z+\frac\kappa2\alpha(\hat l,l)\hat l+ \kappa a\hat l,a+\alpha(\hat l,l),l) \label{Ad2} 
%\end{eqnarray}

\begin{lm}\label{6d}
If a discrete subgroup $\Gamma\subset G$ acts freely on $X$, then it is conjugate to a subgroup of the simply-connected nilpotent subgroup $\hat G \rtimes Z \subset G$, where $Z\subset \SL(2, {\Bbb R})$ denotes the subgroup of unipotent upper triangular matrices.
\end{lm}
\proof Assume that there exists an element $\gamma_0=(z_0,a_0,l_0)S_0\in\Gamma$ such that $1$ is not an eigenvalue of $S_0$. After possibly conjugating $\Gamma$ according to (\ref{cN2}) and (\ref{cN4}), we may assume $z_0=l_0=0$. Then $(v,u)=0$ is a fixed point of $\gamma_0$, which is a contradiction. Hence $\hat\Gamma:=\pro_{\SL(2,{\Bbb R})}(\Gamma)$ consists of maps for which 1 is an eigenvalue (of multiplicity 2). Therefore also the Zariski closure $\hat Z$ of $\hat \Gamma$ contains only such maps. By the Lie-Kolchin Theorem, $\hat Z$ is contained in a subgroup that is conjugate to the group $Z$ of unipotent upper triangular matrices. Consequently, $\Gamma$ is conjugate to a subgroup of $\hat G\rtimes Z \subset G$. \qed

{\sl Proof of Proposition~\ref{cqN}.} Assume $\Gamma$ were a discrete subgroup of the isometry group acting properly and cocompactly. We may assume that $\Gamma$ is contained in the connected component $G$ of the isometry group. The group $G=\hat G \rtimes \SL(2,\RR)$ is linear, since it has a trivial center, hence its adjoint representation is faithful. By Lemma \ref{freeAction} we may therefore assume that $\Gamma$ acts freely. By Lemma \ref{6d}, $\Gamma$ is contained in a simply-connected nilpotent subgroup $N'$ containing $\hat G$. Let $U$ be the syndetic hull of $\Gamma$, which exists by Proposition~\ref{PS}, and let $\fu$ be its Lie algebra. The Lie group $U$ acts properly and cocompactly on $N(\kappa)$ by Proposition~\ref{PKo}. Hence $\fu$ is a vector space complement of $\hat \fn'_+$ in the Lie algebra $\fn'$ of $N'$, see Proposition~\ref{PM}. By Lemma \ref{6d}, we may assume that $\fn'=\hat \fg\rtimes \RR$, where $\hat \fg\cong \fl\oplus\fa\oplus\fl$ as in the proof of Proposition~\ref{PN} and $s\in\RR$ acts trivially on $\fa$ and by $s\cdot n$ for $n=${\small $\begin{pmatrix} 0&1\\ 0&0\end{pmatrix}$} on $\fl\cong \RR^2$. In order to distinguish the two $\fl$-factors, we will rename the first one as $\fz$. We denote the standard basis of $\fz\cong \RR^2$ by $z_1,z_2$ and the one of $\fl\cong\RR^2$ by $l_1,l_2$. Then $z_1,z_2,a=[l_1,l_2], l_1,l_2$ is a basis of $\hat \fg$. Since both $\fn_-'= \hat\fg_-$ and $\fu$ are vector space complements of $\hat\fn_+'=\hat\fg_+\oplus\RR$ in $\fn'$, there is a linear map $\psi:\hat\fg_-\rightarrow \hat\fg_+\oplus\RR$ characterised by $v+\psi(v)\in \fu$ for all $v\in\hat \fg_-$. Then
\begin{eqnarray*}
[v_1+\psi(v_1), v_2+\psi(v_2)]&=&[v_1,v_2]+[\psi(v_1),\psi(v_2)]+[v_1,\psi(v_2)]+[\psi(v_1),v_2]\\
&=&[v_1,v_2]+[v_1,\psi(v_2)]+[\psi(v_1),v_2] \ \in \ \fu
\end{eqnarray*}
since $\hat\fg_+\oplus\RR$ is abelian. This implies 
\begin{equation}\label{section}
\psi\big([v_1,\psi(v_2)]+[\psi(v_1),v_2]\big)=[v_1,v_2]\in\fa
\end{equation} 
for all $v_1,v_2\in\hat\fg_-$ since  $[v_1,v_2]\in\hat\fg_+$ and $[v_1,\psi(v_2)]+[\psi(v_1),v_2]\in\hat\fg_-$.  This is equivalent to 
\begin{eqnarray}
&& \psi\big([\psi(z_i),l_j]\big)=0,\quad i,j=1,2,\label{(1)}\\
&& \psi\big([l_1,\psi(l_2)]+[\psi(l_1),l_2]\big)=a \label{(2)}.
\end{eqnarray}
Note that $\hat\fg_+\oplus\RR=\fa\oplus\RR$. Assume first that $\psi(l_1)\in\fa$.
Then $$[l_1,\psi(l_2)]+[\psi(l_1),l_2]=:z\in\Span\{z_1,z_2\},$$ thus $\psi(z)=a$ by \eqref{(2)}. On the other hand, $\psi([a,l_j])=\psi([\psi(z),l_j])=0$, $j=1,2$, by \eqref{(1)}. This implies $\psi(z_1)=\psi(z_2)=0$, which contradicts $\psi(z)=a$. Therefore, $\psi(l_1)\not\in\fa$. Let $\psi(l_1)=ta+sn$, $s\not=0$. For $\psi(z_1)=t_1 a+s_1 n$, we obtain from \eqref{(1)}
$$\psi\big([\psi(z_1),l_1]\big) = \psi\big([t_1a+s_1n,l_1]\big)=\psi(-\kappa t_1 z_1)=-\kappa t_1(t_1a+s_1n)=0,$$
thus $t_1=0$. Furthermore, $\psi\big([\psi(z_1),l_2]\big)=\psi([s_1 n,l_2])=s_1\psi(l_1)=0$ by \eqref{(1)}, which gives $s_1=0$. 
Hence $\psi(z_1)=0$. We have $[l_1,\psi(l_2)]+[\psi(l_1),l_2]=z+sl_1$ for some $z\in\fz$. By  \eqref{(2)}, we obtain $\psi(z+sl_1)=a$. Therefore,
$$[z+sl_1,\psi(l_2)]+[\psi(z+sl_1),l_2]=[z+sl_1,\psi(l_2)]+[a,l_2]=:z' $$
is in $\fz$. Hence \eqref{section} implies $\psi(z')=[z+sl_1,l_2]=sa$. Since $\psi(z')\not=0$ and $\psi(z_1)=0$, the vectors $z_1$ and $z'$ are linearly independent. Therefore, $\psi(\fz)\subset\fa$. However, $\psi(z+sl_1)=a$ now gives $s\psi(l_1)=a-\psi(z)\in\fa$, which contradicts the condition $\psi(l_1)\not\in\fa$ obtained above. Hence, no $\Gamma$ with the assumed properties can exist. 
\qed

\subsection{The spaces  $Z(\eps,c)$}\label{subsecZepsc}
The spaces $Z(\eps,c)$ are pseudo-Hermitian symmetric spaces. In particular, they are universal covers of spaces that can be embedded as extrinsic symmetric spaces, see Proposition~\ref{APext}. Here we will determine these extrinsic symmetric spaces explicitly. It will turn out that they are simply-connected, thus diffeomorphic to $Z(\eps,c)$. 

Take $\fl\in\{\fsu(2),\fsl(2,\RR)\}$, $\fa=0$, and let the involution $\theta_\fl$ and the 3-form $\gamma$ be given as in item 4 in the list in Subsection~\ref{S31}. 
We denote the symmetric triple obtained by these data by 
$\fd_{c}=(\fd_c,\theta,\ip)$. Note that the underlying vector space, the involution $\theta$ and the inner product $\ip$ do not change with $c$. 

We define
$$\textstyle  \beta_\fl:=-\eps\frac12\kappa_\fl, $$
where $\kappa_\fl$ is the Killing form of $\fl$. 

\begin{re} \label{Von}{\rm The inner product $\beta_\fl$ is positive definite for $\fl=\fsu(2)$, and has signature $(1,2)$ if $\fl=\fsl(2,\RR)$. The basis $e_1$, $e_2$, $e_3$ is orthonormal with respect to $\beta_\fl$. For $\fl=\fsl(2,\RR)$, the vector $e_1$ is time-like and $e_2$ and $e_3$ are space-like. Using this, we will often identify the inner product space $(\fl,\beta_\fl)$ with the standard (pseudo-)Euclidean space $(\RR^3,\diag(\eps,1,1))$. The notation  `$\perp$' will always mean orthogonal with respect to $\beta_\fl$ or the standard inner product, depending on the context.}
\end{re}
We identify 
\begin{equation}\label{beta}
\fl\cong\fl^*,\ 
 l\longmapsto z_l:=\beta_\fl(l,\cdot).
 \end{equation}
 Then $Al\mapsto (A^*)^{-1}z_l$ for any automorphism $A$ of $\fl$ and $Dl\mapsto -D^*z_l$ for any derivation $D$ of $\fl$. 
\begin{lm} \label{dc} The symmetric triple $\fd_c$ 
%with Hermitian structure $J_c$ 
is isomorphic to $(\fl\rtimes\fl,\theta_\fl\oplus\theta_\fl,\ip_c')$, where the first summand $\fl$ of $\fl\rtimes\fl$ is considered as an abelian subalgebra on which the second summand $\fl$ acts by the adjoint representation and the inner product $\ip'_c$ is given by
\begin{equation}\label{ipc} \la (l_1,l_1'), (l_2,l_2')\ra'_c=\beta_\fl(l_1,l_2')+\beta_\fl(l
_2,l_1')- 2c\beta_\fl(l_1',l_2').
\end{equation}
\end{lm}
 \proof  The proof is carried out in two steps. First we will prove that the symmetric triple $\fd_c$ 
is isomorphic to $(\fd_0=\fl^*\rtimes \fl, \theta, \ip_c)$, 
where $\fl$ acts on $\fl^*$ by the coadjoint action and $\ip_c=  \ip-2c \beta_\fl.$
It is easy to check that $\gamma(e_1,e_2,e_3)= c\beta_\fl([e_1,e_2],e_3)$. Hence  
\[F_1:\ \textstyle \fd_c\longrightarrow \fd_0,\quad  z+l\longmapsto  z+c\beta_\fl(l,\cdot)+l\]
for $z\in\fl^*$ and $l\in\fl$ is an isomorphism of Lie algebras, which commutes with $\theta$. Moreover, 
\begin{eqnarray*}
\la F_1(z_1+l_1),F_1(z_2+l_2)\ra_c&=&\textstyle \la z_1+ c\beta_\fl(l_1,\cdot)+l_1,z_2+c\beta_\fl(l_2,\cdot)+l_2 \ra_c\\
&=&z_1(l_2)+z_2(l_1)\\
&=&\la z_1+l_1, z_2+l_2\ra \;  .
\end{eqnarray*}
In the second step we show that $(\fd_0=\fl^*\rtimes \fl, \theta, \ip_c)$ is isomorphic to $(\fl\rtimes\fl,\theta_\fl\oplus\theta_\fl,\ip_c')$ using (\ref{beta}).
Indeed,
\[ F_2:\ \fl\rtimes\fl \longrightarrow \fl^*\rtimes\fl,\quad (l,l')\longmapsto\textstyle z_l+l'\]
is an isomorphism of Lie algebras since the Killing form is invariant under the adjoint representation.
Contrary to our usual notation, we write the elements of $\fl\rtimes\fl$ as pairs here to distinguish the two $\fl$-components. The map $F_2$ is also an isometry since
\begin{eqnarray*}
\la F_2(l_1,l'_1),F_2(l_2,l'_2)\ra_c&=&\textstyle \la z_{l_1}+l'_1,z_{l_2}+l'_2 \ra_c \\
 &=& z_{l_1}(l_2')+z_{l_2}(l_1')-2c\beta_\fl(l_1',l_2')\\
&=&\beta_\fl(l_1,l_2')+\beta_\fl(l
_2,l_1')-2c\beta_\fl(l'_1,l'_2)\\
&=&\la (l_1,l_1'),(l_2,l_2')\ra_c'. 
\end{eqnarray*}
Moreover, $F_2\circ (\theta_\fl\oplus\theta_\fl) =\theta \circ F_2$. 
 \qed

We denote by $\grO^+(1,2)\subset \grO(1,2)$ the subgroup of elements that preserve the time orientation, and by $\grO^-(1,2)\subset\grO(1,2)$ the subset of elements that reverse it.

\begin{pr}
\begin{enumerate}\label{isoz}
\item The transvection group $\hat G$ of the symmetric space $Z(\eps,c)$ is isomorphic to $\RR^3\rtimes \SO(3)$ if $\eps =1$ and to $\RR^3\rtimes \SO_0(1,2)$ if $\eps =-1$. In both cases, the stabiliser $\hat G_+$ is isomorphic to $\RR \times \SO(2)$. 
\item The isometry group of $Z(\eps,c)$ is isomorphic to $\RR^3\rtimes \grO(3)$ if $\eps =1$ and to $\RR^3\rtimes \grO^+(1,2)$ if $\eps =-1$, where an element $A$ of the (pseudo-)orthogonal group acts on $\RR^3$ by $x\mapsto |A|Ax$. 
\end{enumerate}
\end{pr}
\proof 
ad 1. In Lemma \ref{dc} we proved that the symmetric triple $(\hat \fg, \hat \theta, \ip)$ is isomorphic to $(\fl\rtimes\fl,\theta_\fl\oplus\theta_\fl,\ip_c')$, and in this proof we will identify the former with the latter. We understand $\hat G$ as a Lie group with Lie algebra $\fl\rtimes\fl$. First we consider  the case $\eps=1$, i.e., $\fl\rtimes\fl=\fsu(2)\rtimes\fsu(2)$. The simply-connected Lie group associated with this Lie algebra equals $\fsu(2)\rtimes \SU(2)$, where $\SU(2)$ acts on $\fsu(2)$ by the adjoint representation. The connected subgroup with Lie algebra $\hat\fg_+$ equals $\fso(2)\times \SO(2)$, where $\SO(2)$ is considered as a subgroup of $\SU(2)$ via its natural embedding.  Its intersection with the centre of $\fsu(2)\rtimes \SU(2)$ consists of the identity $(0,{\rm I})$ and the element $(0,-{\rm I})$. Since $\SU(2)/\{\pm {\rm I}\}\cong \SO(3)$ and $\fsu(2)\cong \fso(3)$, the transvection group of $Z(1,c)$ is isomorphic to $\fso(3)\rtimes \SO(3)$, see Section~\ref{S2}. We identify $\fl\cong\RR^3$ using the basis $e_1,e_2,e_3$, which is orthonormal with respect to $\beta_\fl$. Then the adjoint representation of $\SO(3)$ becomes the natural representation of $\SO(3)$ on $\RR^3$ and the claim follows. Moreover, the Lie algebra of the stabiliser equals $\RR e_1\rtimes \RR e_1$. Its connected subgroup equals $\RR e_1\rtimes \SO(2)$, where $\SO(2)=\{ A\in\SO(3)\mid Ae_1=e_1\}$. Thus $\hat G_+\cong \RR e_1 \rtimes \SO(2)$.
The case $\eps=-1$ is treated in the same way. Now the simply-connected Lie group with Lie algebra $\fl\rtimes\fl$ equals $\fsl(2,\RR)\rtimes \widetilde{\SL}(2,\RR)$. The connected subgroup with Lie algebra $\hat\fg_+$ is given by $\RR e_1\rtimes\RR$, where $\RR\subset\widetilde{\SL}(2,\RR)$ is the preimage of $\SO(2)\subset \SL(2,\RR)$ under the covering map $\widetilde{\SL}(2,\RR)\to\SL(2,\RR)$. This subgroup contains the centre $Z$ of $\fsl(2,\RR)\rtimes \widetilde{\SL}(2,\RR)$.  We have $\widetilde{\SL}(2,\RR)/Z\cong \PSL(2,\RR)\cong \SO_0(1,2)$. Moreover, we again identify $\fl\cong \RR^3$ using the basis $e_1,e_2,e_3$. We obtain $\hat G\cong \RR^3\rtimes \SO_0(1,2)$ and $\hat G_+\cong \RR\cdot e_1\rtimes \SO(2)$, where $\SO(2)=\{ A\in\SO_0(1,2)\mid Ae_1=e_1\}$.

ad 2. Let us consider the case $\eps=1$. 
First we determine the automorphism group of $(\hat\fg,\theta,\ip)\cong(\fl\rtimes\fl, \theta_\fl \oplus\theta_\fl,\ip'_c)$. If $\ph\in \Aut(\fl\rtimes\fl, \theta_\fl \oplus\theta_\fl,\ip'_c)$, then necessarily
\begin{equation}
    \ph= \begin{pmatrix}
    (\bar \ph^{-1})^* & R \\ 0 & \bar \ph
\end{pmatrix} :\ \fl\rtimes\fl\longrightarrow \fl\rtimes \fl, \label{phi}
\end{equation}
where $^*$ denotes the adjoint with respect to the inner product $\beta_\fl$. The map $\bar\ph$ has to be an automorphism of $ (\fl,\theta_\fl)$.  The automorphism group of $\fl=\fsu(2)\cong\fso(3)$ equals $\SO(3)$, which acts by conjugation. The subgroup of elements that commute with $\theta_\fl$ equals $ \{\diag(|A|,A)\mid A\in\grO(2)\}\cong \grO(2)$ with respect to the basis $e_1,e_2,e_3$ of $\fl$. On the other hand, for $\bar\ph\in \grO(2)$, the map $\bar\ph\oplus\bar\ph: \fl\rtimes\fl\ni(l,l')\mapsto (\bar\ph l,\bar\ph l')$ is an automorphism of $(\fl\rtimes\fl, \theta_\fl \oplus\theta_\fl,\ip'_c)$. Thus it remains to determine those automorphisms $\ph$ of $(\fl\rtimes\fl, \theta_\fl \oplus\theta_\fl,\ip'_c)$ for which $\bar\ph=\id_\fl$. Then, by (\ref{phi}), $\ph (l,l')=(l+R(l'),l')$. Since $\ph$ is a Lie algebra automorphism, $R$ has to be a derivation of $\fl$. All derivations of $\fl$ are of the form $\ad(l)$ for some $l\in\fl$. Since, moreover, $\ph$ has to commute with $\theta_\fl\oplus\theta_\fl$, the vector $l$ has to be in $\fl_+$, thus  $l=re_1$ for some $r\in\RR$. Consequently, $\ph(l,l')=(l+r[e_1,l']_\fl,l')$.
In particular, $\ph$ can be considered as the adjoint action of an element of $\hat G_+$, namely of $\exp(re_1,0)\in\hat G_+$. Furthermore, $\{\bar\ph\oplus\bar\ph\mid\bar\ph=\diag(|A|,A) \mbox{ for } A\in \grO(2)\}\cap \Ad(\hat G_+)$ is equal to  $\{\bar\ph\oplus\bar\ph\mid\bar\ph=\diag(1,A) \mbox{ for } A\in \SO(2)\}$. 

In summary, if $\bar P_0$ denotes the subgroup of $\Aut(\fl\rtimes\fl, \theta_\fl \oplus\theta_\fl,\ip'_c)$ generated by $\bar\ph\oplus\bar\ph$ for $\bar\ph=D:=\diag(-1,-1,1)$, then every automorphism of $(\fl\rtimes\fl, \theta_\fl \oplus\theta_\fl,\ip'_c)$ is of the form $\varphi'\cdot \Ad(g_+)$ for some $\varphi'\in \bar P_0$ and $g_+\in\hat G_+$. Since in addition $\bar P_0\cap \Ad(\hat G_+)=\{e\}$,  we can apply Corollary \ref{co2.4} to $\bar P_0=\la D\oplus D\ra\cong \ZZ_2$. This yields that the isometry group of $Z(1,c)$ is isomorphic to  $\hat G\rtimes \ZZ_2$. For the description of the  $\ZZ_2$-action we use the isomorphism $\hat G\cong\RR^3\rtimes \SO(3)$. Then the generator $D\oplus D$ of $\ZZ_2$ maps $(b,A)\in\hat G$ to $(Db, DAD^{-1})$. Finally we use that $(\RR^3\rtimes \SO(3))\rtimes\ZZ_2$  is isomorphic to $\RR^3\rtimes \grO(3)$, where $\grO(3)$ acts as claimed in item 2 of the proposition. The isomorphism is given by 
\begin{eqnarray}
  (\RR^3\rtimes \SO(3))\rtimes\ZZ_2 &\longrightarrow &  \RR^3\rtimes \grO(3)\nonumber\\
 (b,A,\id)&\longmapsto& (b, A)\label{O3}\\
 (b,A, D \oplus D) &\longmapsto & (b,-AD).\nonumber
\end{eqnarray}

The proof for $\eps=-1$ follows the same lines. 
Here we obtain again that $\bar P_0$ is generated by $\bar \ph \oplus \bar \ph$, where $\bar \ph=\diag(-1,-1,1)$. Now we get an isomorphism $\hat G\rtimes  \ZZ_2 \to \RR^3\rtimes \grO^+(1,2)$ in the same way as in (\ref{O3}).
\qed
\begin{de} We define an antisymmetric bilinear map
$\times_\eps:\RR^3\times\RR^3\to \RR^3$ by 
$$e_1\times_\eps e_2=e_3,\ e_2\times_\eps e_3=\eps e_1,\ e_3\times_\eps e_1=e_2. $$
If $\eps=1$, this is the standard cross product, and we will simply write $\times$. If $\eps = -1$, we will write $\times'$ instead of $\times_{-1}$.
\end{de}
\begin{pr}\label{Zepsc} 
 The set
\[M=\{ (v,u)\in \RR^3\times\RR^3 \mid \eps u_1^2+ u_2^2+ u_3^2=\eps,\, \eps u_1v_1+u_2v_2+u_3v_3=0\} \]
is an extrinsic symmetric space in $$\big(\RR^3\times\RR^3, \ip= 2(\eps du_1dv_1+ du_2dv_2+ du_3 dv_3)- 2c(\eps du_1^2+du_2^2+du_3^2) \big).$$

For $\eps=1$, this extrisic symmetric space is isometric to $Z(1,c)$. For $\eps=-1$ it consists of two connected components, each of which is isometric to $Z(-1,c)$.

 The group of extrinsic isometries of $M$ is isomorphic to $\RR^3\rtimes \grO(3)$ if $\eps=1$, and to $\RR^3\rtimes \grO(1,2)$ if $\eps=-1$, where $\grO(3)$ and $\grO(1,2)$ act on $\RR^3$ by $b\mapsto |A|Ab$. An extrinsic isometry $(b,A)$ acts on $\RR^3\times \RR^3$ by
\begin{equation}
(b,A)(v,u)=(Av+ b\times_\eps Au ,Au). \label{wirkung}
\end{equation}
In particular, for $\eps=1$, the isometry group of $Z(1,c)$ is isomorphic to the extrinsic isometry group of $M$. For $\eps=-1$, it is only isomorphic to the subgroup $\RR^3\rtimes\grO^+(1,2)$, whereas the action of $\grO^-(1,2)$ interchanges the connected components of~$M$.
\end{pr}

\begin{re}{\rm
Proposition \ref{Zepsc} shows that $Z(\eps,c)$ is diffeomorphic to the tangent bundle of the standard sphere $S^2\subset\RR^3$ if $\eps=1$. If $\eps=-1$, it is diffeomorphic to the tangent bundle of $H^2\subset\RR^{1,2}$, where $H^2$ is one sheet of the two-sheeted hyperboloid. When endowed with the metric induced by the ambient Lorentzian inner product, $H^2$ is a model for the Riemannian hyperbolic plane. 

Alternatively, the extrinsic symmetric space $M$ can be identified with the extrinsic symmetric space 
\[\{ x\in \RR^6\mid \eps x_4^2+ x_5^2+x_6^2=\eps, \, \eps x_1x_4+x_2x_5+x_3x_6=-c \eps
\}\]
in $\big(\RR^6, \ip=2(\eps dx_1dx_4+dx_2dx_5+dx_3dx_6)\big)$. An isometry is given by $\RR^3\times\RR^3\to \RR^6$, $(v,u)\mapsto x:=(v-cu,u)$.
}\end{re}
{\sl Proof of Proposition \ref{Zepsc}.\ }
Consider $V=\fl\oplus \fl$ as a vector space with inner product $\ip_c'$ defined by (\ref{ipc}). Elements of the vector space $V$ are denoted by $(v,u)$. The symmetric space $Z(\eps,c)$ is pseudo-Hermitian, as was noted in Proposition \ref{PHerm}. The symmetric triple of $Z(\eps,c)$ is isomorphic to $(\fl\rtimes\fl,\theta_\fl\oplus\theta_\fl,\ip_c')$ by Lemma \ref{dc}, and the pseudo-Hermitian structure on this latter symmetric triple is given by $J_0=\ad((0,e_1))$. Thus there is an extrinsic symmetric space in $V$ associated with this symmetric triple. We will follow the general procedure described in Section~\ref{Sext} to determine this extrinsic symmetric space. We need to compute
%$J_c=\ad(e_1-\frac c2\kappa_\fl(e_1,\cdot))$. 
the connected subgroup $K'\subset \SO_0(V)$ with Lie algebra $\ad(\fl\rtimes \fl)\subset \fso(V)$.
To this end, recall that the embedding $\ad(0\rtimes\fl)\cong\fso(\fl,\beta_\fl)\subset\fso(\fl\oplus \fl)=\fso(V)$ is given by $A(v,u)=(Av,Au)$ for $A\in \fso(\fl,\beta_\fl)$ and $(v,u)\in \fl\oplus \fl=V$. 
Furthermore,
$\ad((b,0))(v,u)=([b,u]_\fl,0)$. We put $K'_0:=\{e^{\ad(b,0)}\mid b\in\fl\}$, which consists of maps 
$(v,u)\mapsto (v+[b,u]_\fl, u)$ for $b\in\fl$. Then 
$K'=K'_0\rtimes \SO_0(\fl,\beta_\fl)$. Using (\ref{ext'}), the desired extrinsic symmetric space is given by
\begin{eqnarray*}
K'\big((0,e_1)\big)&=&K_0'\big({\rm SO}_0(\fl,\beta_\fl) \big((0,e_1)\big)\big)\\
&=&\textstyle \left\{\begin{array}{ll} K'_0(\{0\}\times \{u\in\fl\mid \beta_\fl(u,u)=1 \}),& \mbox{if } \eps=1,\\
K'_0(\{0\}\times \{u\in\fl\mid \beta_\fl(u,u)=-1, \, u_1>0 \}),& \mbox{if } \eps=-1,
\end{array}
\right.\\
&=&\textstyle \left\{\begin{array}{ll}\{(v,u)\mid \beta_\fl(u,u)=1,\ \beta_\fl(v,u)=0\},&\mbox{if } \eps=1,\\
\{(v,u)\mid \beta_\fl(u,u)=-1,\ \beta_\fl(v,u)=0, \, u_1>0\}, & \mbox{if } \eps=-1.
\end{array}
\right.
\end{eqnarray*}
This space is simply-connected, hence isometric to $Z(\eps,c)$. Finally, we identify $\fl \cong \RR^3$ according to Remark~\ref{Von}. Then $(V,\ip_c')$ becomes isometric to $\big(\RR^3\times\RR^3, \ip=2(\eps du_1dv_1+ du_2dv_2+ du_3 dv_3)- 2c(\eps du_1^2+du_2^2+du_3^2) \big)$. Under this identification, $K'((0,e_1))$ coincides with the space $M$ for $\eps=1$, and for $\eps=-1$ it is one of the two connected components of $M$.

We already computed that the isometry group of $Z(\eps,c)$ is isomorphic to $\RR^3\rtimes\mathrm{O}(3)$ if $\eps=1$, and to $\RR^3\rtimes\mathrm{O}^+(1,2)$ if $\eps=-1$. So if $\eps=1$, the group of extrinsic isometries of $M$ is contained in $\RR^3\rtimes\mathrm{O}(3)$. Analogously, if $\eps=-1$, then the group of extrinsic isometries of $M$  preserving the connected components of $M$ is a subgroup of $\RR^3\rtimes\mathrm{O}^+(1,2)$. On the other hand, the action of $\RR^3\rtimes \grO(3)$ or $\RR^3\rtimes\grO(1,2)$ on $\RR^3\times\RR^3$ given by \eqref{wirkung} is isometric with respect to the inner product $\ip$, and maps the subspace $M$ to itself, hence is an extrinsic isometry. Obviously, the action of $\diag(-1,1,1)\in \grO^-(1,2)$ interchanges the connected components of $M$. This proves the claim about the extrinsic isometry group.
\qed

We turn to the existence question for compact quotients of the spaces $Z(\eps,c)$. As a first observation, we prove that in the case $\eps=1$, no element of the transvection group is fixed point free.
\begin{lm}\label{notfree}
Every element of the transvection group of $Z(1,c)$ has a fixed point on $Z(1,c)$. In particular, no (discrete) subgroup of the transvection group acts freely.
\end{lm}
\proof
Recall that an element $(b,A)$ of the transvection group $\RR^3\rtimes \SO(3)$ of $Z(1,c)$ acts on a point $(v,u)\in Z(1,c)$ via
$
(b,A)(v,u) = (Av+ b\times Au ,Au)
$, where `$\times $' is the usual cross product on $\RR^3$.
In order to find a fixed point on $Z(1,c)$, we have to solve two equations:
\begin{align}
    Au &= u \label{eps1FP1} \\
    Av+b\times u&=v  \label{eps1FP2}
\end{align}
for $u,v\in \RR^3$ with $\|u\|=1$ and $u\perp v$ with respect to the standard Euclidean product.

Assume first that $A=I$ is the identity matrix. If $b=0$, then $(b,A)$ is the identity in $\RR^3\rtimes \SO(3)$, thus it fixes every element of $Z(1,c)$. If $b\not=0$, then put $u:=b/\|b\|\in S^2$. In this case, e.g., $(0,u)$ is a solution of (\ref{eps1FP1}) and  (\ref{eps1FP2}). 

Now suppose that $A\not=I$. Equation \eqref{eps1FP1} always has a solution. We fix such a solution $u\in S^2$ with $Au=u$. Then $A$ maps $u^\perp$ to $u^\perp.$ Let $A'$ denote the restriction of $A$ to $u^\perp$. Equation (\ref{eps1FP2}) now
reduces to the equation 
\begin{equation} \label{eps1FP2'}
    (A'-I)v=-b\times u
\end{equation}    
in $u^\perp$. Since $A\in \SO(3)$ fixes $u$ but is not the identity, $A'$ is a non-trivial rotation in $u^\perp$. In particular, $1$ is not an eigenvalue of $A'$. Hence Equation \eqref{eps1FP2'} has a solution $v\in u^\perp$ for every $b$. Thus $(b,A)$ has a fixed point. 
\qed

The symmetric space $Z(\eps,c)$ has no compact quotients by discrete subgroups of its isometry group. In fact, as in the previous sections the condition that the action providing the compact quotient be free may be dropped.

\begin{pr}  
There exists no discrete subgroup of the isometry group of $Z(\eps,c)$ that acts properly and cocompactly on $Z(\eps,c)$.
\end{pr}

\proof 
Let $\Gamma$ be a discrete subgroup of the isometry group of $Z(\eps,c)$, and assume that $\Gamma$ acts properly and cocompactly. The isometry group contains an affine subgroup of finite index, namely $\RR^3\rtimes\SO(3)$ for $\eps=1$ and $\RR^3\rtimes \SO_0(1,2)$ for $\eps=-1$. We may assume that $\Gamma$ is contained in this affine subgroup. Since affine groups are linear, we may therefore assume that $\Gamma$ acts freely by Lemma \ref{freeAction}. If $\eps=1$, the assertion then follows from Lemma~\ref{notfree}. Indeed, each non-finite subgroup of the isometry group contains elements of the transvection group different from the identity, thus it does not act freely.

Now let $\Gamma$ be a discrete subgroup of the isometry group of $Z(-1,c)$ acting freely and properly on $Z(-1,c)$. As above we assume that $\Gamma$ is contained in $\RR^3\rtimes \SO_0(1,2)$. We define
$$\Gamma_0:=\Gamma\cap \RR^3,\quad \hat \Gamma =\{A\in {\rm SO}_0(1,2)\mid \exists\, b\in\RR^3: (b,A)\in \Gamma\}.$$ 
The set $\Gamma_0$ is a discrete subgroup of $\RR^3$, and $\hat\Gamma$ is a subgroup of $\SO_0(1,2)$, which need not be discrete. Since $\SO_0(1,2)\cong\PSL(2,\RR)$, results and terminology for elements and subgroups of $\PSL(2,\RR)$ can also be applied to elements and subgroups of of $\SO_0(1,2)$. We will do so in the following, see Appendix \ref{PSL} for a summary of the relevant definition and facts.

Step 1: $\hat \Gamma\subset\SO_0(1,2)$ does not contain elliptic elements.

Assume  that $A\not=I$ were an elliptic element in $\hat \Gamma$, and choose $b\in\RR^3$ such that $(b,A)\in\Gamma$. We will construct a fixed point of $(b,A)$ on $Z(-1,c)$, which will be a contradiction to the freeness of the action of $\Gamma$. Because $A$ is elliptic, it admits a fixed point $u\in H^2\subset\RR^{1,2}$. Then $A\in\SO_0(1,2)$ maps $u^\perp$ to $u^\perp.$ Let $A'$ denote the restriction of $A$ to $u^\perp$. Note that the restriction of the metric on $\RR^{1,2}$ to $u^\perp\subset \RR^{1,2}$ is positive definite. Since $A\in \SO(1,2)$ fixes $u$ but is not the identity, $A'$ is a non-trivial rotation in $u^\perp$. In particular, $1$ is not an eigenvalue of $A'$. Thus 
$(A'-I)v=-b\times' u $ has a solution $v$
in $u^\perp$. Hence $(b,A)(v,u)=(Av+b\times' Au,Au)=(v,u)$. 

Step 2: The elements of $\hat \Gamma\subset \SO_0(1,2)$ do not have a non-vanishing common fixed vector in~$\RR^3$.

Assume that $f\in\RR^3$, $f\not=0$, were fixed by all $A\in\hat \Gamma$.  Then $F: \, \Gamma\backslash Z(-1,c)\to \RR$, $\Gamma\cdot(v,u)\mapsto \beta_\fl(u,f)$ would be a well-defined continuous map, because $\beta_\fl(Au,f) = \beta_\fl(Au,Af) =  \beta_\fl(u,f)$ for $(v,u)\in Z(-1,c)$, $A\in \hat \Gamma$. The map $F$ has non-compact image, which is a contradiction to the compactness of $\Gamma\backslash Z(-1,c)$.  

Step 3: $\Gamma_0\subset \RR^3$ does not contain time-like elements. 

Let $(b,I)\in\Gamma_0$. If $b$ were time-like, then $u:=b/(-\la b,b\ra)^{1/2}$ would be in $H^2$. But since $b\times' u =0$, the point $(0,u)\in Z(-1,c)$ would then be a fixed point of $(b,I)$, which is a contradiction. 

Step 4: $\dim \Span_{\Bbb R}(\Gamma_0)< 2$ 

It is clear from Step 3 that $\dim \Span_{\Bbb R}(\Gamma_0)\not=3$. Suppose $V_0:=\Span_{\Bbb R} (\Gamma_0)$ were two-dimensional. Each $A\in\hat\Gamma$ satisfies $A(V_0)=V_0$ and preserves the lattice $\Gamma_0$ of $V_0$. If the indefinite scalar product on $\RR^3$ restricted to $V_0$ were positive definite, this would imply that $\hat \Gamma$ is finite. But then $\Gamma\backslash Z(-1,c)$ would not be compact. Hence the scalar product on $V_0$ is not positive definite. Also, $V_0$ cannot contain a time-like vector because otherwise $\Gamma_0$ would also contain a time-like vector, which contradicts Step 3. Thus there is a basis $b_1, b_2$ of $V_0$ such that $b_1\perp V_0$ and $b_2$ is space-like and belongs to $\Gamma_0$. Any element of $A\in\hat \Gamma$ maps $V_0$ to $V_0$ thus also $V_0^\perp$ to $V_0^\perp$. Hence $A(b_1)=a b_1$ for some $a\in\RR$.  We have $a>0$ because $b_1$ is light-like and $A\in\SO_0(1,2)$  preserves the time orientation and hence each half of the light cone.  Furthermore, $A(b_2)=\pm b_2+ \mu  b_1$ since $A$ is an orthogonal map. On the other hand, $A(\Gamma_0)=\Gamma_0$, hence the matrix of $A|_{V_0}$ with respect to a basis consisting of elements of $\Gamma_0$ is an integer matrix. Consequently, $\big|A|_{V_0}\big|=\pm1$. This implies $a=1$. Thus $Ab_1=b_1$ for all $A\in\hat\Gamma$. This is impossible by Step 2.

Step 5: $\Gamma_0=\{0\}$ 

Assume that $\Gamma_0=\ZZ f$. The elements of $\hat \Gamma$ restrict to automorphisms of $\Gamma_0$, hence map $f$ to $\pm f$. Those elements $(b,A)\in \Gamma$ that satisfy $Af=f$ form a subgroup of index at most 2 in $\Gamma$, whence we may assume that $f$ is a common fixed point of $\hat\Gamma$. This is impossible by Step 2.

Step 6: $\hat \Gamma$ is not discrete. 

If $\hat \Gamma$ were discrete, it would act properly on the hyperboloid $H^2$ by Proposition \ref{Proper}. By Step 1 no element of $\hat\Gamma$ has a fixed point on $H^2$, so the action on $H^2$ is also free. We will show that the map $\Gamma\backslash Z(-1,c) \to \hat\Gamma\backslash H^2$ induced by the projection $Z(-1,c)\to H^2$ is a locally trivial fibration with fibre $\RR^2$. Since the fibre is non-compact, this contradicts the assumed compactness of $\Gamma\backslash Z(-1,c)$. Take $\Gamma\cdot u$ in $\hat\Gamma\backslash H^2$. Since $\hat \Gamma$ acts freely and properly on $H^2$, we can choose an open neighbourhood $U$ of $u$ in $H^2$ such that $\{\hat\gamma \in\hat \Gamma\mid \hat\gamma U\cap U\not=\emptyset\}=\{I\}$.  Now we consider the preimage $\tilde U$ of $U$ under the projection $Z(-1,c)\to H^2$. Since $\Gamma_0$ is trivial, the set $\{\gamma\in\Gamma\mid \gamma \tilde U\cap \tilde U\not=\emptyset\}$
also contains only the identity. Consequently, the restriction of the projection $H^2\to \hat \Gamma\backslash H^2$ to $U$ as well as the restriction of the projection $Z(-1,c)\to \Gamma\backslash Z(-1,c)$ to $\tilde U$ are homeomorphisms. So we may consider $U$ and $\tilde U$ as subsets $U\subset \hat \Gamma\backslash H^2$, $\tilde U \subset \Gamma\backslash Z(-1,c)$. Then $\tilde U$ is the preimage of $U$ under $\Gamma\backslash Z(-1,c) \to  \hat \Gamma\backslash H^2$. Since $\tilde U\cong U\times\RR^2$, the claim follows.

Step 7: Up to conjugation, $\hat \Gamma$ is virtually contained in 
the parabolic subgroup $\bP \subset \SO_0(1,2)$. 
 
Recall the definition of $\bP$ from equation \eqref{parSubgroup} in Appendix \ref{PSL}. By Step 1, $\hat\Gamma$ does not contain elliptic elements, and it is not discrete by Step 6. Hence it is elementary \cite[Theorem 8.3.1]{B}, see also Proposition \ref{NoEll}. The claim follows now from the description of the elementary groups, see Proposition \ref{ClassElementary}. Indeed, this is immediate if $\hat\Gamma$ falls into cases (b) or (c) of Proposition \ref{ClassElementary}. If $\hat\Gamma$ falls into case (a), it is abelian and purely parabolic or purely hyperbolic. In either case, one may pick a non-identity element $A\in\hat\Gamma$. This element is conjugate to an element of $\bP$, i.e. there exists $B\in\SO_0(1,2)$ such that $BAB^{-1}\in \bP$. Since $\hat\Gamma$ is abelian, any element of $B\hat\Gamma B^{-1}$ commutes with $BAB^{-1}$. Now, any element of $\SO_0(1,2)$ that commutes with a non-trivial element of $\bP$ is itself in $\bP$. It follows that $B\hat\Gamma B^{-1}\subset \bP$.

Step 8: $Z(-1,c)=(\RR^3\rtimes \bP)/\RR$.

It can be shown easily that $\bP$ acts simply transitively on the hyperboloid $H^2$. Moreover, $(b,I)(0,u)=( b\times' u ,u)$ for $b\in\RR^3$. Now the identity $u^\perp=\{b\times' u\mid b\in\RR^3\}$ proves that $\RR^3\rtimes \bP$ acts transitively on $Z(-1,c)$. In order to compute the stabiliser group we consider the base point $(0,e_1)\in Z(-1,c)$, where $e_1=(1,0,0)\in H^2\subset\RR^3$. Since $\bP$ acts freely on $H^2$, the stabiliser of this point is equal to $\{ (b,I)\mid(b\times'e_1,e_1)= (0,e_1)\}=\RR e_1\times\{I\}\subset \RR^3\rtimes \bP$.

Step 9: Final contradiction.

By the above steps, we may assume that $\Gamma \subset \RR^3 \rtimes \bP$. Instead of the standard basis $e_1,e_2,e_3$ of $\RR^3$ we will use here the basis $f_1, e_2, f_3$, where $f_1=(e_1+e_3)/\sqrt2$ and $f_3=(-e_1+e_3)/\sqrt{2}$. This is exactly the basis that we used to define $\bP$ in \eqref{parSubgroup}. Since $\RR^3 \rtimes \bP$ is completely solvable, $\Gamma$ has a unique syndetic hull $S\subset \RR^3 \rtimes \bP$, which acts properly and cocompactly on $Z(-1,c)=(\RR^3 \rtimes \bP)/\RR$, see Appendix \ref{Syndetic}. From Proposition \ref{PM} it follows that $\RR^3 \rtimes \fp = \RR \oplus \fs$ as vector spaces (where $\fp$ and $\fs$ are the Lie algebras of $\bP$ and $S$, respectively, and $\RR$ is spanned by $f_1-f_3 \in \RR^3$).
We denote by $q: \, \fs \to \fp$ the projection onto $\fp$. Note that this is a Lie algebra homomorphism. Because $\RR\oplus \fs = \RR^3 \rtimes \fp$, the map $q$ is surjective. For dimensional reasons its kernel is two-dimensional. Moreover, as the kernel of a Lie algebra homomorphism, $\ker(q)$ is an ideal in $\fs$. Since both $\ker(q)$ and the vector $f_1-f_3$ lie in the abelian subalgebra $\R^3$, it holds that $[f_1-f_3,\ker(q)]=0 $. We obtain that $\ker(q)$ is an ideal in $\RR^3\rtimes \fp = \RR\oplus\fs$. It follows that $\ker (q)$ is a two-dimensional $\bP$-invariant subspace of $\RR^3$, hence $\ker(q)=\Span\{f_1,e_2\}$. 
The Lie algebra $\fp$ has a basis $(X,Y)$ such that
$$
[X,f_1] = f_1 \; , \; [X,e_2]=0 \; , \; [Y,f_1]=0 \; , \; [Y,e_2] = f_1. 
$$
In particular, $[X,Y]=Y$. We have $\fs=\Span\{f_1,e_2,\tilde X,\tilde Y\}$, where $q(\tilde X)=X$ and $q(\tilde Y)=Y$. Here $\tilde X$ and $\tilde Y$ act on $\Span\{f_1,e_2\}$ in the same way as $X$ and $Y$. However, the adjoint action of $\tilde X$ on $\fs$ is not trace-free. This shows that $S$ is not unimodular, which contradicts the existence of a lattice in $S$.   We conclude that $\Gamma$ is not contained in $\RR^3\rtimes \bP$ giving a final contradiction. \qed

\subsection{The space $Z'(c)$}

Denote by $\fd'_c$ the symmetric triple obtained from the data given in item 5 in the list from subsection \ref{S31}. Observe that the Lie algebra and inner product of the symmetric triple $\fd_c'$ coincide with those of the symmetric triple $\fd_c$ from item 4 for the case $\eps=-1$, only the involution $\theta_\fl$ now has $e_2$ as fixed vector as opposed to $e_1$. Hence, we again introduce the inner product $\beta_\fl=\frac{1}{2} \kappa_\fl$, and identify $\fl^*\cong \fl$ via this inner product. Moreover, in the same way as in Lemma \ref{dc}, one can prove:

\begin{lm} \label{dc'} The symmetric triple $\fd_c'$ is isomorphic to $(\fl\rtimes\fl,\theta_\fl\oplus\theta_\fl,\ip_c')$, where the first summand $\fl$ of $\fl\rtimes\fl$ is considered as an abelian subalgebra on which the second summand $\fl$ acts by the adjoint representation and the inner product $\ip'_c$ is given by
\begin{equation}\la (l_1,l_1'), (l_2,l_2')\ra'_c=\beta_\fl(l_1,l_2')+\beta_\fl(l
_2,l_1')- 2c\beta_\fl(l_1',l_2').
\end{equation}
\end{lm}

The following description of the transvection group of $Z'(c)$ will make use of the universal cover $\widetilde{\SO}_0(1,2)$ of $\SO_0(1,2)$. Let $p:\widetilde{\SO}_0(1,2)\to \SO_0(1,2)$ be the covering map. We consider $\SO_0(1,1)$ as subgroup  of $\SO_0(1,2)$ that stabilises the vector $e_2$. Then  $p^{-1}(\SO_0(1,1))$ has a component group that is isomorphic to $\ZZ$. Let $(p^{-1}(\SO_0(1,1))_0\cong \SO_0(1,1)\cong \RR$ denote the identity component. Furthermore, let $\widetilde{\grO}(1,2)$ be the universal cover of $\grO(1,2)$ and denote the covering map by $p$ as well.
\begin{pr}\label{isoz'}
\begin{enumerate}
\item The transvection group $\hat G$ of the symmetric space $Z'(c)$ is isomorphic to  $\RR^3\rtimes \widetilde {\SO}_0(1,2)$. The stabiliser $\hat G_+$ is isomorphic to the subgroup $\RR e_2\times (p^{-1}(\SO_0(1,1))_0\ \cong \ \RR \times \RR$.
\item The isometry group of $Z'(c)$ is isomorphic to $\RR^3\rtimes \widetilde{\grO}(1,2)$, where an element $A$ of $\widetilde{\grO}(1,2)$ acts on $\RR^3$ by $x\mapsto |p(A)|\cdot p(A)(x)$.
\end{enumerate}
\end{pr}
\proof The Lie algebra of the transvection group of $Z'(c)$ can be identified with $\fl\rtimes\fl$. As in the case of $Z(-1,c)$ we have $\fl\rtimes\fl\cong \RR^3\rtimes \fso(1,2)$. Under this isomorphism, $\hat\fg_+$ corresponds to the Lie subalgebra $\RR e_2\times \fso(1,1)$. The simply-connected Lie group with Lie algebra $\RR^3\rtimes \fso(1,2)$ equals $\RR^3\rtimes \widetilde{\SO}_0(1,2)$. The connected subgroup with Lie algebra $\RR e_2\times \fso(1,1)$ is equal to $\RR e_2\times (p^{-1}(\SO_0(1,1))_0$. It does not contain elements of the centre of $\RR^3\rtimes \widetilde{\SO}_0(1,2)$ besides the identity. Thus the group $\RR^3\rtimes \widetilde{\SO}_0(1,2)$ acts effectively on $Z'(c)$, and is therefore isomorphic to the transvection group. This proves the first assertion.

In order to determine the isometry group, we proceed as in Proposition ~\ref{isoz}. 
We identify $(\hat\fg,\theta,\ip)$ with $(\fl\rtimes\fl, \theta_\fl \oplus\theta_\fl,\ip'_c)$. The automorphism group of $(\fl\rtimes\fl, \theta_\fl \oplus\theta_\fl,\ip'_c)$ consists of maps of the form 
\[
    \ph= \begin{pmatrix}
    \bar \ph& R \\ 0 & \bar \ph
\end{pmatrix} :\ \fl\rtimes\fl\longrightarrow \fl\rtimes \fl, 
\]
where $\bar \ph \in \{A\in \SO(1,2)\mid Ae_2=\pm1\}\cong \grO(1,1)$ and $R=\ad(te_2)$ for some $t\in\RR$. A map of this kind belongs to $\Ad(\hat G_+)$ if and only if $\bar \ph$ is in $\SO_0(1,1)\subset \grO(1,1)$. Now we define $\bar P_0$ to apply Corollary~\ref{co2.4}. For $D=\diag(\delta_1,\delta_2,\delta_3)\in\SO(1,2)$, let $\ph_D$ denote the automorphism $\bar\ph_D\oplus \bar\ph_D:\fl\rtimes\fl\to\fl\rtimes\fl$, where $\bar\ph_D$ is given by the matrix $D$ with respect to the basis $e_1,e_2,e_3$. Now we put \[\bar P_0=\{\ph_D\mid D=\diag(\delta_1,\delta_2,\delta_3)\in \SO(1,2)\}\cong \ZZ_2\times\ZZ_2.\] 
Then $\bar P_0$ satisfies the assumptions in Corollary~\ref{co2.4}. Each map $\ph_D$ defines an automorphism $F_D$ of $\RR^3\rtimes \SO_0(1,2)$ with differential $(F_D)_*=\ph_D$ by $F_D(b,A)=(DA,DAD^{-1})$. Let $\hat P_0$ denote the group of all these automorphisms $F_D$. Note that $F_D$ lifts to an automorphism $\tilde F_D$ of $\hat G=\RR^3\rtimes \widetilde\SO_0(1,2)$. These lifts constitute the group $\tilde P_0=\Psi_2^{-1}(\bar P_0)$. The map
\begin{eqnarray*}
 (\RR^3\rtimes{\rm SO}_0(1,2))\rtimes \hat P_0 &\longrightarrow &  \RR^3\rtimes \grO(1,2)\\
 (b,A,F_D)&\longmapsto& (b, \delta_2AD)
\end{eqnarray*}
for $D=\diag(\delta_1,\delta_2,\delta_3)$  is an isomorphism, where an element $A$ of ${\grO}(1,2)$ acts on $\RR^3$ by $x\mapsto |A|\cdot A(x)$. This implies that also the universal coverings $\Iso(Z'(c))=(\RR^3\rtimes \widetilde {\SO}_0(1,2))\rtimes \tilde P_0$ and $\RR^3\rtimes \widetilde{\grO}(1,2)$ are isomorphic. 
 \qed

\begin{pr}\label{Z'}
 The set
\[M=\{ (v,u)\in \RR^3\times\RR^3 \mid  -u_1^2+ u_2^2+ u_3^2=1,\, -u_1v_1+u_2v_2+u_3v_3=0\} \]
is an extrinsic symmetric space in 
\[\big(\RR^3\times\RR^3, \ip=- 2(- du_1dv_1+ du_2dv_2+ du_3 dv_3)+ 2c(-du_1^2+du_2^2+du_3^2) \big).\]
The universal cover of $M$ is isometric to $Z'(c)$. The group of extrinsic isometries of $M$ is isomorphic to $\RR^3\rtimes \grO(1,2)$, where $\grO(1,2)$ acts on $\RR^3$ by $b\mapsto |A|Ab$. An extrinsic isometry $(b,A)$ acts on $\RR^3\times \RR^3$ by
\[(b,A)(v,u)=(Av+ b\times' Au ,Au). \]
\end{pr}

\begin{re}{\rm 
Proposition \ref{Z'} shows that $Z'(c)$ is diffeomorphic to the tangent bundle of the universal cover $\tilde S^{1,1}$ of the pseudo-Riemannian sphere $S^{1,1}\subset \RR^{1,2}$. An element $(b,A)\in\RR^3\rtimes \widetilde\grO(1,2)$ of the isometry group of $Z'(c)$ acts on $Z'(c)=\{(v,u)\in \RR^3\times \tilde S^{1,1}\mid  v\perp \pi(u) \mbox{ in } \RR^{1,2}\} $ by
$$(b,A)(v,u)=(p(A)v+ b\times' \pi(Au) ,Au),$$
where $\pi: \, \tilde S^{1,1} \to S^{1,1}$ denotes the covering map. 

Alternatively, the extrinsic symmetric space $M$ can be identified with the extrinsic symmetric space 
\[\{ x\in \RR^6\mid -x_4^2+ x_5^2+x_6^2=1, \,  -x_1x_4+x_2x_5+x_3x_6=-c 
\}\]
in $\big(\RR^6, \ip=2(dx_1dx_4-dx_2dx_5-dx_3dx_6)\big)$. An isometry is given by $\RR^3\times\RR^3\to \RR^6$, $(v,u)\mapsto x:=(v-cu,u)$.
}\end{re}
{\sl Proof of Proposition \ref{Z'}.} 
The symmetric space $Z'(c)$ is para-Hermitian by Proposition \ref{Ppara}. Thus, the symmetric triple $(\fl\rtimes\fl,\theta_\fl\oplus\theta_\fl,\ip_c')$, being isomorphic to that of $Z'(c)$, admits a para-K\"ahler structure, namely $J=\ad\left((0,e_2)\right)$. According to Section~\ref{Sext}, there is an extrinsic symmetric space with $(\fl\rtimes\fl,\theta_\fl\oplus\theta_\fl,\ip_c')$ as its symmetric triple. To construct that extrinsic symmetric space, let $V$ be the vector space $\fl\oplus\fl$, endowed with the inner product $-\ip'_c$. We define $K'$ to be the closed connected subgroup of $\mathrm{SO}_0(V)$ whose Lie algebra is $\ad(\fl\rtimes\fl)\subset \fso(V)$. Precisely as in the proof of Proposition \ref{Zepsc}, one determines that $K'=K'_0 \rtimes \SO_0(\fl,\beta_\fl)$, where $K'_0$ consists of maps $V=\fl\oplus\fl \ni (v,u)\mapsto (v+[b,u]_\fl,u)  $ for $b\in\fl$, and $A\in\SO_0(\fl,\beta_\fl)$ acts as $(v,u)\mapsto(Av,Au)$. By \eqref{ext'}, the extrinsic symmetric space is given by 
\begin{align*}
    K'\big((0,e_2)\big) &= K'_0\left( {\SO}_0(\fl,\beta_\fl)\right)((0,e_2)) \\
    &= K'_0 \left( \{0\} \times \{ u\in\fl \, | \, \beta_\fl(u,u)=1 \} \right)\\
    &= \{ (v,u) \, | \, \beta_\fl(u,u)=1, \, \beta_\fl(u,v)=0 \} \; .
\end{align*}
Under the identification $\fl\cong\RR^3$, $\beta_\fl$ is identified with $\ip_{1,2}$  and $K'\big((0,e_2)\big)$ becomes
$$
M = \{ (v,u)\in\RR^3\times\RR^3 \, | \, \langle u,u\rangle_{1,2}=1 , \, \langle u,v\rangle_{1,2}=0 \} \; .
$$ 
This space is extrinsic symmetric in $\big(\RR^3\times\RR^3, \ip=- 2(- du_1dv_1+ du_2dv_2+ du_3 dv_3)+ 2c(-du_1^2+du_2^2+du_3^2) \big)$. Note that $M$ is \textit{not} simply connected, thus $Z'(c)$ is not given by $M$ itself, but its universal covering. The claim is proved.
\qed

We turn to the existence of compact quotients of $Z'(c)$. As for $Z(-1,c)$, our arguments use concepts and results regarding elements and subgroups of $\PSL(2,\RR)\cong\SO_0(1,2)$. These are recalled in Appendix \ref{PSL}.

Recall that $p: \, \widetilde{\SO}_0(1,2)\to \SO_0(1,2)$ and $\pi: \, \tilde S^{1,1} \to S^{1,1}$ denote the respective covering maps. Let $\mathrm{Hyp}\subset \SO_0(1,2)$ denote the set of hyperbolic elements, and $\widetilde{\mathrm{Hyp}} \subset \widetilde{\SO}_0(1,2)$ its preimage in the universal cover. Elements of $\widetilde{\mathrm{Hyp}}$ will also be called hyperbolic. The subset $\widetilde{\mathrm{Hyp}}$ is open, and has a countably infinite number of connected components. There is exactly one connected component, denote it by $\widetilde{\mathrm{Hyp}}_0$, which has the identity $I\in\widetilde{\SO}_0(1,2)$ as an accumulation point. A hyperbolic element $A\in\widetilde{\SO}_0(1,2)$ belongs to $A\in\widetilde{\mathrm{Hyp}}_0$ if and only if there is a one-parameter subgroup in $\widetilde\SO_0(1,2)$ containing $A$. For any hyperbolic $A\in\widetilde{\SO}_0(1,2)$ there exists a unique $k\in\ZZ$ such that $z^{-k}A\in\widetilde{\mathrm{Hyp}}_0$, in this case we write $A\in\widetilde{\mathrm{Hyp}}_k$. The subsets $\widetilde{\mathrm{Hyp}}_k$, $k\in \ZZ$, are precisely the connected components of $\widetilde{\mathrm{Hyp}}$ \cite{G}.
Every element $A$ of $\widetilde{\mathrm{Hyp}}_0$ has a fixed point on $\tilde S^{1,1}$: Since $p(A)$ is hyperbolic, it has a fixed point $x\in S^{1,1}$. In fact, $x$ is a common fixed point of a one-parameter subgroup through $p(A)$. That one-parameter subgroup lifts to a one-parameter subgroup $A(t)$ through $A$. Then, $A(t)$ preserves $\pi^{-1}\{x\}$. Because the preimage of $x$ is discrete, it follows for any $\tilde x \in \pi^{-1}\{x\}$ that $A\tilde x = A(0)\tilde x = \tilde x$.

In Appendix \ref{PSL} we introduce subgroups $\bN$, $\bH$, $\bP$ of $\SO_0(1,2)$, which we will use in the following. We denote by $\tilde\bN$, $\tilde \bH$, $\tilde \bP$ the preimages of $\bN$, $\bH$, $\bP$ in the universal cover $\widetilde\SO_0(1,2)$. We again use the standard basis $e_1,e_2,e_3$ of $\RR^3$ and the basis $f_1, e_2, f_3$, where $f_1=(e_1+e_3)/\sqrt2$ and $f_3=(-e_1+e_3)/\sqrt{2}$.

\begin{pr}\label{PZ'}
There exists no discrete subgroup of the isometry group of $Z'(c)$ that acts properly and cocompactly on $Z'(c)$.
\end{pr}

{\sl Proof of Prop.~\ref{PZ'}.} 
Let $\Gamma$ be a discrete subgroup of the isometry group of $Z'(c)$ acting properly and cocompactly on $Z'(c)$. In Proposition \ref{isoz'} it was shown that this isometry group is isomorphic to $\RR^3 \rtimes \widetilde \grO(1,2)$. We may assume that $\Gamma$ is contained in the finite index subgroup $\RR^3\rtimes \widetilde{\SO}_0(1,2)$. We may also assume that $\Gamma$ acts freely: Denote by $\Gamma'$ the projection of $\Gamma$ to $\RR^3\rtimes\SO_0(1,2)$. By Lemma \ref{finGen}, $\Gamma$ is finitely generated, hence so is $\Gamma'$. Since the affine group $\RR^3\rtimes \SO_0(1,2)$ is linear, we may apply Selberg's lemma to obtain a torsion-free finite index subgroup of $\Gamma'$. The preimage of this subgroup in $\Gamma$ is still torsion-free and has finite index. The proper action of a torsion-free group is free, see the proof of Lemma \ref{FreeFinIndex}.

Again, we define
$$
\Gamma_0 := \Gamma \cap \RR^3, \quad  \quad \hat \Gamma := \{ A \in \widetilde{\SO}_0(1,2) \, | \, \exists\, b\in\RR^3: \, (b,A)\in\Gamma \} \; . 
$$
Recall that $p: \, \widetilde{\SO}_0(1,2)\to \SO_0(1,2)$ and $\pi: \, \tilde S^{1,1} \to S^{1,1}$ denote the respective covering maps, and that $\pi(Au)=p(A)\pi(u)$ for all $A\in\widetilde{\SO}_0(1,2)$ and $u\in\tilde S^{1,1}$.

Step 1: There is no non-zero element $f\in\RR^3$ such that $p(\hat \Gamma)(f)\subset\{\pm f\}$, i.e., the elements of $p(\hat \Gamma)$ do not have a common non-zero fixed point up to sign in $\RR^3$.

If $f\in\RR^3$ satisfied $p(\hat \Gamma)(f)\subset\{\pm f\}$, then $\Gamma\setminus Z'(c) \to [0,\infty), \, \Gamma\cdot(v,u)\mapsto  |\beta_\fl(\pi(u),f) |$ would be a continuous surjection. Because $[0,\infty)$ is non-compact, the quotient $\Gamma\setminus Z'(c)$ could not be compact either, contradicting our assumption.

Step 2: If $A\in\hat \Gamma$ is hyperbolic, then $A$ has no fixed point on $\tilde S^{1,1}$. 

Let $A\in\hat\Gamma$, and $(b_A,A)$ be a corresponding element of $\Gamma$. Suppose $u\in\tilde S^{1,1}$ were a fixed point of $A$. Then $\pi(u)$ is an eigenvector of $p(A)$ with eigenvalue 1, and $p(A)$ preserves the two-dimensional subspace $\pi(u)^\perp$. Since the metric on $\pi(u)^\perp$ has signature $(1,1)$, $p(A)$ cannot have a further fixed vector on $\pi(u)^\perp$. Otherwise $p(A)$ would be the identity, which would contradict the hyperbolicity of $A$. Thus the restriction of $p(A)-I$ to $\pi(u)^\perp$ is bijective. In particular there exists a solution $v\in \pi(u)^\perp$ of $(p(A)-I)v=\pi(u)\times'b_A$, since the right side of this equation is an element of $\pi(u)^\perp$. The point $(v,u)\in Z'(c)$ would then be a fixed point of $(b_A,A)$ since $(b_A,A)(v,u)=(p(A)v+b_A\times'\pi(u),Au)=(v,u)$. We obtain a contradiction.

Step 3: The closure $C$ of $p(\hat \Gamma)$ in $\SO_0(1,2)$ is neither discrete nor equal to $\SO_0(1,2)$.

Here we view $C$ as a subgroup of $\PSL(2,\RR)$, and $\hat\Gamma$ as a subgroup of $\widetilde\PSL(2,\RR)$. There is a subset $\widetilde{\mathrm{Hyp}}_0 \subset\widetilde\PSL(2,\RR)$ consisting of those elements whose projection to $\PSL(2,\RR)$ is hyperbolic, and which belong to a one-parameter subgroup. The isomorphism $\PSL(2,\RR)\cong \SO_0(1,2)$ induces an isomorphism between the universal covers $\widetilde\PSL(2,\RR)$ and $\widetilde\SO_0(1,2)$, and this isomorphism identifies the subsets $\widetilde{\mathrm{Hyp}}_0$ in $\widetilde\PSL(2,\RR)$ and $\widetilde\SO_0(1,2)$. Thus, if there existed an element of $\hat\Gamma \subset \widetilde\PSL(2,\RR)$ that lies in $\widetilde{\mathrm{Hyp}}_0 \subset \widetilde\PSL(2,\RR)$, that would contradict Step 2, because elements of $\widetilde{\mathrm{Hyp}}_0 \subset \widetilde\SO_0(1,2)$ have fixed points on $\tilde S^{1,1}$. 

Assume first that $C$ is discrete. We apply Proposition \ref{discrSG}. Obviously, $C$ cannot be finite. Furthermore, we can also exclude the possibilities in items (b) and (c) of Proposition \ref{discrSG}. Indeed, up to conjugation $p(\hat\Gamma)$ would be contained in $\bN$ or $\bH$. For each of these groups, there exists a vector $f\in\RR^3$ that gets mapped to $\pm f$ by each group element. This is impossible by Step 1. Hence $C$ contains a subgroup which is isomorphic to a free group with two generators, and consists of hyperbolic elements only. Let $B_1,B_2$ be the two hyperbolic generators of this free subgroup. Let $\tilde B_1, \tilde B_2 \in \hat\Gamma\subset \widetilde\PSL(2,\RR)$ be preimages of $B_1,B_2 \in C$. The following arguments can be found in the paper \cite{Fa}, which includes a review of results by Goldman and by Matelski. The commutator $[B_1,B_2]$ is again hyperbolic. Thus, either $\tr [B_1,B_2]>2$, which implies $[\tilde B_1,\tilde B_2]\in\widetilde{\mathrm{Hyp}}_0$ \cite[Corollary 2.19]{Fa}, or otherwise $\tr [B_1,B_2]<-2$. In the latter case, $\tr [ B_1,[B_1,B_2]]>2$ \cite[Lemma 2.10 and Lemma 2.8]{Fa}. Thus $[\tilde B_1,[\tilde B_1,\tilde B_2]]\in\widetilde{\mathrm{Hyp}}_0$ again by \cite[Corollary 2.19]{Fa}. In either case, we get a contradiction to Step 2. Consequently, $C$ is not discrete. The group $C$ cannot be equal to $\PSL(2,\RR)$ since otherwise we would find two hyperbolic elements $A_1,A_2\in p(\hat \Gamma)$ with $\tr [A_1,A_2]>2$. If $\tilde A_1$ and $\tilde A_2$ are elements of $\hat \Gamma\subset \widetilde\PSL(2,\RR)$ with $p(\tilde A_j)=A_j$, $j=1,2$, then $[\tilde A_1,\tilde A_2]$ belongs to $\widetilde{\mathrm{Hyp}}_0$ \cite[Corollary 2.19]{Fa}. Again we get a contradiction to Step 2.

Step 4: $C$ is conjugate to a subgroup of $\bP$.

If $C$ is two-dimensional, it is conjugate to $\bP$. So it remains to consider the case where $C$ is one-dimensional. Then the identity component $C_0$ of $C$ is conjugate to one of the groups $\SO(2)$, $\bN$, $\bH_0$. If $C_0$ is conjugate to $\bN$, then $C$ is contained in $\bP$ and we are done. Otherwise $C$ is conjugate to one of the groups $\SO(2)$, $\bH_0$ or $\bH$. Again, each of these groups has a fixed point in $\RR^3$ up to a sign, which is impossible by Step~1.

Step 5: Final contradiction.

The result of Step 4 implies that $\Gamma$ may be assumed to be contained in $\fl^*\rtimes \tilde \bP$. We consider the submanifold $Y=(\fl^*\rtimes \tilde \bP)/(\RR\times\RR)$ of $Z'(c)$, where $\RR\times\RR$ is the stabilizer of a point in $Z'(c)$ as explained in Proposition \ref{isoz'}. The connected component $Y_0$ of $Y$ containing the base point is given by $(\fl^*\rtimes \tilde \bP_0)/(\RR\times\RR)$, where $\tilde \bP_0$ denotes the identity component of $\tilde \bP$. The subgroup $\Gamma_w:=\Gamma\cap (\fl^*\rtimes \tilde \bP_0)$ of $\Gamma$ acts properly and cocompactly on $Y_0$. Because $\fl^*\rtimes\tilde\bP_0$ is completely solvable, $\Gamma_w$ has a unique syndetic hull $S_0\subset \fl^*\rtimes\tilde\bP_0$, see Prop. \ref{PS}, which also acts properly and cocompactly on $Y_0$. Proposition \ref{PM} then implies that $\fl^* \rtimes \fp = \fs \oplus (\RR\cdot e_2\oplus\fh)$, where $\fs$ is the Lie algebra of $S_0$, and $\fh \subset \fp$ is the lie algebra of $\bH_0$.

Let $\hat p:\, \fl^*\rtimes\tilde \bP_0 \to \tilde\bP_0$ and $q: \fl^*\rtimes\fp \to \fp$ denote the projections of the groups and Lie algebras, respectively. Then $q(\fs)$ is a subalgebra of $\fp$. Since $\fl^* \rtimes \fp = \fs \oplus (\RR\cdot e_2\oplus\fh)$, the sum $q(\fs)+\fh$ must equal $\fp$. There are two cases: Either $\hat p(S_0)$ contains hyperbolic elements, or $q(\fs)=\fn$. In the case that $q(\fs)=\fn$, the kernel of $q$ restricted to $\fs$ is a two-dimensional $\fn$-invariant subspace of $\fl^*$, hence equal to $\Span\{f_1,e_2\}$. This is a contradiction, because $e_2$ is also contained in the direct summand $\RR\cdot e_2\oplus \fh$. Thus it is left to exclude the case that $\hat p(S_0)$ contains hyperbolic elements. Since $S_0$ is the syndetic hull of $\Gamma_w$, the quotient $\Gamma_w\backslash S_0$ is compact. Then $\hat p(\Gamma_w)\backslash \hat p(S_0)$ is also compact. It follows that if $\hat p(S_0)$ contains hyperbolic elements, then $\hat p(\Gamma_w)$ contains hyperbolic elements as well. Indeed, were $\hat p(\Gamma_w)\subset \tilde\bP_0$ to only contain elements of $\tilde\bN_0$, then the continuous map $\hat p( S_0) \to \RR, \, \tilde A\mapsto \langle f_3, p(\tilde A)f_1 \rangle$ would descend to a well-defined map $\hat p(\Gamma_w)\backslash \hat p(S_0) \to \RR$. If $\hat p(S_0)$ contains hyperbolic elements, these maps have unbounded image, as can be seen by examining the images of $\tilde A^n$ for $n\in\ZZ$ and $p(\tilde A)=P(a,b)$ with $a\neq 1$. This contradicts the compactness of $\hat p(\Gamma_w)\backslash \hat p(S_0)$. Now, any element in $\tilde\bP_0$ is part of a one-parameter family in $\tilde\bP_0$, in particular $\tilde\bP_0 \cap \widetilde{\mathrm{Hyp}}$ is contained in $\widetilde{\mathrm{Hyp}}_0$. Thus, if $\hat p(S_0)$ contains hyperbolic elements, then $\hat\Gamma$ contains elements of $\widetilde{\mathrm{Hyp}}_0$, contradicting Step 2. Overall, we conclude that $\Gamma$ cannot be contained in $\fl^*\rtimes\tilde\bP$, which gives the final contradiction.
\qed

\begin{appendix}  
\section{Proper actions and syndetic hulls}\label{Syndetic}

Recall that the action of a topological group $\Gamma$ on a space $X$ is proper if for every compact $K\subset X$, the set $\{ \gamma\in \Gamma \mid \gamma K\cap K\neq\emptyset \}$ is compact in $\Gamma$. In particular, if $\Gamma$ is a discrete group, then its action is proper if for every compact subset $K\subset X$ the set of $\gamma\in G$ such that $K \cap \gamma(K)\not=\emptyset$ is finite. 

In the context of this paper, $X$ is a smooth manifold. We consider a Lie group $G$ acting smoothly on $X$ by diffeomorphisms. If $\Gamma$ is a discrete subgroup of $G$ such that the induced action of $\Gamma$ is proper, free, and cocompact, then $\Gamma\backslash X$ is compact manifold and is called compact quotient of $X$. The assumption that the $\Gamma$-action be free is necessary for the quotient $\Gamma\backslash X$ to be a manifold. Often, however, this assumption is not essential if we are only interested in the existence or non-existence of compact quotients. In favorable circumstances, the freeness of the action is in fact a consequence of properness, at least up to passing to a finite index subgroup. This is the content of the next lemmas. 

Before we do so, let us note some simple facts that will be used implicitly throughout this paper. In the situation above, suppose that $\Gamma$ is a discrete subgroup of the Lie group $G$. Let $g$ be in $G$. Then $\Gamma$ acts properly, freely and cocompactly on $X$ if and only if $g\Gamma g^{-1}$ does. Moreover, suppose $\Gamma' \subset \Gamma$ is a finite-index subgroup. Then, $\Gamma'\backslash X$ is a finite cover of $\Gamma\backslash X$, thus $\Gamma'$ acts cocompactly if and only if $\Gamma$ does.

We turn to the aforementioned fact that the freeness of a group action is implied by properness in some cases, at least up to passing to a finite index subgroup. This is essentially due to Selbergs's lemma. Though well known, we give the arguments for convenience. 

We call a group $\Gamma$ linear if there exists an injective homomorphism $\Gamma\to\GL(n,\RR)$ for some $n\in\NN$.  
\begin{lm}\label{FreeFinIndex}
   Let $\Gamma$ be a linear finitely generated discrete group, which acts properly on a Hausdorff space $X$. Then, $\Gamma$ has a finite-index subgroup that acts freely. 
\end{lm}
\proof
Suppose $x\in X$ is a fixed point of $\gamma\in\Gamma$. Then, $x$ is also fixed by $\gamma^k$ for all $k\in\NN$. By properness of the action, the stabiliser $\Gamma_x\subset \Gamma$ of $x$ is finite.  Thus, $\gamma$ must be a torsion element. The claim now follows from Selberg's lemma, which states that every finitely generated linear group has a torsion-free subgroup of finite index.
\qed

\begin{lm}\label{finGen}
    Let $\Gamma$ be a discrete group acting properly on a locally compact, connected space $X$ by homeomorphisms, such that the quotient $\Gamma/X$ is compact. Then, $\Gamma$ is finitely generated.
\end{lm}

\proof
Let $A\subset X$ be such that $\Gamma\cdot A=X$ and $\Gamma_A\cdot A$ is a neighborhood of $A$, where $\Gamma_A=\{ \gamma\in \Gamma \mid \gamma\cdot A \cap A\neq 0\}$. It is proved in \cite{K} that $\Gamma_A$ then generates $\Gamma$. We will show that there exists a set $A$ of this kind such that $\Gamma_A$ is finite. Cover $X$ by open subsets $U_\alpha$ with compact closure. Denote by $\pi: \, X\to \Gamma\backslash X$ the quotient map. The sets $\pi(U_\alpha)$ form an open cover of $\Gamma\backslash X$, hence there exists a finite subcover $\pi(U_{\alpha_1}),\cdots,\pi(U_{\alpha_k})$. Let $A:= U_{\alpha_1}\cup\cdots \cup U_{\alpha_k}$. Because $\pi(U_{\alpha_1}),\cdots,\pi(U_{\alpha_k})$ is an open cover of $\Gamma\backslash X$, every point in $X$ lies in the orbit of some point in $A$, hence $\Gamma\cdot A= X$. Moreover, $\Gamma_A\cdot A$ is the union of the open sets $\gamma\cdot A$ for $\gamma\in\Gamma_A$, and because $A=e\cdot A \subset \Gamma_A\cdot A$, the set $\Gamma_A\cdot A$ is a neighborhood of $A$. Thus, $\Gamma_A$ generates $\Gamma$. To see that it is finite, note that $\Gamma_A$ is a subset of $\Gamma_{\bar{A}}$. The closure $\bar A$ is compact by assumption. Because $\Gamma$ is discrete and acts properly, $\Gamma_{\bar A}$ and thus $\Gamma_A$ is finite.   
\qed 

We summarise our discussion in the following lemma. Note that it applies to connected manifolds in particular.

\begin{lm}\label{freeAction}
    Let $\Gamma$ be a linear discrete group that acts properly and cocompactly on a connected, locally compact Hausdorff space $X$ by homeomorphisms. Then, $\Gamma$ contains a finite index subgroup that acts freely on $X$.
\end{lm}

Next, we recall the notion of a syndetic hull. It is one of the main technical tools in our study of the existence of compact quotients. Let $\Gamma$ be a closed subgroup of a Lie group $G$. A connected Lie subgroup $S$ of $G$ with $\Gamma \subset S$ and $\Gamma \backslash S$ compact is called a syndetic hull of $\Gamma$ in $G$. In general, syndetic hulls do not need to exist or be unique if they do. However, for completely solvable Lie groups existence and uniqueness is in fact guaranteed:

\begin{pr}{\rm \cite{S}} \label{PS}
    Every closed subgroup of a completely solvable Lie group possesses a unique syndetic hull.
\end{pr}

The properness of the action of a discrete subgroup is equivalent to the properness of the action of its syndetic hull, as implied by the following proposition.

\begin{pr}{\rm\cite{Ko1}} \label{PKo}
    Suppose a locally compact group $S$ acts on a locally compact Hausdorff space $X$. Let $\Gamma$ be a co-compact discrete subgroup of $S$. Then the $S$-action on $X$ is proper if and only if the $\Gamma$-action on $X$ is proper.
\end{pr}

Lastly, we state a condition for the compactness of a quotient of a $G$-homogeneous space for solvable $G$.

\begin{pr}{\rm\cite[Prop.~3.7]{M}}\label{PM}
Let $G$ be a simply-connected solvable Lie group and $L$ and $H$ be connected closed subgroups of $G$. Assume that the $L$-action on $G/H$ is proper. Then the following conditions are equivalent: 
\begin{enumerate}[label=(\roman*)]
\item the space $L\backslash G/H$ is compact,
\item $G=LH$,
\item $\fg=\fl\oplus\fh$ as a linear space.
\end{enumerate}
\end{pr}

\section[Subgroups of $\protect\mathrm{PSL}(2,\RR)$ and $\SO_0(1,2)$]{Subgroups of $\protect\mathrm{PSL}(2,\RRR)$ and $\SO_0(1,2)$}\label{PSL}

The group $\PSL(2,\RR)$, its elements and subgroups have been studied extensively in the context of hyperbolic geometry. The Lie groups $\SO_0(1,2)$ and $\PSL(2,\RR)$ are isomorphic, whence these results may be translated to results on $\SO_0(1,2)$. The transvection and isometry groups of some of the spaces examined in this paper are extensions of $\SO_0(1,2)$. Subgroups of $\SO_0(1,2)$ thus appear naturally in the study of compact quotients of these spaces. In this appendix we recall some relevant notions for elements and subgroups of $\PSL(2,\RR)$, and translate them to statements for $\SO_0(1,2)$ as they will be needed in the paper.

Recall that $\PSL(2,\RR)$ is the quotient of $\SL(2,\RR)$ by the subgroup $\{ \pm I \}$. It acts as the connected component of the isometry group of the hyperbolic plane $\HH$ via Möbius transformations. Consider the Lie algebra $\fsl(2,\RR) \cong \fso(1,2)$. Its Killing form has signature $(1,2)$, and is preserved by the adjoint action of $\PSL(2,\RR)$. The mapping $A \mapsto \Ad(A)$ establishes an isomorphism from $\PSL(2,\RR)$ to the connected component $\SO_0(1,2)$ of the Lorentz group. 

Elements of $\PSL(2,\RR)$ are categorised according to their trace. Namely, $A \in \PSL(2,\RR)$ is called elliptic if $|\tr(A)|<2$, parabolic if $|\tr(A)|=2$, and hyperbolic if $|\tr(A)|>2$. Moreover, the absolute value of the trace is closely tied to conjugacy classes: If two elements of $\PSL(2,\RR)$ are conjugate, then the absolute value of their trace is equal. If $A\in\PSL(2,\RR)$ is hyperbolic, then the conjugacy class of $A$ is determined by $|\tr(A)|$, i.e. $|\tr(B)|=|\tr(A)|$ implies that $B$ and $A$ are conjugate. There are three conjugacy classes of parabolic elements, one of them containing only the identity. If $A\in\PSL(2,\RR)$ is elliptic, then the set of $B\in\PSL(2,\RR)$ with $|\tr(B)|=|\tr(A)|$ consists of two conjugacy classes.

We will call an element of $\SO_0(1,2)$ elliptic/parabolic/hyperbolic if its image under the isomorphism $\SO_0(1,2)\cong \PSL(2,\RR)$ is elliptic/parabolic/hyperbolic.
Elliptic elements correspond to rotations, hyperbolic ones to hyperbolic rotations. Elliptic elements have a fixed point in the hyperbolic plane $H^2\subset \RR^{1,2}$, hyperbolic elements have a fixed point on the hyperboloid $S^{1,1}\subset \RR^{1,2}$. Parabolic elements fix a light-like vector.

The action of $\PSL(2,\RR)$ on $\HH$ can be extended to an action on $\widehat{\HH}=\HH \cup \partial \HH \cup \{\infty\}$. A subgroup $\Gamma$ of $\mathrm{PSL}(2,\RR)$ is called elementary if the action of $\Gamma$ on $\widehat{\HH}$ has a finite orbit. We call a subgroup $\Gamma \subset \SO_0(1,2)$ elementary if its image under the isomorphism $\SO_0(1,2)\cong \PSL(2,\RR)$ is elementary. Elementary subgroups are classified, indeed any elementary subgroup of $\PSL(2,\RR)$ is of one of the following forms:
\begin{itemize}
    \item[(a)] abelian, and containing only elliptic, or only parabolic, or only hyperbolic elements,
    \item[(b)] conjugate in $\mathrm{PSL}(2,\RR)$ to a subgroup of the image of the subgroup 
    $$
    \left\{ \begin{pmatrix}
        a & b \\ 0 & a^{-1}
    \end{pmatrix} \mid a\in\RR\setminus \{0\} , \, b\in\RR \right\} \; \subset \SL(2,\RR)\; 
    $$
    in $\PSL(2,\RR)$,
    \item[(c)] or conjugate in $\mathrm{PSL}(2,\RR)$ to a subgroup of the image of the subgroup
    $$
    \left\{ \begin{pmatrix}
        a & 0 \\ 0 & a^{-1}
    \end{pmatrix} \mid a\in\RR\setminus \{0\} \right\} \cup \left\{ \begin{pmatrix}
        0 & a \\ -a^{-1} & 0
    \end{pmatrix} \mid a\in\RR\setminus \{0\} \right\}\; \subset \SL(2,\RR) 
    $$
    in $\PSL(2,\RR)$.
\end{itemize}
A proof can be found in {\rm\cite[Section~5.1]{B}}. The corresponding statement in $\SO_0(1,2)$ is as follows:

\begin{pr} \label{ClassElementary}
Any elementary subgroup of $\SO_0(1,2)$ is of one of the following forms:
\begin{itemize}
    \item[(a)] abelian, and containing only elliptic, or only parabolic, or only hyperbolic elements,
    \item[(b)] conjugate in $\SO_0(1,2)$ to a subgroup of
    \begin{equation}
\bP:=\left\{ P(a,b):=\begin{pmatrix}
    a^2 & -\sqrt{2}ab & -b^2 \\ 0 &1& \sqrt{2}ba^{-1} \\ 0&0&a^{-2}
\end{pmatrix} \mid a\neq 0 , \, b\in\RR \right\} \; \subset \mathrm{SO}(1,2)_0 \; , \label{parSubgroup}
\end{equation}
    \item[(c)] or conjugate in $\SO_0(1,2)$ to a subgroup of 
    $$\bH:=\bH_0\ \cup\  \bH_0 \cdot \left(\begin{array}{ccc} 0&0&-1\\ 0&-1&0\\-1&0&0\end{array}\right)$$ for $\bH_0:=\{P(a,0)=\diag (a^2,1,a^{-2})\mid a\neq 0\}$.
\end{itemize}
\end{pr}
Here, the matrices are represented with respect to a basis $(f_1,e_2,f_3)$ of $\RR^{1,2}$, where $f_1,f_3$ are isotropic and $\la f_1,f_3 \ra_{1,2}=1$, and $e_2$ is spacelike and orthonormal to $f_1,f_3$. Indeed, if $(e_1,e_2,e_3)$ denotes the standard orthonormal basis, we let $f_1=(e_1+e_3)/\sqrt2$ and $f_3=(-e_1+e_3)/\sqrt{2}$. We also need the subgroup
$$
\bN:=\left\{ P(1,b) \mid  b\in\RR \right\} \; \subset \mathrm{SO}(1,2)_0 \; .
$$
The elements of $\bN$ are parabolic, and any parabolic element of $\SO_0(1,2)$ is conjugate to an element of $\bN$. The elements of $\bP$ are either hyperbolic or in $\bN$. Any hyperbolic element of $\SO_0(1,2)$ is conjugate to an element of $\bP$, in fact even conjugate to an element of the subgroup $\bH_0\subset \bP$.

The following proposition describes subgroups of $\SO_0(1,2)$, which contain no elliptic elements.

\begin{pr}{\rm\cite[Thm.~8.3.1]{B}}\label{NoEll}
If a subgroup $\Gamma$ of $\SO_0(1,2)$ contains no elliptic elements, then it is discrete or elementary.
\end{pr}

Any discrete subgroup of $\mathrm{PSL}(2,\RR)$ acts properly on the hyperbolic plane $\HH$ \cite[Thm.~8.4.1]{B}. Under the isomorphism $\PSL(2,\RR)\cong \SO_0(1,2)$, the action of $\PSL(2,\RR)$ on $\HH$ is translated to the natural action of $\SO_0(1,2)$ on the hyperboloid $H^2\subset \RR^{1,2}$. Thus, the following holds:

\begin{pr}\label{Proper}
Any discrete subgroup of $\mathrm{SO}_0(1,2)$ acts properly on the hyperboloid $H^2\subset \RR^{1,2}$.
\end{pr}

Moreover, the following can be said about discrete subgroups.

\begin{pr}{\rm\cite[Proposition 3.1.2]{Hub}}\label{discrSG}
    For any discrete subgroup $\Gamma\subset\SO_0(1,2)$, one of the following holds:
    \begin{itemize}
        \item[(a)] $\Gamma$ is finite and consists only of elliptic elements.
        \item[(b)] $\Gamma$ is infinite cyclic generated by a parabolic element.
        \item[(c)] $\Gamma$ is infinite cyclic generated by a hyperbolic element, or $\Gamma$ contains an infinite cyclic subgroup of index 2 generated by a hyperbolic element. 
        \item[(d)] $\Gamma$ contains a subgroup consisting entirely of hyperbolic elements, which is isomorphic to the free group on two generators. 
    \end{itemize}
    In the cases (a), (b) and (c), $\Gamma$ is elementary.
\end{pr}

Note that in case (b) of the above proposition, $\Gamma$ is conjugate to a subgroup of $\bN$. In case (c), it is conjugate to a subgroup of $\bH$. \\[2ex]

\end{appendix}

{\bf Conflict of interest statement.} The authors have no conflicts of interest to declare that are relevant to this article.

\medskip

\noindent {\bf Data availability statement.}
    Data sharing not applicable to this article as no datasets were generated or analyzed during the current study.

%\vspace{1.6cm}
%{\footnotesize 
%
%Ines Kath (corresponding author)\\ 
%Institut f\"ur Mathematik und Informatik\\ 
%der Universit\"at Greifswald\\
%Walther-Rathenau-Str.\,47, D-17489 Greifswald, Germany\\ 
%email: ines.kath@uni-greifswald.de\\[2ex] 
%}
\end{document}